\documentstyle[myart,amsfonts,12pt]{article}
\normalbaselines
\font\tenbm=cmmib10
\font\sevenbm=cmmib7
\font\fivebm=cmmib5
\newfam\bmfam
\textfont\bmfam=\tenbm \scriptfont\bmfam=\sevenbm
\scriptscriptfont\bmfam=\fivebm
{\count0=\number\bmfam \multiply\count0 by "100
\def\defbgreek#1#2#3{{\count1=\count0 \advance\count1 by "#2#3
  \global\mathchardef#1=\count1 }}
\defbgreek\balpha  0B \defbgreek\brho       1A
\defbgreek\bbeta   0C \defbgreek\bsigma     1B
\defbgreek\bgamma  0D \defbgreek\btau       1C
\defbgreek\bdelta  0E \defbgreek\bupsilon   1D
\defbgreek\bepsilon0F \defbgreek\bphi       1E
\defbgreek\bzeta   10 \defbgreek\bchi       1F
\defbgreek\bmeta   11 \defbgreek\bpsi       20
\defbgreek\btheta  12 \defbgreek\bomega     21
\defbgreek\biota   13 \defbgreek\bvarepsilon22
\defbgreek\bkappa  14 \defbgreek\bvartheta  23
\defbgreek\blambda 15 \defbgreek\bvarpi     24
\defbgreek\bmu     16 \defbgreek\bvarrho    25
\defbgreek\bnu     17 \defbgreek\bvarsigma  26
        \defbgreek\bxi     18 \defbgreek\bvarphi    27
\defbgreek\bpi     19}
\newtheorem{opred}{Definition}[section]
\input tcilatex.tex

\begin{document}

\title{Asymmetric Nondegenerate Geometry }
\author{Yuri A. Rylov}
\date{Institute for Problems in Mechanics, Russian Academy of Sciences, \\
101-1, Vernadskii Ave., Moscow, 119526, Russia. \\
email: rylov@ipmnet.ru\\
Web site: {$http://rsfq1.physics.sunysb.edu/\symbol{126}rylov/yrylov.htm$}.\\
or mirror Web site: {$http://194.190.131.172/\symbol{126}rylov/yrylov.htm$}.}
\maketitle

\begin{abstract}
Nondegenerate geometry (T-geometry) with nonsymmetric world function is
considered. In application to the space-time geometry the asymmetry of world
function means that the past and the future are not equivalent
geometrically. T-geometry is described in terms of finite point subspaces
and world function between pairs of points of these subsets, i.e. in the
language which is immanent to geometry and free of external means of
description (coordinates, curves). Such a description appears to be simple
and effective even in the case of complicated T-geometry. Antisymmetric
component of the world function generates appearance of additional metric
fields. This leads to appearance of three sorts of Christoffel symbols and
three sorts of geodesics. Three sorts of the first order tubes (future, past
and neutral) appear. If the fields connected with the antisymmetric
component are strong enough, the timelike first order tube has a finite
length in the timelike direction. It was shown earler that the symmetric
T-geometry explains non-relativistic quantum effects without a reference to
principles of quantum mechanics. One should expect that nonsymmetric
space-time T-geometry is also characteristic for microcosm, and it will be
useful in the elementary particle theory, because there is a series of
unexpected associations. First order tubes are associated with world tubes
of closed strings. Antisymmetry is associated with supersymmetry. Three
sorts of tubes are associated with three sorts of quarks. Limitation of the
tube in time direction is associated with confinement. At any rate, the
space-time T-geometry with additional parameters has more capacities, than
usual Riemannian geometry.
\end{abstract}

\newpage

\section{Introduction}

The most general geometry is obtained, if for its construction one uses the
purely metric conception of geometry. The conception of geometry (CG) is
defined as a method of construction of the standard (Euclidean) geometry.
There are several different conceptions of geometry, (i.e. methods of the
Euclidean geometry construction). They differ in amounts of numerical
information on geometry.

\bigskip

\noindent 
\centerline{$
\begin{array}{|c|c|c|}
\hline
& & \\
\text{\it title of CG} &
\begin{array}{c}
\text{\it non-numerical}  \\
\text{\it information}
\end{array}
&
\text{\it numerical information}
\\ & &
\\
\hline
\text{Euclidean CG} & \text{Euclidean axioms} & \emptyset
\\
\hline
\text{Riemannian CG}\label{tab} &
\begin{array}{c}
\text{Manifold, curve} \\
\text{coordinate system}
\end{array}
& n,\;\;g_{ik}\left( x\right)
\\
\hline
\begin{array}{c}
\text{topology-} \\ \text{metric CG}
\end{array}
&
\begin{array}{c}
\text{topological space,} \\
\text{curve}
\end{array}
&
\begin{array}{c}
\rho \left( P,Q\right) \geq 0,\\
\rho \left( P,Q\right) =0,\text{\ \ iff \
\thinspace }P=Q \\
\rho \left( P,Q\right) +\rho \left( Q,R\right) \geq \rho \left( P,R\right)
\end{array}  \\
\hline
\begin{array}{c}
\text{purely} \\ \text{metric CG}
\end{array}
& \emptyset  & \sigma \left( P,Q\right) =\frac{1}{2}\rho ^{2}
\left( P,Q\right) \in {\Bbb R} \\
\hline
\end{array}
$} \bigskip

Varying the numerical information at the fixed non-numerical one, we obtain
different geometries in the scope of the same conception of geometry.
Varying continuously numbers and functions, constituting numerical
information of CG, one obtains a continuous set of geometries, each of them
differs slightly from the narrow one. Any admissible value of numerical
information is attributed to some geometry in the scope of the given CG. The
purely metric CG contains only numerical information, and it generates the
most general class of geometries, suitable for applying them as the
space-time geometries.

Information on geometry is contained in the world function $\Sigma \left(
P,Q\right) =\frac{1}{2}\rho ^{2}\left( P,Q\right) $, where $\rho \left(
P,Q\right) $ is the distance between the points $P$ and $Q$. Conventionally
the world function $\Sigma $ is considered to be a symmetric real function
of points $P$ and $Q$ of the point set $\Omega $. 
\begin{equation}
\Sigma \left( P,Q\right) =\Sigma \left( Q,P\right) ,\qquad P,Q\in \Omega
\label{g1.1}
\end{equation}

In the papers \cite{R00,R01,R02} a nondegenerate geometry with the symmetric
world function was considered. It is a very general conception of geometry.
It admits one to construct such a space-time geometry, which explains
quantum phenomena without referring to principles of quantum mechanics, i.e.
as a geometric properties of the space-time.

The purely metric CG generates a set of geometries each of them is labelled
by the form of the world function $\Sigma $. Any such a geometry is called
nondegenerate geometry. T-geometry (tubular geometry) is another name of
nondegenerate geometry. In general, T-geometries are not axiomatized. Only
one of T-geometries (the standard geometry) has been axiomatized. The proper
Euclidean geometry ${\cal G}_{{\rm E}}$ plays the role of this standard
geometry. It is labelled by the world function $\Sigma _{{\rm E}}$. The
standard geometry has been studied very well. All geometric objects and all
statements of the standard geometry ${\cal G}_{{\rm E}}$ are expressed in
terms and only in terms of the world function $\Sigma _{{\rm E}}$. To obtain
geometric objects and statements of a T-geometry ${\cal G}$, labelled by the
world function $\Sigma $, it is sufficient one to take geometric objects and
statements of the standard geometry ${\cal G}_{{\rm E}}$, labelled by $%
\Sigma _{{\rm E}}$, and to replace $\Sigma _{{\rm E}}$ by $\Sigma $ there.
As far as all statements of ${\cal G}$ are expressed in terms of $\Sigma $,
one obtains statements of ${\cal G}$ as a modification of statements of $%
{\cal G}_{{\rm E}}$.

Such an approach admits one to describe any T-geometry ${\cal G}$ without
introducing axiomatization of ${\cal G}$. Let us yield necessary definitions.

\begin{opred}
T-geometry is the set of all statements about properties of all geometric
objects .
\end{opred}

The T-geometry is constructed on the point set $\Omega $ by giving the world
function $\Sigma $. The $\Sigma $-space $V=\{\Sigma ,\Omega \}$ is obtained
from the metric space after removal of the constraints, imposed on the
metric and introduction of the world function $\Sigma $ 
\begin{equation}
\Sigma \left( P,Q\right) =\frac{1}{2}\rho ^{2}\left( P,Q\right) ,\qquad
P,Q\in \Omega  \label{g2.1}
\end{equation}
instead of the metric $\rho $:

\begin{opred}
\label{dg2.1}$\Sigma $-space $V=\{\Sigma ,\Omega \}$ is nonempty set $\Omega 
$ of points $P$ with given on $\Omega \times \Omega $ real function $\Sigma $
\begin{equation}
\Sigma :\quad \Omega \times \Omega \rightarrow {\Bbb R},\qquad \Sigma
(P,P)=0,\qquad \forall P\in \Omega .  \label{g2.2}
\end{equation}
\end{opred}

The function $\Sigma $ is known as the world function \cite{S60}, or $\Sigma 
$-function. The metric $\rho $ may be introduced in $\Sigma $-space by means
of the relation (\ref{g2.1}). If $\Sigma $ is positive, the metric $\rho $
is also positive, but if $\Sigma $ is negative, the metric is imaginary.

\begin{opred}
\label{dg2.2}. Nonempty point set $\Omega ^{\prime }\subset \Omega $ of $%
\Sigma $-space $V=\{\Sigma ,\Omega \}$ with the world function $\Sigma
^{\prime }=\Sigma |_{\Omega ^{\prime }\times \Omega ^{\prime }}$, which is a
contraction $\Sigma $ on $\Omega ^{\prime }\times \Omega ^{\prime }$, is
called $\Sigma $-subspace $V^{\prime }=\{\Sigma ^{\prime },\Omega ^{\prime
}\}$ of $\Sigma $-space $V=\{\Sigma ,\Omega \}$.
\end{opred}

Further the world function $\Sigma ^{\prime }=\Sigma |_{\Omega ^{\prime
}\times \Omega ^{\prime }}$, which is a contraction of $\Sigma $ will be
denoted as $\Sigma $. Any $\Sigma $-subspace of $\Sigma $-space is a $\Sigma 
$-space. In T-geometry a geometric object ${\cal O}$ is described by means
of skeleton-envelope method. It means that any geometric object ${\cal O}$
is defined as follows.

\begin{opred}
\label{d1.7} Geometric object ${\cal O}$ is some $\Sigma $-subspace of $%
\Sigma $-space, which can be represented as a set of intersections and joins
of elementary geometric objects (EGO).
\end{opred}

\begin{opred}
\label{dd3.1}Elementary geometric object ${\cal E}\subset \Omega $ is a set
of zeros of the envelope function 
\begin{equation}
f_{{\cal P}^{n}}:\qquad \Omega \rightarrow {\Bbb R},\qquad {\cal P}%
^{n}\equiv \left\{ P_{0},P_{1},...P_{n}\right\} \in \Omega ^{n+1}
\label{b1.4}
\end{equation}
i.e. 
\begin{equation}
{\cal E}={\cal E}_{f}\left( {\cal P}^{n}\right) =\left\{ R|f_{{\cal P}%
^{n}}\left( R\right) =0\right\}  \label{b1.4a}
\end{equation}
The finite set ${\cal P}^{n}\subset \Omega $ of parameters of the envelope
function $f_{{\cal P}^{n}}$ is the skeleton of elementary geometric object
(EGO). The set ${\cal E}\subset \Omega $ of points forming EGO is called the
envelope of its skeleton ${\cal P}^{n}$. The envelope function $f_{{\cal P}%
^{n}}$ is an algebraic function of $s$ arguments $w=\left\{
w_{1},w_{2},...w_{s}\right\} $, $s=(n+2)(n+1)$. Each of arguments $%
w_{k}=\Sigma \left( Q_{k},L_{k}\right) $ is a $\Sigma $-function of two
arguments $Q_{k},L_{k}\in \left\{ R,{\cal P}^{n}\right\} $.
\end{opred}

For continuous T-geometry the envelope ${\cal E}$ is usually a continual set
of points. The envelope function $f_{{\cal P}^{n}}$, determining EGO is a
function of the running point $R\in \Omega $ and of parameters ${\cal P}%
^{n}\in \Omega ^{n+1}$. Thus, any elementary geometric object is determined
by its skeleton ${\cal P}^{n}$ and by the form of the envelope function $f_{%
{\cal P}^{n}}$.

Let one investigates T-geometry on the $\Sigma $-space $V=\left\{ \Sigma
,\Omega \right\} $. For some special choice $\Sigma _{{\rm E}}$ of $\Sigma $%
-function, the $\Sigma $-space $V$ turns to a $\Sigma $-subspace $V_{{\rm E}%
}^{\prime }=\left\{ \Sigma _{{\rm E}},\Omega \right\} $ of a $n$-dimensional
proper Euclidean space $V_{{\rm E}}=\left\{ \Sigma _{{\rm E}},\Omega _{{\rm E%
}}\right\} $, \ $\Omega \subset \Omega _{{\rm E}}$. (It will be shown). Then
all relations between geometric objects in $V_{{\rm E}}^{\prime }$ are
relations of proper Euclidean geometry. Replacement of $\Sigma _{{\rm E}}$
by $\Sigma $ means a deformation of $V_{{\rm E}}^{\prime }$, because world
function $\Sigma $ describes distances between two points, and change of
these distances is a deformation of the space. We shall use concept of
deformation in a wide meaning, including in this term any increase and any
reduction of number of points in the set $\Omega $. Then any transition from 
$\left\{ \Sigma _{{\rm E}},\Omega _{{\rm E}}\right\} $ to $\left\{ \Sigma
,\Omega \right\} $ is a deformation of $\left\{ \Sigma _{{\rm E}},\Omega _{%
{\rm E}}\right\} $.

Let us write Euclidean relations between geometric objects in $V_{{\rm E}%
}^{\prime }$ in the $\sigma $-immanent form (i.e. in the form, which
contains references only to geometrical objects and $\Sigma $-function).
Changing the world function $\Sigma _{{\rm E}}$ by $\Sigma $ in these
relations, one obtains the relations between the geometric objects in the $%
\Sigma $-space $V=\left\{ \Sigma ,\Omega \right\} $.

Thus, geometry in the proper Euclidean space is known very well, and one
uses deformation, described by world function, to establish T-geometry in
arbitrary $\Sigma $-space. Considering deformations of Euclidean space, one
goes around the problem of axiomatics in the $\Sigma $-space $V=\left\{
\Sigma ,\Omega \right\} $. One uses only Euclidean axiomatics. T-geometry of
arbitrary $\Sigma $-space is obtained as a result of ''deformation of proper
Euclidean geometry''. This point is very important, because axiomatics of
arbitrary T-geometry is very complicated. It is relatively simple only for
highly symmetric spaces. Investigation of arbitrary deformations is much
simpler, than investigations of arbitrary axiomatics. Formally, a work with
deformations of $\Sigma $-spaces is manipulations with the world function.
These manipulations may be carried out without mention of space deformations.

Description of EGOs by means (\ref{b1.4}) is carried out in the
deform-invariant form (invariant with respect to $\Sigma $-space
deformations). The envelope function $f_{{\cal P}^{n}}$ as a function of
arguments $w_{k}=\Sigma \left( Q_{k},L_{k}\right) ,$ $Q_{k},L_{k}\in \left\{
R,{\cal P}^{n}\right\} $ does not depend on the form of the world function $%
\Sigma $. Thus, definition of the envelope function is invariant with
respect to deformations (deform-invariant), and the envelope function
determines any EGO in all $\Sigma $-spaces at once.

Let ${\cal E}_{{\rm E}}$ be EGO in the Euclidean geometry ${\cal G}_{{\rm E}%
} $. Let ${\cal E}_{{\rm E}}$ be described by the skeleton ${\cal P}^{n}$
and the envelope function $f_{{\cal P}^{n}}$ in the $\Sigma $-space $V_{{\rm %
E}}=\left\{ \Sigma _{{\rm E}},\Omega \right\} $. Then the EGO ${\cal E}$ in
the T-geometry ${\cal G}$, described by the same skeleton ${\cal P}^{n}$ and
the same envelope function $f_{{\cal P}^{n}}$ in the $\Sigma $-space $%
V=\left\{ \Sigma ,\Omega \right\} $, is an analog in ${\cal G}$ of the
Euclidean EGO ${\cal E}_{{\rm E}}$. T-geometry ${\cal G}$ may be considered
to be a result of deformation of the Euclidean geometry ${\cal G}_{{\rm E}}$%
, when distances $\sqrt{\Sigma \left( P,Q\right) +\Sigma \left( Q,P\right) }$
between the pairs of points $P$ and $Q$ are changed. At such a deformation
the Euclidean EGO ${\cal E}_{{\rm E}}$ transforms to its analog ${\cal E}$.

The Euclidean space has the most powerful group of motion, and the same
envelope ${\cal E}_{{\rm E}}$ may be generated by the envelope function $f_{%
{\cal P}^{n}}$ with different values ${\cal P}_{(1)}^{n},$ ${\cal P}%
_{(2)}^{n},...$ of the skeleton ${\cal P}^{n}$, or even by another envelope
function $f_{\left( 1\right) {\cal Q}^{m}}$. It means that the Euclidean EGO 
${\cal E}_{{\rm E}}$ may have several analogs ${\cal E}_{(1)},{\cal E}%
_{(2)},...$ in the geometry ${\cal G}$. In other words, deformation of the
Euclidean space may split EGOs, (but not only deform them). Note that the
splitting may be interpreted as a kind of deformation.

In the topologo-metric conception of geometry \cite{T59,ABN86,BGP92} one
uses concept of curve as the continuous mapping 
\begin{equation}
{\cal L}:\;\;\left[ 0,1\right] \rightarrow \Omega ,\qquad \left[ 0,1\right]
\subset {\Bbb R},  \label{g1.4}
\end{equation}
The curve ${\cal L}\left( \left[ 0,1\right] \right) \subset \Omega $ is
considered to be an important geometrical object of geometry. From point of
view of T-geometry the set of points ${\cal L}\left( \left[ 0,1\right]
\right) \subset \Omega $ cannot be considered to be EGO, because the mapping
(\ref{g1.4}) is not deform-invariant. Indeed, let us consider a sphere $%
{\cal S}_{P_{0}P_{1}}$, passing through the point $P_{1}$ and having its
center at the point $P_{0}$. It is described by the envelope function 
\begin{equation}
f_{P_{0}P_{1}}\left( R\right) =\sqrt{\Sigma \left( P_{0},R\right) +\Sigma
\left( R,P_{0}\right) }-\sqrt{\Sigma \left( P_{1},R\right) +\Sigma \left(
R,P_{1}\right) }  \label{g1.4a}
\end{equation}
In the two-dimensional proper Euclidean space the envelope function (\ref
{g1.4a}) describes a one-dimensional circumference ${\cal L}_{1}$, whereas
in the three-dimensional proper Euclidean space the envelope function (\ref
{g1.4a}) describes a two-dimensional sphere ${\cal S}_{2}$. The point set $%
{\cal L}_{1}$ can be represented as the continuous mapping (\ref{g1.4}),
whereas the surface ${\cal S}_{2}$ cannot. Transition from two-dimensional
Euclidean space to three-dimensional Euclidean space is a space deformation.
Thus, deformation of the $\Sigma $-space may destroy the property of EGO of
being a curve (\ref{g1.4}).

Application of objects, defined by the property (\ref{g1.4}) for
investigation of T-geometries is inconvenient, because the T-geometry
investigation is founded on deform-invariant methods. Formally, one cannot
choose appropriate envelope function for description of the set (\ref{g1.4}%
), because the envelope function is deform-invariant, whereas the set (\ref
{g1.4}) is not. Hence, (\ref{g1.4}) is incompatible with the definition \ref
{dd3.1} of EGO. In other words, the point set (\ref{g1.4}) is not
deform-invariant.

There are at least two approach to geometry:

(1) The geometry is supposed to be completely described by the world
function. The standard (proper Euclidean) geometry is such a kind of
geometry. Then any possible geometry is one of T-geometries.

(2) The geometry is supposed to be described by the world function
incompletely. There are functions of three and more number of points, which
describe such a geometry (for instance, functions $\Sigma _{3}\left(
P_{0},P_{2},P_{2}\right) $,...). We shall refer to such a geometry as the
fortified geometry (FG). The standard (proper Euclidean) geometry is
described by the world function $\Sigma $. In the Euclidean geometry the
function $\Sigma _{3}$ and other functions have some trivial form (for
instance, $\Sigma _{3{\rm E}}\equiv 0$).

The T-geometry pretends only to description of usual geometry, but not to
fortified geometry (FG). It is not quite clear, whether the papers \cite
{T59,ABN86,BGP92} pretend to description of the fortified geometry, or they
are restricted by the usual one. If they pretend only to a description of
usual geometry, then the concept of the curve (\ref{g1.4}) is only a method
for description of geometry (something like a coordinate system), but not an
object of geometry. All real properties of geometry (not those of the
description method) should be formulated in the form invariant with respect
to a change of the description method (in particular, in the form, where the
concept of the curve is not mentioned at all).

If the papers \cite{T59,ABN86,BGP92} pretend to a description of some
special cases of FG, which cannot be described by T-geometry, then a use of
the concept (\ref{g1.4}) may be justified. But such an approach seems to be
too complicated. It seems to be simpler one to have investigated the usual
geometry completely, and thereafter to begin investigation of a more
complicated case of fortified geometry.

The next step in the investigation of the usual geometry is a refuse from
the constraint (\ref{g1.1}). Such a geometry will be referred to as a
nonsymmetric T-geometry. The nonsymmetric T-geometry can be investigated by
the same methods, as the symmetric one. The world function $\Sigma $ in the
nonsymmetric T-geometry is presented in the form 
\begin{eqnarray}
\Sigma \left( P,Q\right) &=&G\left( P,Q\right) +A\left( P,Q\right) ,\qquad
P,Q\in \Omega  \label{b2.2} \\
G\left( P,Q\right) &=&G\left( Q,P\right) ,\qquad A\left( P,Q\right)
=-A\left( Q,P\right)  \label{g1.2} \\
G\left( P,Q\right) &=&\frac{1}{2}\left( \Sigma \left( P,Q\right) +\Sigma
\left( Q,P\right) \right) ,  \label{c2.3a} \\
A\left( P,Q\right) &=&\frac{1}{2}\left( \Sigma \left( P,Q\right) -\Sigma
\left( Q,P\right) \right)  \label{g1.3}
\end{eqnarray}
where $G$ denotes the symmetric part of the world function $\Sigma $,
whereas $A$ denotes its antisymmetric part.

Motives for consideration of nonsymmetric T-geometry are as follows. In the
symmetric T-geometry the distance from the point $P$ to the point $Q$ is the
same as the distance from the point $Q$ to the point $P$. In the asymmetric
T-geometry it is not so. Apparently, it is not important for spacelike
distances in the space-time, because it can be tested experimentally for
spacelike distances. In the case, when interval between points $P$ and $Q$
is timelike, one uses watch to measure this interval. But the watch can
measure the time interval only in one direction, and one cannot be sure that
the time interval is the same in opposite direction.

If the antisymmetric part $A$ of the world function does not vanish, it
means that the future and the past are not equivalent geometrically. We do
not insist that this fact takes place, but we admit this. It is useful to
construct a nonsymmetric T-geometry, to apply it to the space-time and to
obtain the corollaries of asymmetry which could be tested experimentally.
The symmetrical part of the world function generates the field of the metric
tensor $g_{ik}$. In a like way the antisymmetric part generates some vector
force filed $a_{i}$. Maybe, existence of this field can be tested
experimentally. For construction of nonsymmetric T-geometry one does not
need to make any additional supposition. It is sufficient to remove the
constraint (\ref{g1.1}) and to apply mathematical technique developed for
the symmetric T-geometry with necessary modifications.

Besides, there is a hope that nonsymmetric T-geometry will be useful in the
elementary particle theory, where the main object is a superstring. The
first order tubes (main objects of T-geometry) are associated with world
tubes of strings. In the nonsymmetric T-geometry antisymmetric variables
appear. They are absent in the conventional symmetric T-geometry, but
antisymmetric variables are characteristic for the superstring theory.

Two important general remarks.

1. Nonsymmetric T-geometry, as well as the symmetric one, is considered on
an arbitrary set $\Omega $ of points $P$. It is formulated in the scope of
the purely metric conception of geometry, which is very simple, because it
uses only very simple tools for the geometry description. The T-geometry
formulated in terms of the world function $\Sigma $ and finite subsets $%
{\cal P}^{n}\equiv \left\{ P_{0},P_{1},...,P_{n}\right\} $ of the set $%
\Omega $. Mathematically it means, that the purely metric CG uses only
mappings

\begin{equation}
m_{n}:\;\;I_{n}\rightarrow \Omega ,\qquad I_{n}\equiv \left\{
0,1,...n\right\} \subset \left\{ 0\right\} \cup {\Bbb N}  \label{c2.8c}
\end{equation}
whereas the topology-metric CG uses much more complicated mappings (\ref
{g1.4}), known as curves ${\cal L}$. Both mappings (\ref{c2.8c}) and (\ref
{g1.4}) are methods of the geometry description (and construction). But the
method (\ref{c2.8c}) is much simpler. It can be studied exhaustively,
whereas the set of mappings (\ref{g1.4}) cannot.

2. The nonsymmetric T-geometry will be mainly interpreted as a symmetric
T-geometry determined by the two-point scalar $G\left( P,Q\right) $ with
some additional metric structures, introduced to the symmetric geometry by
means of the additional two-point scalar $A\left( P,Q\right) $. For
instance, in the symmetric space-time T-geometry the world line of a free
particle is described by a geodesic. In the nonsymmetric space-time
T-geometry there are, in general, several different types of geodesics. This
fact may be interpreted in the sense, that a free particle has some internal
degrees of freedom, and it may be found in different states. In these
different states the free particle interacts with the force fields,
generated by the two-point scalar $A\left( P,Q\right) $, differently.
Several different types of geodesics are results of this interaction.

In the second section formulation of the asymmetric T-geometry in a
coordinate-free form is presented. In the third section properties of tubes
in the $\Sigma $-space are considered. The fourth section is devoted to
consideration of $\Sigma $-space, given on a manifold. Derivatives of the
world function at the coincidence points $x=x^{\prime }$ are considered in
the fifth section. The two-point curvature tensor is investigated in the
sixth section. The seventh and eighth sections are devoted to gradient lines
and their $\sigma $-immanent description. Conditions of the first order tube
degeneration are considered in the ninth section. In the tenth section
examples of the first order tubes are considered.

\section{$\Sigma $-space and T-geometry. Coordinate-free description}

\begin{opred}
\label{dg2.3}. $\Sigma $-space $V=\{\Sigma ,\Omega \}$ is called
isometrically embeddable in $\Sigma $-space $V^{\prime }=\{\Sigma ^{\prime
},\Omega ^{\prime }\}$, if there exists such a monomorphism $f:\Omega
\rightarrow \Omega ^{\prime }$, that $\Sigma (P,Q)=\Sigma ^{\prime
}(f(P),f(Q))$,\quad $\forall P,\forall Q\in \Omega ,\quad f(P),f(Q)\in
\Omega ^{\prime }$,
\end{opred}

Any $\Sigma $-subspace $V^{\prime }$ of $\Sigma $-space $V=\{\Sigma ,\Omega
\}$ is isometrically embeddable in it.

\begin{opred}
\label{dd1.1}. Two $\Sigma $-spaces $V=\{\Sigma ,\Omega \}$ and $V^{\prime
}=\{\Sigma ^{\prime },\Omega ^{\prime }\}$ are called to be isometric
(equivalent), if $V$ is isometrically embeddable in $V^{\prime }$, and $%
V^{\prime }$ is isometrically embeddable in $V$.
\end{opred}

\begin{opred}
\label{dd1.2}The $\Sigma $-space $M=\{\Sigma ,\Omega \}$ is called a finite $%
\Sigma $-space, if the set $\Omega $ contains a finite number of points.
\end{opred}

\begin{opred}
\label{dd1.3}. The $\Sigma $-subspace $M_{n}({\cal P}^{n})=\{\Sigma ,{\cal P}%
^{n}\}$ of the $\Sigma $-space $V=\{\Sigma ,\Omega \}$, consisting of $n+1$
points ${\cal P}^{n}=\left\{ P_{0},P_{1},...,P_{n}\right\} $ is called the $%
n $th order $\Sigma $-subspace .
\end{opred}

The T-geometry is a set of all propositions on properties of $\Sigma $%
-subspaces of $\Sigma $-space $V=\{\Sigma ,\Omega \}$. Presentation of
T-geometry is produced in the language, containing only references to $%
\Sigma $-function and constituents of $\Sigma $-space, i.e. to its $\Sigma $%
-subspaces.

\begin{opred}
\label{d3} A description is called $\sigma $-immanent, if it does not
contain any references to objects or concepts other, than finite subspaces
of the $\Sigma $-space and its world function (metric).
\end{opred}

$\sigma $-immanence of description provides independence of the description
on the method of description. In this sense the $\sigma $-immanence of a
description in T-geometry reminds the concept of covariance in Riemannian
geometry. Covariance of some geometrical relation in Riemannian geometry
means that the considered relation is valid in all coordinate systems and,
hence, describes only the properties of the Riemannian geometry in itself.
Covariant description provides cutting-off from the coordinate system
properties, considering the relation in all coordinate systems at once. The $%
\sigma $-immanence provides truncation from the methods of description by
absence of a reference to objects, which do not relate to geometry in itself
(coordinate system, concept of curve, dimension).

The idea of constructing the T-geometry is very simple. All relations of
proper Euclidean geometry are written in the $\sigma $-immanent form and
declared to be valid for any $\Sigma $-function. This results that any
relation of proper Euclidean geometry corresponds to some relation of
T-geometry. It is important that in the relations, declared to be relations
of T-geometry, only the properties (\ref{g2.2}) were used. The special
properties of the Euclidean $\Sigma $-function are not to be taken into
account. The metric part of these relations was formulated and proved by
K.~Menger \cite{M28}. Let us present this result in our designations in the
form of the theorem 
\begin{theorem}
\label{t1} The symmetric $\Sigma  $-space $V=\{\Sigma  ,\Omega \}$ is 
isometrically
embeddable in $n$-dimensional proper Euclidean space $E_n$, if and only
if any $(n+2)$th order $\Sigma $-subspace $M({\cal P}^{n+2})\subset \Omega $
is isometrically embeddable in $E_n$.
\end{theorem}

Unfortunately, the formulation of this theorem is not $\sigma $-immanent, as
far as it contains a reference to $n$-dimensional Euclidean space $E_{n}$
which is not defined $\sigma $-immanently. A more constructive version of
the $\Sigma $-space Euclideaness conditions will be formulated, as soon as
we construct a necessary mathematical technique.

The basic elements of T-geometry are finite $\Sigma $-subspaces $M_{n}({\cal %
P}^{n})$, i.e. finite sets 
\begin{equation}
{\cal P}^{n}=\{P_{0},P_{1},\ldots ,P_{n}\}\subset \Omega  \label{a1.9}
\end{equation}

The simplest finite subset is a nonzero vector $\overrightarrow{{\cal P}^{1}}%
={\bf P}_{0}{\bf P}_{1}\equiv \overrightarrow{P_{0}P_{1}}$. The vector $%
\overrightarrow{P_{0}P_{1}}$ is an ordered set of two points $\left\{
P_{0},P_{1}\right\} $. The scalar product $\left( {\bf P}_{0}{\bf P}_{1}.%
{\bf Q}_{0}{\bf Q}_{1}\right) $ of two vectors ${\bf P}_{0}{\bf P}_{1}$ and $%
{\bf Q}_{0}{\bf Q}_{1}$

\begin{equation}
\left( {\bf P}_{0}{\bf P}_{1}.{\bf Q}_{0}{\bf Q}_{1}\right) =\Sigma \left(
P_{0},Q_{1}\right) -\Sigma \left( P_{1},Q_{1}\right) -\Sigma \left(
P_{0},Q_{0}\right) +\Sigma \left( P_{1},Q_{0}\right)  \label{a2.1a}
\end{equation}
is the main construction of T-geometry, and we substantiate this definition.

$\sigma $-immanent expression for scalar product $\left( {\bf P}_{0}{\bf P}%
_{1}.{\bf Q}_{0}{\bf Q}_{1}\right) $ of two vectors ${\bf P}_{0}{\bf P}_{1}$
and ${\bf Q}_{0}{\bf Q}_{1}$ in the proper Euclidean space has the form (\ref
{a2.1a}). This relation can be easily proved as follows.

In the proper Euclidean space three vectors ${\bf P}_{0}{\bf P}_{1}$, ${\bf P%
}_{0}{\bf Q}_{1}$, and ${\bf P}_{1}{\bf Q}_{1}$ are coupled by the relation 
\begin{equation}
|{\bf P}_{1}{\bf Q}_{1}|^{2}=|{\bf P}_{0}{\bf Q}_{1}-{\bf P}_{0}{\bf P}%
_{1}|^{2}=|{\bf P}_{0}{\bf P}_{1}|^{2}+|{\bf P}_{0}{\bf Q}_{1}|^{2}-2({\bf P}%
_{0}{\bf P}_{1}.{\bf P}_{0}{\bf Q}_{1})  \label{g2.3}
\end{equation}
where $({\bf P}_{0}{\bf P}_{1}.{\bf P}_{0}{\bf Q}_{1})$ denotes the scalar
product of two vectors ${\bf P}_{0}{\bf P}_{1}$ and ${\bf P}_{0}{\bf Q}_{1}$
in the proper Euclidean space, and $|{\bf P}_{0}{\bf P}_{1}|^{2}\equiv ({\bf %
P}_{0}{\bf P}_{1}.{\bf P}_{0}{\bf P}_{1})$. It follows from (\ref{g2.3}) 
\begin{equation}
({\bf P}_{0}{\bf P}_{1}.{\bf P}_{0}{\bf Q}_{1})={\frac{1}{2}}\left( |{\bf P}%
_{0}{\bf Q}_{1}|^{2}+|{\bf P}_{0}{\bf P}_{1}|^{2}-|{\bf P}_{1}{\bf Q}%
_{1}|^{2}\right)  \label{g2.4}
\end{equation}
Substituting the point $Q_{1}$ by $Q_{0}$ in (\ref{g2.4}), one obtains 
\begin{equation}
({\bf P}_{0}{\bf P}_{1}.{\bf P}_{0}{\bf Q}_{0})={\frac{1}{2}}\{|{\bf P}_{0}%
{\bf Q}_{0}|^{2}+|{\bf P}_{0}{\bf P}_{1}|^{2}-|{\bf P}_{1}{\bf Q}_{0}|^{2}\}
\label{g2.5}
\end{equation}
Subtracting (\ref{g2.5}) from (\ref{g2.4}) and using the properties of the
scalar product in the proper Euclidean space, one obtains 
\begin{equation}
({\bf P}_{0}{\bf P}_{1}.{\bf Q}_{0}{\bf Q}_{1})={\frac{1}{2}}\{|{\bf P}_{0}%
{\bf Q}_{1}|^{2}+|{\bf P}_{1}{\bf Q}_{0}|^{2}-|{\bf P}_{0}{\bf Q}_{0}|^{2}-|%
{\bf P}_{1}{\bf Q}_{1}|^{2}\}  \label{g2.6}
\end{equation}
Taking into account that in the proper Euclidean geometry $|{\bf P}_{0}{\bf Q%
}_{1}|^{2}=2\Sigma \left( P_{0},Q_{1}\right) =2G\left( P_{0},Q_{1}\right) $,
one obtains the relation\ (\ref{a2.1a}) from the relation (\ref{g2.6}).

In the Euclidean geometry the world function is symmetric, and the order of
arguments in the rhs of (\ref{a2.1a}) is not essential. In the asymmetric
T-geometry the order of arguments in the rhs of (\ref{a2.1a}) is essential.
The order has been chosen in such a way that 
\begin{eqnarray}
\left( {\bf P}_{0}{\bf P}_{1}.{\bf Q}_{0}{\bf Q}_{1}\right) _{{\rm s}}
&\equiv &\frac{1}{2}\left( \left( {\bf P}_{0}{\bf P}_{1}.{\bf Q}_{0}{\bf Q}%
_{1}\right) +\left( {\bf Q}_{0}{\bf Q}_{1}.{\bf P}_{0}{\bf P}_{1}\right)
\right)  \nonumber \\
&=&G\left( P_{0},Q_{1}\right) -G\left( P_{1},Q_{1}\right) -G\left(
P_{0},Q_{0}\right) +G\left( P_{1},Q_{0}\right)  \label{g2.7}
\end{eqnarray}
\begin{eqnarray}
\left( {\bf P}_{0}{\bf P}_{1}.{\bf Q}_{0}{\bf Q}_{1}\right) _{{\rm a}}
&\equiv &\frac{1}{2}\left( \left( {\bf P}_{0}{\bf P}_{1}.{\bf Q}_{0}{\bf Q}%
_{1}\right) -\left( {\bf Q}_{0}{\bf Q}_{1}.{\bf P}_{0}{\bf P}_{1}\right)
\right)  \nonumber \\
&=&A\left( P_{0},Q_{1}\right) -A\left( P_{1},Q_{1}\right) -A\left(
P_{0},Q_{0}\right) +A\left( P_{1},Q_{0}\right)  \label{g2.8}
\end{eqnarray}
It follows from (\ref{a2.1a}) that 
\begin{equation}
\left( {\bf P}_{0}{\bf P}_{1}.{\bf Q}_{0}{\bf Q}_{1}\right) =-\left( {\bf P}%
_{1}{\bf P}_{0}.{\bf Q}_{0}{\bf Q}_{1}\right) ,\qquad \left( {\bf P}_{0}{\bf %
P}_{1}.{\bf Q}_{0}{\bf Q}_{1}\right) =-\left( {\bf P}_{0}{\bf P}_{1}.{\bf Q}%
_{1}{\bf Q}_{0}\right)  \label{g2.9}
\end{equation}
Thus, the scalar product $\left( {\bf P}_{0}{\bf P}_{1}.{\bf Q}_{0}{\bf Q}%
_{1}\right) $ of two vectors ${\bf P}_{0}{\bf P}_{1}$ and ${\bf Q}_{0}{\bf Q}%
_{1}$ is antisymmetric with respect to permutation $P_{0}\leftrightarrow
P_{1}$ of points determining the vector ${\bf P}_{0}{\bf P}_{1}$, as well as
with respect to permutation $Q_{0}\leftrightarrow Q_{1}$.

\begin{opred}
\label{d.2.6} The finite $\Sigma $-space $M_{n}({\cal P}^{n})=\{\Sigma ,%
{\cal P}^{n}\}$ is called oriented $\overrightarrow{M_{n}({\cal P}^{n})}$,
if the order of its points ${\cal P}^{n}=\{P_{0},P_{1},\ldots P_{n}\}$ is
fixed.
\end{opred}

\begin{opred}
\label{d3.1.2b}. The $n$th order multivector $m_{n}$ is the mapping 
\begin{equation}
m_{n}:\qquad I_{n}\rightarrow \Omega ,\qquad I_{n}\equiv \left\{
0,1,...,n\right\}  \label{g2.10}
\end{equation}
\end{opred}

The set $I_{n}$ has a natural ordering, which generates an ordering of
images $m_{n}(k)\in \Omega $ of points $k\in I_{n}$. The ordered list of
images of points in $I_{n}$ has one-to-one connection with the multivector
and may be used as the multivector descriptor. Different versions of the
point list will be used for writing the $n$th order multivector descriptor: 
\[
\overrightarrow{P_{0}P_{1}...P_{n}}\equiv {\bf P}_{0}{\bf P}_{1}...{\bf P}%
_{n}\equiv \overrightarrow{{\cal P}^{n}} 
\]
Originals of points $P_{k}$ in $I_{n}$ are determined by the order of the
point $P_{k}$ in the list of descriptor. Index of the point $P_{k}$ has
nothing to do with the original of $P_{k}$. Further we shall use descriptor $%
\overrightarrow{P_{0}P_{1}...P_{n}}$ of the multivector instead of the
multivector. In this sense the $n$th order multivector $\overrightarrow{%
P_{0}P_{1}...P_{n}}$ in the $\Sigma $-space $V=\{\Sigma ,\Omega \}$ may be
defined as the ordered set $\{P_{l}\},\quad l=0,1,\ldots n$ of $n+1$ points $%
P_{0},P_{1},...,P_{n}$, belonging to the $\Sigma $-space $V$. Some points
may be identical. The point $P_{0}$ is the origin of the multivector $%
\overrightarrow{P_{0}P_{1}...P_{n}}$. Image $m_{n}\left( I_{n}\right) $ of
the set $I_{n}$ contains $k$ points ($k\leq n+1).$ The set of all $n$th
order multivectors $m_{n}$ constitutes the set $\Omega
^{n+1}=\bigotimes\limits_{k=1}^{n+1}\Omega $, and any multivector $%
\overrightarrow{{\cal P}^{n}}\in \Omega ^{n+1}$.

\begin{opred}
\label{d3.1.6b}. The scalar $\Sigma $-product $(\overrightarrow{{\cal P}^{n}}%
.\overrightarrow{{\cal Q}^{n}})$ of two $n$th order multivectors $%
\overrightarrow{{\cal P}^{n}}$ and $\overrightarrow{{\cal Q}^{n}}$ is the
real number 
\begin{equation}
(\overrightarrow{{\cal P}^{n}}.\overrightarrow{{\cal Q}^{n}})=\det \Vert (%
{\bf P}_{0}{\bf P}_{i}.{\bf Q}_{0}{\bf Q}_{k})\Vert
,\,\,\,\,\,\,\,\,\,\,\,i,k=1,2,...n  \label{g2.11}
\end{equation}
\begin{eqnarray}
({\bf P}_{0}{\bf P}_{i}.{\bf Q}_{0}{\bf Q}_{k}) &\equiv &\Sigma
(P_{0},Q_{k})+\Sigma (P_{i},Q_{0})-\Sigma (P_{0},Q_{0})-\Sigma (P_{i},Q_{k}),
\label{g2.12} \\
P_{0},P_{i},Q_{0},Q_{k} &\in &\Omega ,\qquad \overrightarrow{{\cal P}^{n}},%
\overrightarrow{{\cal Q}^{n}}\in \Omega ^{n+1}  \nonumber
\end{eqnarray}
\end{opred}

Operation of permutation of the multivector points can be effectively
defined in the $\Sigma $-space. Let us consider two $n$th order multivectors 
$\overrightarrow{{\cal P}^{n}}=\overrightarrow{P_{0}P_{1}P_{2}...P_{n}}$ and 
$\overrightarrow{{\cal P}_{(k\leftrightarrow l)}^{n}}=\overrightarrow{%
P_{0}P_{1}...P_{k-1}P_{l}P_{k+1}...P_{l-1}P_{k}P_{l+1}...P_{n}},\;\;(n\geq
1) $, which is a result of permutation of points $P_{k},$ $P_{l},\;\;(k<l)$.
The scalar $\Sigma $-product $(\overrightarrow{{\cal P}^{n}}.\overrightarrow{%
{\cal Q}^{n}})$ is defined by the relation (\ref{g2.11}). If $k=0$, then the
permutation $P_{0}\leftrightarrow P_{l}$ changes the sign of the $l$th row
of the determinant (\ref{g2.11}), and 
\begin{equation}
(\overrightarrow{{\cal P}^{n}}.\overrightarrow{{\cal Q}^{n}})=-(%
\overrightarrow{{\cal P}_{(0\leftrightarrow l)}^{n}}.\overrightarrow{{\cal Q}%
^{n}})\qquad l=1,2,...n,\qquad \forall \overrightarrow{{\cal Q}^{n}}\in
\Omega ^{n+1},  \label{g2.13}
\end{equation}
Let us show this for the case $l=1$. Making transposition $%
P_{0}\leftrightarrow P_{1}$ in the determinant (\ref{g2.11}), one obtains
for elements of the first row 
\begin{equation}
({\bf P}_{0}{\bf P}_{1}.{\bf Q}_{0}{\bf Q}_{k})\rightarrow ({\bf P}_{1}{\bf P%
}_{0}.{\bf Q}_{0}{\bf Q}_{k})=-({\bf P}_{0}{\bf P}_{1}.{\bf Q}_{0}{\bf Q}%
_{k}),\qquad k=1,2,...,n  \label{g2.13a}
\end{equation}
For remaining $n-1$ rows one obtains 
\begin{equation}
({\bf P}_{0}{\bf P}_{i}.{\bf Q}_{0}{\bf Q}_{k})\rightarrow ({\bf P}_{1}{\bf P%
}_{i}.{\bf Q}_{0}{\bf Q}_{k}),\qquad i=2,3,...n,\qquad k=1,2,...,n
\label{g2.13b}
\end{equation}
Let us take into account identity 
\begin{equation}
({\bf P}_{1}{\bf P}_{i}.{\bf Q}_{0}{\bf Q}_{k})-({\bf P}_{1}{\bf P}_{0}.{\bf %
Q}_{0}{\bf Q}_{k})\equiv ({\bf P}_{0}{\bf P}_{i}.{\bf Q}_{0}{\bf Q}_{k})
\label{g2.13c}
\end{equation}
which follows from definition (\ref{g2.12}) of the scalar $\Sigma $-product.
Then subtracting the first row from all other rows, one obtains that $n-1$
rows $i=2,3,...n$ of the transformed determinant coincide with corresponding 
$n-1$ rows of the determinant (\ref{g2.11}). This proves the relation (\ref
{g2.13}) for $l=1.$ In the same way one can prove (\ref{g2.13}) for other
values of $l$.

If $k\neq 0$, and $k<l$, the permutation $P_{k}\leftrightarrow P_{l}$
transposes the $l$th row and the $k$th row of the determinant (\ref{g2.11}).
The determinant changes its sign and, taking into account (\ref{g2.13}), one
obtains 
\begin{equation}
(\overrightarrow{{\cal P}^{n}}.\overrightarrow{{\cal Q}^{n}})=-(%
\overrightarrow{{\cal P}_{(k\leftrightarrow l)}^{n}}.\overrightarrow{{\cal Q}%
^{n}})\qquad k\neq l,\qquad l,k=0,1,2,...n,\qquad \forall \overrightarrow{%
{\cal Q}^{n}}\in \Omega ^{n+1},  \label{g2.14}
\end{equation}
As far as the relation (\ref{g2.14}) is valid for permutation of any two
points of the multivector $\overrightarrow{{\cal P}^{n}}$ and for any
multivector $\forall \overrightarrow{{\cal Q}^{n}}\in \Omega ^{n+1},$ one
may write 
\begin{equation}
\overrightarrow{{\cal P}_{(i\leftrightarrow k)}^{n}}=-\overrightarrow{{\cal P%
}^{n}},\qquad i,k=0,1,...n,\qquad i\neq k,\qquad n\geq 1.  \label{g2.15}
\end{equation}
Thus, a change of the $n$th order multivector sign $(n\geq 1)$
(multiplication by the number $a=-1$) may be always defined as an odd
permutation of points.

Let us consider the relation 
\begin{equation}
\overrightarrow{{\cal P}^{n}}T\overrightarrow{{\cal R}^{n}}:\qquad (%
\overrightarrow{{\cal P}^{n}}.\overrightarrow{{\cal Q}^{n}})=(%
\overrightarrow{{\cal R}^{n}}.\overrightarrow{{\cal Q}^{n}})\wedge (%
\overrightarrow{{\cal Q}^{n}}.\overrightarrow{{\cal P}^{n}})=(%
\overrightarrow{{\cal Q}^{n}}.\overrightarrow{{\cal R}^{n}}),\qquad \forall 
\overrightarrow{{\cal Q}^{n}}\in \Omega ^{n+1},  \label{g2.16}
\end{equation}
between two $n$th order multivectors $\overrightarrow{{\cal P}^{n}}\in
\Omega ^{n+1}$ and $\overrightarrow{{\cal R}^{n}}\in \Omega ^{n+1}.$ The
relation (\ref{g2.16}) is reflexive, symmetric and transitive, and it may be
considered as an equivalence relation.

\begin{opred}
\label{d3.1.6a}. Two $n$th order multivectors $\overrightarrow{{\cal P}^n}%
\in \Omega ^{n+1}$ and $\overrightarrow{{\cal R}^n}\in \Omega ^{n+1}$ are
equivalent $\overrightarrow{{\cal P}^n}=\overrightarrow{{\cal R}^n}$, if the
relations (\ref{g2.16}) takes place.
\end{opred}

\begin{opred}
\label{d3.1.6b0}. If the $n$th order multivector $\overrightarrow{{\cal N}%
^{n}}$ satisfies the relations 
\begin{equation}
(\overrightarrow{{\cal N}^{n}}.\overrightarrow{{\cal Q}^{n}})=0\wedge (%
\overrightarrow{{\cal Q}^{n}}.\overrightarrow{{\cal N}})=0,\qquad \forall 
\overrightarrow{{\cal Q}^{n}}\in \Omega ^{n+1},  \label{g2.17}
\end{equation}
$\overrightarrow{{\cal N}^{n}}$ is the null $n$th order multivector.
\end{opred}

Any $n$th order multivector $\overrightarrow{{\cal P}^{n}}$, having at least
two similar points $P_{k}=P_{l},\;\;(l\neq k)$, is the null $n$th order
multivector $\overrightarrow{{\cal N}^{n}}$, because in this case the
determinant in (\ref{g2.11}) has similar rows ($k$th and $l$th) and vanishes
for $\forall \overrightarrow{{\cal Q}^{n}}\in \Omega ^{n+1}$. The same is
valid for determinant $\left( \overrightarrow{{\cal Q}^{n}}.\overrightarrow{%
{\cal N}^{n}}\right) $, but in this case the determinant (\ref{g2.11}) has
similar columns ($k$th and $l$th) and vanishes for $\forall \overrightarrow{%
{\cal Q}^{n}}\in \Omega ^{n+1}$.

Any $n$th order multivector $\overrightarrow{{\cal P}^{n}}$ has either at
least two similar points and belongs to the equivalence class $\xi \left( 
{\cal N}^{n}\right) $ of null multivector $\overrightarrow{{\cal N}^{n}}$,
or has all different points and belongs to the equivalence class $\xi \left( 
\overrightarrow{M_{n}}({\cal P}^{n})\right) $, where $\overrightarrow{M_{n}}(%
{\cal P}^{n})=\left\{ \varepsilon ,\Sigma ,{\cal P}^{n}\right\} $ is the $n$%
th order oriented $\Sigma $-subspace of $V=\left\{ \Sigma ,\Omega \right\} $%
, consisting of $n+1$ points ${\cal P}^{n}$, and $\varepsilon $ is the
orientation of this subset which takes values $\varepsilon =\pm 1$. Thus,
the set $\Omega ^{n+1}$ of all $n$th order multivectors may be presented in
the form 
\begin{equation}
\Omega ^{n+1}=\xi \left( {\cal N}^{n}\right) \bigcup\limits_{\overrightarrow{%
M_{n}}({\cal P}^{n})}\xi \left( \overrightarrow{M_{n}}({\cal P}^{n})\right)
\label{g2.18}
\end{equation}
Any equivalence class $\xi \left( \overrightarrow{M_{n}}({\cal P}%
^{n})\right) $ contains $\frac{1}{2}(n+1)!$ elements. The equivalence class $%
\xi \left( +1,\Sigma ,{\cal P}^{n}\right) $ can be obtained from the
multivector $\overrightarrow{{\cal P}^{n}}$ by means of even permutation of
points. Multivectors which are obtained from the multivector $%
\overrightarrow{{\cal P}^{n}}$ by means of the odd permutation belong to the
equivalence class $\xi \left( -1,\Sigma ,{\cal P}^{n}\right) $. The factor
set $\Omega ^{n+1}/T$ is formed by the equivalence classes $\xi \left( {\cal %
N}^{n}\right) $ of null multivectors $\overrightarrow{{\cal N}^{n}}$ and $n$%
th order oriented $\Sigma $-subspaces $\overrightarrow{M_{n}}({\cal P}^{n})$%
. It is easy to see that the scalar $\Sigma $-product $(\overrightarrow{%
{\cal P}^{n}}.\overrightarrow{{\cal Q}^{n}})$ of two non-null $n$th order
multivectors $\overrightarrow{{\cal P}^{n}}\in \xi \left( \overrightarrow{%
M_{n}}({\cal P}^{n})\right) $ and $\overrightarrow{{\cal Q}^{n}}\in \xi
\left( \overrightarrow{M_{n}}({\cal Q}^{n})\right) $ depends only on $n$th
order oriented $\Sigma $-subspaces $\overrightarrow{M_{n}}({\cal P}^{n})$, $%
\overrightarrow{M_{n}}({\cal Q}^{n})$. Thus, in reality the scalar $\Sigma $%
-product $(\overrightarrow{{\cal P}^{n}}.\overrightarrow{{\cal Q}^{n}})$
describes mutual disposition of $n$th order oriented $\Sigma $-subspaces $%
\overrightarrow{M_{n}}({\cal P}^{n})$ and $\overrightarrow{M_{n}}({\cal Q}%
^{n}),$ which are constituents of the $\Sigma $-space $V=\left\{ \Sigma
,\Omega \right\} $. This result shows that as well as $\Sigma $-subspaces, 
{\it the multivectors are natural constituents of the }$\Sigma ${\it -space}%
, and description in terms of multivectors is equivalent to a description in
terms of $\Sigma $-subspaces.

This fact puts the question, if it is possible to express the scalar $\Sigma 
$-product of $\overrightarrow{M_{n}}({\cal P}^{n})$ and $\overrightarrow{%
M_{n}}({\cal Q}^{n})$ directly in terms of $n$th order oriented $\Sigma $%
-subspaces. In principle, this problem can be solved, but the solution will
be very complicated and ineffective. The fact is that the multivectors are
associated with integer numbers, whereas the oriented $\Sigma $-subspaces
are associated with natural numbers which do not contain zero. Descriptions
in terms of only natural numbers seems ineffective. Even a record of a
number in Roman numerals looks very complicated, because Roman numerals are
associated with natural numbers and do not contain zero. Apparently, the
rule for multiplication of two numbers, written in Roman numerals, looks
very complicated. Multivectors are associated with Arabic numerals. They
contain null multivectors. Reduced to $\Sigma $-subspaces, the multivectors
generate the equivalence class $\xi \left( +1,M_{n}\right) ,$ which is an
analog of positive number, the equivalence class $\xi \left( -1,M_{n}\right) 
$, which is an analog of negative numbers and the equivalence class $\xi
\left( {\cal N}^{n}\right) ,$ which is an analog of zero.

\begin{opred}
\label{d3.1.6c}. The length $|\overrightarrow{{\cal P}^{n}}|$ of the
multivector $\overrightarrow{{\cal P}^{n}}$ is the number 
\begin{equation}
|\overrightarrow{{\cal P}^{n}}|=\left\{ 
\begin{array}{c}
\mid \sqrt{(\overrightarrow{{\cal P}^{n}}.\overrightarrow{{\cal P}^{n}})}%
\mid =|\sqrt{F_{n}({\cal P}^{n})}|,\quad (\overrightarrow{{\cal P}^{n}}.%
\overrightarrow{{\cal P}^{n}})\geq 0 \\ 
i\mid \sqrt{(\overrightarrow{{\cal P}^{n}}.\overrightarrow{{\cal P}^{n}})}%
\mid =i|\sqrt{F_{n}({\cal P}^{n})}|,\quad (\overrightarrow{{\cal P}^{n}}.%
\overrightarrow{{\cal P}^{n}})<0
\end{array}
\right. \qquad \overrightarrow{{\cal P}^{n}}\in \Omega ^{n+1}  \label{g2.19}
\end{equation}
where the quantity $F_{n}({\cal P}^{n})$ is defined by the relations 
\begin{equation}
F_{n}:\quad \Omega ^{n+1}\rightarrow {\Bbb R},\qquad \Omega
^{n+1}=\bigotimes\limits_{k=1}^{n+1}\Omega ,\qquad n=1,2,\ldots
\label{g2.20}
\end{equation}
\begin{equation}
F_{n}\left( \overrightarrow{{\cal P}^{n}}\right) =\det ||\left( {\bf P}_{0}%
{\bf P}_{i}.{\bf P}_{0}{\bf P}_{k}\right) ||,\qquad P_{0},P_{i},P_{k}\in
\Omega ,\qquad i,k=1,2,...n  \label{g2.21}
\end{equation}
\begin{equation}
\left( {\bf P}_{0}{\bf P}_{i}.{\bf P}_{0}{\bf P}_{k}\right) \equiv \Sigma
\left( P_{i},P_{0}\right) +\Sigma \left( P_{0},P_{k}\right) -\Sigma \left(
P_{i},P_{k}\right) ,\qquad i,k=1,2,...n,  \label{a1.7}
\end{equation}
\end{opred}

The function (\ref{g2.20}) is a symmetric function of all its arguments $%
{\cal P}^{n}=\{ P_{0},P_{1},...,$ $P_{n}\} $, i.e. it is invariant with
respect to permutation of any points $P_{i},$ $P_{k}$, \ $i,k=0,1,...n$. It
follows from representation 
\[
F_{n}\left( \overrightarrow{{\cal P}^{n}}\right) =F_{n}\left( {\cal P}%
^{n}\right) =\left( \overrightarrow{{\cal P}^{n}}.\overrightarrow{{\cal P}%
^{n}}\right) 
\]
and the relation (\ref{g2.15}). It means that the squared length $|%
\overrightarrow{{\cal P}^{n}}|^{2}=\left| M\left( {\cal P}^{n}\right)
\right| ^{2}$ of any multivector $\overrightarrow{{\cal P}^{n}}$ does not
depend on the order of points. The squared length of any finite subset $%
{\cal P}^{n}$ is unique.

In the case, when multivector $\overrightarrow{{\cal P}^{n}}$ does not
contain similar points, it coincides with the oriented finite $\Sigma $%
-subspace $\overrightarrow{M_{n}({\cal P}^{n})}$, and it is a constituent of 
$\Sigma $-space. In the case, when at least two points of multivector
coincide, the multivector length vanishes, and the multivector is considered
to be a null multivector. The null multivector $\overrightarrow{{\cal P}^{n}}
$ is not a finite $\Sigma $-subspace $M_{n}({\cal P}^{n})$, or an oriented
finite $\Sigma $-subspace $\overrightarrow{M_{n}({\cal P}^{n})}$, but a use
of null multivectors assists in creation of a more simple technique, because
the null multivectors $\overrightarrow{{\cal P}^{n}}$ play a role of zeros.
Essentially, the multivectors are basic objects of T-geometry. As to
continual geometric objects, which are analogs of planes, sphere, ellipsoid,
etc., they are constructed by means of skeleton-envelope method (see \cite
{R02}) with multivectors, or finite $\Sigma $-subspaces used as skeletons.
As a consequence the T-geometry is presented $\sigma $-immanently, i.e.
without references to objects, external with respect to $\Sigma $-space.

The usual vector ${\bf P}_{0}{\bf P}_{1}\equiv \overrightarrow{P_{0}P_{1}}%
\equiv \overrightarrow{{\cal P}^{1}}=\left\{ P_{0},P_{1}\right\}
,\;\;P_{0},P_{1}\in \Omega $ is a special case of multivector. The squared
length $|{\bf P}_{0}{\bf P}_{1}|^{2}$ of the vector ${\bf P}_{0}{\bf P}_{1}$
is defined by the relation (\ref{a2.1a}). This gives 
\begin{equation}
|{\bf P}_{0}{\bf P}_{1}|^{2}\equiv \left( {\bf P}_{0}{\bf P}_{1}.{\bf P}_{0}%
{\bf P}_{1}\right) =\Sigma \left( P_{0},P_{1}\right) +\Sigma \left(
P_{1},P_{0}\right) =2G\left( P_{0},P_{1}\right)  \label{g2.22}
\end{equation}
The following quantities are also associated with the vector ${\bf P}_{0}%
{\bf P}_{1}$%
\begin{eqnarray}
|{\bf P}_{1}{\bf P}_{0}|^{2} &\equiv &\left( {\bf P}_{1}{\bf P}_{0}.{\bf P}%
_{1}{\bf P}_{0}\right) =2G\left( P_{1},P_{0}\right) ,  \label{g2.23} \\
\left( {\bf P}_{0}{\bf P}_{1}.{\bf P}_{1}{\bf P}_{0}\right) &=&-\Sigma
\left( P_{0},P_{1}\right) -\Sigma \left( P_{1},P_{0}\right) =-2G\left(
P_{0},P_{1}\right) ,  \label{g2.24}
\end{eqnarray}
It is rather unexpected that $|{\bf P}_{0}{\bf P}_{1}|^{2}=2G\left(
P_{0},P_{1}\right) $, but it is well that the vector ${\bf P}_{0}{\bf P}_{1}$
has only one length, but not two $\sqrt{2\Sigma \left( P_{0},P_{1}\right) }$
and $\sqrt{2\Sigma \left( P_{1},P_{0}\right) }$

\begin{opred}
\label{d1.13} The squared length $\left| M\left( {\cal P}^{n}\right) \right|
^{2}$ of the $n$th order $\Sigma $-subspace $M\left( {\cal P}^{n}\right)
\subset \Omega $ of the $\Sigma $-space $V=\left\{ \Sigma ,\Omega \right\} $
is the real number. 
\[
\left| M\left( {\cal P}^{n}\right) \right| ^{2}=F_{n}\left( {\cal P}%
^{n}\right) , 
\]
where $M\left( {\cal P}^{n}\right) =\left\{ P_{0},P_{1},...,P_{n},\right\}
\subset \Omega $ with all different $P_{i}\in \Omega $, \ $i=0,1,...n$,\ \ $%
\overrightarrow{{\cal P}^{n}}\in \Omega ^{n+1},$ and the quantity $F_{n}(%
{\cal P}^{n})$ is defined by the relations (\ref{g2.21}) -- (\ref{a1.7}).
\end{opred}

The meaning of the written relations is as follows. In the special case,
when the $\Sigma $-space is Euclidean space, its $\Sigma $-function is
symmetric and coincides with $\Sigma $-function of Euclidean space, any two
points $P_{0},P_{i}$ determine the vector ${\bf P}_{0}{\bf P}_{i}$, and the
relation (\ref{a1.7}) is a $\sigma $-immanent expression for the scalar $%
\Sigma $-product $\left( {\bf P}_{0}{\bf P}_{i}.{\bf P}_{0}{\bf P}%
_{k}\right) $ of two vectors. Then the relation (\ref{g2.21}) is the Gram's
determinant for $n$ vectors ${\bf P}_{0}{\bf P}_{i},\quad i=1,2,\ldots n$,
and $\sqrt{F_{n}({\cal P}^{n})}/n!$ is the Euclidean volume of the $(n+1)$%
-edr with vertices at the points ${\cal P}^{n}$.

Now we enable to formulate in terms of the world function the necessary and
sufficient condition of that the $\Sigma $-space is the $n$-dimensional
Euclidean space

\noindent I. 
\begin{equation}
\Sigma \left( P,Q\right) =\Sigma \left( Q,P\right) ,\qquad P,Q\in \Omega
\label{a3.4a}
\end{equation}

\noindent II. 
\begin{equation}
\exists {\cal P}^{n}\subset \Omega ,\qquad F_{n}({\cal P}^{n})\neq 0,\qquad
F_{n+1}(\Omega ^{n+2})=0,  \label{a3.4}
\end{equation}

\noindent III. 
\begin{equation}
\Sigma (P,Q)={\frac{1}{2}}\sum_{i,k=1}^{n}g^{ik}({\cal P}^{n})[x_{i}\left(
P\right) -x_{i}\left( Q\right) ][x_{k}\left( P\right) -x_{k}\left( Q\right)
],\qquad \forall P,Q\in \Omega ,  \label{a3.5}
\end{equation}
where the quantities $x_{i}\left( P\right) $, $x_{i}\left( Q\right) $ are
defined by the relations 
\begin{equation}
x_{i}\left( P\right) =\left( {\bf P}_{0}{\bf P}_{i}.{\bf P}_{0}{\bf P}%
\right) ,\qquad x_{i}\left( Q\right) =\left( {\bf P}_{0}{\bf P}_{i}.{\bf P}%
_{0}{\bf Q}\right) ,\qquad i=1,2,...n  \label{a3.5a}
\end{equation}
The contravariant components $g^{ik}({\cal P}^{n}),$ $(i,k=1,2,\ldots n)$ of
metric tensor are defined by its covariant components $g_{ik}({\cal P}^{n}),$
$(i,k=1,2,\ldots n)$ by means of relations 
\begin{equation}
\sum_{k=1}^{n}g_{ik}({\cal P}^{n})g^{kl}({\cal P}^{n})=\delta
_{i}^{l},\qquad i,l=1,2,\ldots n,  \label{a3.11}
\end{equation}
where covariant components $g_{ik}({\cal P}^{n})$ are defined by relations 
\begin{equation}
g_{ik}({\cal P}^{n})=\left( {\bf P}_{0}{\bf P}_{i}.{\bf P}_{0}{\bf P}%
_{k}\right) ,\qquad i,k=1,2,\ldots n  \label{a3.9}
\end{equation}

\noindent IV.\quad The relations 
\begin{equation}
\left( {\bf P}_{0}{\bf P}_{i}.{\bf P}_{0}{\bf P}\right) =x_{i},\qquad
x_{i}\in {\Bbb R},\qquad i=1,2,\ldots n,  \label{a3.12}
\end{equation}
considered to be equations for determination of $P\in \Omega $, have always
one and only one solution.

\noindent IVa. The relations (\ref{a3.12}), considered to be equations for
determination of $P\in \Omega $, have always not more than one solution.

\begin{remark}
\label{r3} The condition (\ref{a3.4}) is a corollary of the condition
(\ref{a3.5}). It is formulated in the form of a special condition, in order
that a determination of dimension were separated from determination of
coordinate system.
\end{remark}

The condition II determines the space dimension. The condition III describes 
$\sigma $-immanently the scalar $\Sigma $-product properties of the proper
Euclidean space. Setting $n+1$ points ${\cal P}^{n}$, satisfying the
condition II, one determines $n$-dimensional basis of vectors in Euclidean
space. Relations (\ref{a3.9}), (\ref{a3.11}) determine covariant and
contravariant components of the metric tensor, and the relations (\ref{a3.5a}%
) determine covariant coordinates of points $P$ and $Q$ at this basis. The
relation (\ref{a3.5}) determines the expression for $\Sigma $-function for
two arbitrary points in terms of coordinates of these points. Finally, the
condition IV describes continuity of the set $\Omega $ and a possibility of
the manifold construction on it. Necessity of conditions I -- IV for
Euclideaness of $\Sigma $-space is evident. One can prove their sufficiency 
\cite{R01}. The connection of conditions I -- IV with the Euclideaness of
the $\Sigma $-space can be formulated in the form of a theorem. 
\begin{theorem}
\label{c2}The $\Sigma  $-space $V=\{\Sigma  ,\Omega \}$ is the 
$n$-dimensional Euclidean  space, if and only if $\sigma $-immanent
conditions I -- IV are fulfilled.
\end{theorem}
\begin{remark}
\label{r2} For the $\sigma$-space were proper Euclidean, the eigenvalues
of the matrix $g_{ik}({\cal P}^n),\quad i,k=1,2,\ldots n$ must have the same
sign, otherwise it is pseudoeuclidean.
\end{remark}
The theorem states that it is sufficient to know metric (world function) to
construct the Euclidean geometry. Concepts of topological space and curve,
which are used usually in metric geometry for increasing its informativity,
appear to be excess in the sense that they are not needed for construction
of geometry.

Proof of this theorem can be found in \cite{R01}. A similar theorem for
another (but close) necessary and sufficient conditions has been proved in
ref. \cite{R90}. Here we show only constructive character of conditions I --
IV for proper Euclidean space. It means that starting from an abstract $%
\Sigma $-space, satisfying conditions I -- IV, one can determine dimension $%
n $ and construct a rectilinear coordinate system with conventional
description of the proper Euclidean space in it. One constructs sequentially
straight, two-dimensional plane, etc...up to $n$-dimensional plane coincide
with the set $\Omega $. To construct all these objects, one needs to develop
technique of T-geometry.

\begin{opred}
\label{d3.1.5c} Two $n$th order multivectors $\overrightarrow{{\cal P}^{n}}$%
, $\overrightarrow{{\cal Q}^{n}}$ are neutrally collinear ($n$-collinear) $%
\overrightarrow{{\cal P}^{n}}\parallel _{\left( {\rm n}\right) }%
\overrightarrow{{\cal Q}^{n}}$, if 
\begin{equation}
(\overrightarrow{{\cal P}^{n}}.\overrightarrow{{\cal Q}^{n}})(%
\overrightarrow{{\cal Q}^{n}}.\overrightarrow{{\cal P}^{n}})=|%
\overrightarrow{{\cal P}^{n}}|^{2}\cdot |\overrightarrow{{\cal Q}^{n}}|^{2}
\label{a2.10}
\end{equation}
\end{opred}

\begin{opred}
\label{d3.6} The $n$th order multivector $\overrightarrow{{\cal P}^{n}}$ is $%
f$-collinear to $n$th order multivector $\overrightarrow{{\cal Q}^{n}} $ $%
\left( \overrightarrow{{\cal P}^{n}}\parallel _{\left( {\rm f}\right) }%
\overrightarrow{{\cal Q}^{n}}\right) $, if 
\begin{equation}
(\overrightarrow{{\cal P}^{n}}.\overrightarrow{{\cal Q}^{n}})^{2}=|%
\overrightarrow{{\cal P}^{n}}|^{2}\cdot |\overrightarrow{{\cal Q}^{n}}|^{2}
\label{a2.10a}
\end{equation}
\end{opred}

\begin{opred}
\label{d3.7} The $n$th order multivector $\overrightarrow{{\cal P}^{n}}$ is $%
p$-collinear to $n$th order multivector $\overrightarrow{{\cal Q}^{n}}$ $%
\left( \overrightarrow{{\cal P}^{n}}\parallel _{\left( {\rm p}\right) }%
\overrightarrow{{\cal Q}^{n}}\right) $, if 
\begin{equation}
(\overrightarrow{{\cal Q}^{n}}.\overrightarrow{{\cal P}^{n}})^{2}=|%
\overrightarrow{{\cal P}^{n}}|^{2}\cdot |\overrightarrow{{\cal Q}^{n}}|^{2}
\label{a2.10b}
\end{equation}
\end{opred}

Here indices ''f'' and ''p'' are associated with the terms ''future'' and
''past'' respectively.

In the symmetric T-geometry there is only one type of collinearity, because
the three mentioned types of collinearity coincide in the symmetric
T-geometry. The property of the neutral collinearity is commutative, i.e. if 
$\overrightarrow{{\cal P}^{n}}\parallel _{\left( {\rm n}\right) }%
\overrightarrow{{\cal Q}^{n}}$, then $\overrightarrow{{\cal Q}^{n}}\parallel
_{\left( {\rm n}\right) }\overrightarrow{{\cal P}^{n}}$. The property of $p$%
-collinearity and $f$-collinearity are not commutative, in general. Instead,
one has according to (\ref{a2.10a}) and (\ref{a2.10b}) that, if $%
\overrightarrow{{\cal P}^{n}}\parallel _{\left( {\rm p}\right) }%
\overrightarrow{{\cal Q}^{n}}$, then $\overrightarrow{{\cal Q}^{n}}\parallel
_{\left( {\rm f}\right) }\overrightarrow{{\cal P}^{n}}$.

\begin{opred}
\label{d3.15e}. The $n$th order multivector $\overrightarrow{{\cal P}^{n}} $
is $f$-parallel to the $n$th order multivector $\overrightarrow{{\cal Q}^{n}}
$ $\left( \overrightarrow{{\cal P}^{n}}\uparrow \uparrow _{({\rm f})}%
\overrightarrow{{\cal Q}^{n}}\right) $, if 
\begin{equation}
(\overrightarrow{{\cal P}^{n}}.\overrightarrow{{\cal Q}^{n}})=|%
\overrightarrow{{\cal P}^{n}}|\cdot |\overrightarrow{{\cal Q}^{n}}|
\label{a2.11}
\end{equation}
The $n$th order multivector $\overrightarrow{{\cal P}^{n}}$ is $f$%
-antiparallel to the $n$th order multivector $\overrightarrow{{\cal Q}^{n}}$ 
$\left( \overrightarrow{{\cal P}^{n}}\uparrow \downarrow _{({\rm f})}%
\overrightarrow{{\cal Q}^{n}}\right) $, if 
\begin{equation}
(\overrightarrow{{\cal P}^{n}}.\overrightarrow{{\cal Q}^{n}})=-|%
\overrightarrow{{\cal P}^{n}}|\cdot |\overrightarrow{{\cal Q}^{n}}|
\label{a2.12}
\end{equation}
\end{opred}

\begin{opred}
\label{d3.15ee}The $n$th order multivector $\overrightarrow{{\cal P}^{n}}$
is $p$-parallel to the $n$th order multivector $\overrightarrow{{\cal Q}^{n}}
$ $\left( \overrightarrow{{\cal P}^{n}}\uparrow \uparrow _{({\rm p})}%
\overrightarrow{{\cal Q}^{n}}\right) $, if 
\begin{equation}
(\overrightarrow{{\cal Q}^{n}}.\overrightarrow{{\cal P}^{n}})=|%
\overrightarrow{{\cal P}^{n}}|\cdot |\overrightarrow{{\cal Q}^{n}}|
\label{a2.12a}
\end{equation}
The $n$th order multivector $\overrightarrow{{\cal P}^{n}}$ is $p$%
-antiparallel to the $n$th order multivector $\overrightarrow{{\cal Q}^{n}}$ 
$\left( \overrightarrow{{\cal P}^{n}}\uparrow \uparrow _{({\rm p})}%
\overrightarrow{{\cal Q}^{n}}\right) $, if 
\begin{equation}
(\overrightarrow{{\cal Q}^{n}}.\overrightarrow{{\cal P}^{n}})=-|%
\overrightarrow{{\cal P}^{n}}|\cdot |\overrightarrow{{\cal Q}^{n}}|
\label{a2.12b}
\end{equation}
\end{opred}

The $f$-parallelism and the $p$-parallelism are connected as follows. If $%
\overrightarrow{{\cal P}^{n}}\uparrow \uparrow _{({\rm p})}\overrightarrow{%
{\cal Q}^{n}}$, then $\overrightarrow{{\cal Q}^{n}}\uparrow \uparrow _{({\rm %
f})}\overrightarrow{{\cal P}^{n}}$ and vice versa.

Vector ${\bf P}_{0}{\bf P}_{1}=\overrightarrow{{\cal P}^{1}}$ as well as the
vector ${\bf Q}_{0}{\bf Q}_{1}=\overrightarrow{{\cal Q}^{1}}$ are the first
order multivectors. If ${\bf P}_{0}{\bf P}_{1}\uparrow \uparrow _{({\rm f})}%
{\bf Q}_{0}{\bf Q}_{1}$, then ${\bf P}_{1}{\bf P}_{0}\uparrow \downarrow _{(%
{\rm f})}{\bf Q}_{0}{\bf Q}_{1}$ and ${\bf P}_{0}{\bf P}_{1}\uparrow
\downarrow _{({\rm f})}{\bf Q}_{1}{\bf Q}_{0}$

\section{Geometrical objects in $\Sigma $-space. Tubes and their properties.}

The simplest geometrical object in T-geometry is the $n$th order tube ${\cal %
T}\left( {\cal P}^{n}\right) $, which is determined by its skeleton ${\cal P}%
^{n}$. The tube is an analog of Euclidean $n$-dimensional plane, which is
also determined by $n+1$ points ${\cal P}^{n},$ not belonging to a $\left(
n-1\right) $-dimensional plane.

\begin{opred}
\label{d1.14} $n$th order $\Sigma $-subspace $M\left( {\cal P}^{n}\right) =%
{\cal P}^{n}$\ of nonzero length $\;\left| M\left( {\cal P}^{n}\right)
\right| ^{2}=\left| {\cal P}^{n}\right| ^{2}=F_{n}\left( {\cal P}^{n}\right)
\neq 0$ determines geometrical object (set of points) ${\cal T}$ $\left( 
{\cal P}^{n}\right) $, called the $n$th order tube, by means of relation 
\begin{equation}
{\cal T}\left( {\cal P}^{n}\right) \equiv {\cal T}_{{\cal P}^{n}}=\left\{
P_{n+1}|F_{n+1}\left( {\cal P}^{n+1}\right) =0\right\} ,\qquad P_{i}\in
\Omega ,\qquad i=0,1\ldots n+1,  \label{b1.3}
\end{equation}
where the function $F_{n}$ is defined by the relations (\ref{g2.20}) -- (\ref
{a1.7})
\end{opred}

The shape of the tube ${\cal T}\left( {\cal P}^{n}\right) $ does not depend
on the order of points of the multivector $\overrightarrow{{\cal P}^{n}}$.
The basic point ${\cal P}^{n}$, determining the tube ${\cal T}_{{\cal P}%
^{n}} $ belong to ${\cal T}_{{\cal P}^{n}}$.

The first order tube ${\cal T}_{P_{0}P_{1}}$ can be defined by means of
concept of $n$-collinearity (\ref{a2.10}) 
\begin{eqnarray}
{\cal T}\left( {\cal P}^{1}\right) &\equiv &{\cal T}_{({\rm n}%
)P_{0}P_{1}}=\left\{ P_{2}|F_{2}\left( {\cal P}^{2}\right) =0\right\} \equiv
\left\{ P_{2}\left| \overrightarrow{P_{0}P_{1}}||_{({\rm n)}}\overrightarrow{%
P_{0}P_{2}}\right. \right\}  \nonumber \\
&\equiv &\left\{ P_{2}\left| \;|\overrightarrow{P_{0}P_{1}}|^{2}|%
\overrightarrow{P_{0}P_{2}}|^{2}-\left( \overrightarrow{P_{0}P_{1}}.%
\overrightarrow{P_{0}P_{2}}\right) \left( \overrightarrow{P_{0}P_{2}}.%
\overrightarrow{P_{0}P_{1}}\right) =0\right. \right\}  \label{f2.3}
\end{eqnarray}
As far as there are concepts of $f$-collinearity and of $p$-collinearity,
one can define also the first order $f$-tube and $p$-tube on the basis of
these collinearities. The first order $f$-tube is defined by the relation 
\begin{equation}
{\cal T}_{({\rm f})P_{0}P_{1}}=\left\{ R\left| \overrightarrow{P_{0}P_{1}}%
||_{({\rm f})}\overrightarrow{P_{0}R}\right. \right\} =\left\{ R\left| \;|%
\overrightarrow{P_{0}P_{1}}|^{2}|\overrightarrow{P_{0}R}|^{2}-\left( 
\overrightarrow{P_{0}P_{1}}.\overrightarrow{P_{0}R}\right) ^{2}=0\right.
\right\}  \label{f2.4}
\end{equation}
The first order $p$-tube is defined as follows 
\begin{equation}
{\cal T}_{({\rm p})P_{0}P_{1}}=\left\{ R\left| \overrightarrow{P_{0}P_{1}}%
||_{({\rm p})}\overrightarrow{P_{0}R}\right. \right\} =\left\{ R\left| \;|%
\overrightarrow{P_{0}P_{1}}|^{2}|\overrightarrow{P_{0}R}|^{2}-\left( 
\overrightarrow{P_{0}R}.\overrightarrow{P_{0}P_{1}}\right) ^{2}=0\right.
\right\}  \label{f2.5}
\end{equation}

In the symmetric T-geometry all three tubes (\ref{f2.3}) -- (\ref{f2.5})
coincide. In the nonsymmetric T-geometry they are different, in general. The
tubes (\ref{f2.3}), (\ref{f2.4}), (\ref{f2.5}) can be divided into segments,
each of them is determined by one of factors of expressions (\ref{f2.3}) -- (%
\ref{f2.5}).

In all cases the factorization of the expressions 
\[
F_{({\rm f})}\left( P_{0},P_{1},P_{2}\right) =\left| \overrightarrow{%
P_{0}P_{1}}\right| ^{2}\left| \overrightarrow{P_{0}P_{2}}\right| ^{2}-\left( 
\overrightarrow{P_{0}P_{1}}.\overrightarrow{P_{0}P_{2}}\right) ^{2}
\]
\[
F_{({\rm p})}\left( P_{0},P_{1},P_{2}\right) =\left| \overrightarrow{%
P_{0}P_{1}}\right| ^{2}\left| \overrightarrow{P_{0}P_{2}}\right| ^{2}-\left( 
\overrightarrow{P_{0}P_{2}}.\overrightarrow{P_{0}P_{1}}\right) ^{2},
\]
\begin{equation}
F_{({\rm n})}\left( P_{0},P_{1},P_{2}\right) =\left| \overrightarrow{%
P_{0}P_{1}}\right| ^{2}\left| \overrightarrow{P_{0}P_{2}}\right| ^{2}-\left( 
\overrightarrow{P_{0}P_{1}}.\overrightarrow{P_{0}P_{2}}\right) \left( 
\overrightarrow{P_{0}P_{2}}.\overrightarrow{P_{0}P_{1}}\right) 
\label{f2.2b}
\end{equation}
have similar form 
\begin{equation}
F_{({\rm q})}\left( P_{0},P_{1},P_{2}\right) =-F_{({\rm q}0)}F_{({\rm q}%
1)}F_{({\rm q}2)}F_{({\rm q}3)}  \label{f2.0}
\end{equation}
Here index $q$ runs values $f,p,n$, and factorization of expressions $F_{(%
{\rm q})}\left( P_{0},P_{1},P_{2}\right) ,$ $q=f,p,n$ has a similar form 
\[
F_{({\rm q}0)}=F_{({\rm q}0)}\left( P_{0},P_{1},P_{2}\right) =\sqrt{G_{02}}+%
\sqrt{G_{10}}+\sqrt{G_{12}-\eta _{{\rm q}}}
\]
\begin{eqnarray}
F_{({\rm q}1)} &=&F_{({\rm q}1)}\left( P_{0},P_{1},P_{2}\right) =\sqrt{G_{02}%
}-\sqrt{G_{10}}+\sqrt{G_{12}-\alpha _{{\rm q}}\eta _{{\rm q}}}  \nonumber \\
F_{({\rm q}2)} &=&F_{({\rm q}2)}\left( P_{0},P_{1},P_{2}\right) =\sqrt{G_{02}%
}+\sqrt{G_{10}}-\sqrt{G_{12}-\eta _{{\rm q}}}  \label{f2.1} \\
F_{({\rm q}3)} &=&F_{({\rm q}3)}\left( P_{0},P_{1},P_{2}\right) =\sqrt{G_{02}%
}-\sqrt{G_{10}}-\sqrt{G_{12}-\alpha _{{\rm q}}\eta _{{\rm q}}}  \nonumber
\end{eqnarray}
where for brevity one uses designations 
\[
G_{ik}=G\left( P_{i},P_{k}\right) ,\qquad A_{ik}=A\left( P_{i},P_{k}\right)
,\qquad i,k=0,1,2
\]
\begin{eqnarray}
\eta _{{\rm f}} &=&-\eta _{{\rm p}}=A_{10}+A_{02}+A_{21},\qquad \eta _{{\rm n%
}}=\frac{\eta _{{\rm f}}^{2}}{\sqrt{4G_{01}G_{02}+\eta _{{\rm f}}^{2}}+2%
\sqrt{G_{01}G_{02}}}  \label{f2.2ee} \\
\alpha _{{\rm p}} &=&\alpha _{{\rm f}}=1,\qquad \alpha _{{\rm n}}=-1 
\nonumber
\end{eqnarray}
In the symmetric geometry, when $A\left( P,Q\right) =0$,\ \ $\forall P,Q\in
\Omega $, and $\eta =0$, all expressions (\ref{f2.1}), for $F_{\left( {\rm n}%
i\right) },F_{\left( {\rm f}i\right) },F_{\left( {\rm p}i\right)
},\;\;i=0,1,2,3$ coincide.

Factorizations (\ref{f2.0}), (\ref{f2.1}) determine division of the tubes
into segments. For instance, segments ${\cal T}_{({\rm f})[P_{0}P_{1}]},%
{\cal T}_{({\rm p})[P_{0}P_{1}]},{\cal T}_{({\rm n})[P_{0}P_{1}]}$ of the
tubes (\ref{f2.4}), (\ref{f2.5}), (\ref{f2.3}) between the points $P_{0}$
and $P_{1}$ are determined respectively by the relations 
\begin{equation}
{\cal T}_{({\rm q})[P_{0}P_{1}]}=\left\{ P_{2}\left| \sqrt{G_{02}}-\sqrt{%
G_{10}}+\sqrt{G_{12}-\alpha _{{\rm q}}\eta _{{\rm q}}}=0\right. \right\} 
\label{f2.6}
\end{equation}
where index $q$ runs values $f,p,n$. Values of $\eta _{{\rm q}}$, $\alpha _{%
{\rm q}}$ are determined by the relations (\ref{f2.2ee}).

\begin{opred}
\label{d3.1.9}. Section ${\cal S}_{n;P}$ of the tube ${\cal T}({\cal P}^{n})$
at the point $P\in {\cal T}({\cal P}^{n})$ is the set ${\cal S}_{n;P}({\cal T%
}({\cal P}^{n}))$ of points, belonging to the tube ${\cal T}({\cal P}^{n})$ 
\begin{equation}
{\cal S}_{n;P}({\cal T}({\cal P}^{n}))=\{P^{\prime }\mid
\bigwedge_{l=0}^{l=n}\Sigma (P_{l},P^{\prime })=\Sigma (P_{l},P)\},\qquad
P\in {\cal T}({\cal P}^{n}),\qquad P^{\prime }\in \Omega .  \label{a2.38}
\end{equation}
\end{opred}

Let us note that ${\cal S}_{n;P}({\cal T}({\cal P}^{n}))\subset {\cal T}(%
{\cal P}^{n})$, because $P\in {\cal T}({\cal P}^{n})$. Indeed, whether the
point $P$ belongs to ${\cal T}({\cal P}^{n})$ depends only on values of $n+1$
quantities $\Sigma (P_{l},P),\;\;l=0,1,...n$. In accordance with (\ref{a2.38}%
) these quantities are the same for both points $P$ and $P^{\prime }$.
Hence, the running point $P^{\prime }\in {\cal T}({\cal P}^{n})$, if $P\in 
{\cal T}({\cal P}^{n})$.

In the proper Euclidean space the $n$th order tube is $n$-dimensional plane,
containing points ${\cal P}^{n}$, and its section ${\cal S}_{n;P}({\cal T}(%
{\cal P}^{n}))$ at the point $P$ consists of one point $P$.

\begin{opred}
\label{d3.10}Section ${\cal S}_{n;P}$ of the tube ${\cal T}({\cal P}^{n})$
at the point $P\in {\cal T}({\cal P}^{n})$ is minimal, if ${\cal S}%
_{n;P}\left( {\cal T}({\cal P}^{n})\right) =\left\{ P\right\} $.
\end{opred}

\begin{opred}
\label{d3.11}The first order tube ${\cal T}_{P_{0}P_{1}}$ is degenerate, if
its section at any point $P\in {\cal T}_{P_{0}P_{1}}$ is minimal.
\end{opred}

Minimality of the first tube section means that the first order tube
degenerates to a curve, and any section of the tube consists of one point.
It means that there is only one vector $\overrightarrow{P_{0}R}$, $R\in 
{\cal T}_{P_{0}P_{1}}$ of fixed length, which is parallel, or antiparallel
to the vector $\overrightarrow{P_{0}P_{1}}$. As far as in the nonsymmetric
T-geometry there is several types of parallelism, there is several types of
degeneration, in general. In the symmetric T-geometry there is only one type
of the degeneration.

\begin{opred}
\label{d3.12}The $\Sigma $-space $V=\left\{ \Sigma ,\Omega \right\} $ is
degenerate on the set ${\cal T}$ of the first order tubes, if the set ${\cal %
T}$ contains only degenerate tubes ${\cal T}_{P_{0}P_{1}}$.
\end{opred}

\begin{opred}
\label{d3.12a}The $\Sigma $-space $V=\left\{ \Sigma ,\Omega \right\} $ is
locally $f$-degenerate, if all first order tubes ${\cal T}_{({\rm f}%
)P_{0}P_{1}}$ are degenerate.
\end{opred}

\begin{opred}
\label{d3.13}The $\Sigma $-space $V=\left\{ \Sigma ,\Omega \right\} $ is
locally $p$-degenerate, if all first order tubes ${\cal T}_{({\rm p}%
)P_{0}P_{1}}$ are degenerate.
\end{opred}

\begin{opred}
\label{d3.14}The $\Sigma $-space $V=\left\{ \Sigma ,\Omega \right\} $ is
locally $n$-degenerate, if all first order tubes ${\cal T}_{({\rm n}%
)P_{0}P_{1}}\equiv {\cal T}_{P_{0}P_{1}}$ are degenerate.
\end{opred}

Note that the Riemannian space considered to be a $\Sigma $-space is locally
degenerate.

\section{ Asymmetric T-geometry on manifold}

We have considered T-geometry in the coordinate-free form. But to discover a
connection between the T-geometry and usual differential geometry, one needs
to introduce coordinates and to consider the T-geometry on a manifold. It is
important also from the viewpoint of the asymmetric T-geometry application
as a possible space-time geometry. The asymmetric T-geometry on the manifold
may be considered to be a conventional symmetric geometry (for instance,
Riemannian) with additional force fields $a_{i}\left( x\right) $, $%
a_{ikl}\left( x\right) $, generated on the manifold by the antisymmetric
component $A$ of the world function. Testing experimentally existence of
these force fields, one can conclude whether the antisymmetric component $A$
exists and how large it is.

Let us suppose that there exist a set of one-to-one mappings 
\begin{equation}
\mu _{P}:\qquad V\rightarrow D_{P,n}\subset E_{P,n},\qquad \forall P\in
\Omega  \label{aa2.1}
\end{equation}
where $V=\left\{ \Sigma ,\Omega \right\} $ is the $\Sigma $-space and $%
D_{P,n}=\left\{ G_{\left( {\rm E}\right) },{\cal M}_{P,n}\right\} $ is a
continuous region of the $n$-dimensional Euclidean space $E_{P,n}=\left\{
G_{\left( {\rm E}\right) },{\Bbb R}^{n}\right\} ,$ ${\cal M}_{P,n}\subset 
{\Bbb R}^{n}$. The world function of $D_{P,n}$ has the form 
\begin{equation}
G_{\left( {\rm E}\right) }\left( x,x^{\prime }\right) =\frac{1}{2}%
g_{ik}\left( x^{i}-x^{\prime i}\right) \left( x^{k}-x^{\prime k}\right)
,\qquad x^{i},x^{\prime i}\in {\Bbb R},\qquad i,k=1,2,...n  \label{aa2.2}
\end{equation}
Here $g_{ik}=$const, $\det ||g_{ik}||\neq 0$ is the metric tensor in the
rectilinear coordinate system $K$, introduced in the Euclidean space $%
E_{P,n} $. It is easy to verify, that 
\begin{equation}
G_{\left( {\rm E}\right) ,i}g^{ik}G_{\left( {\rm E}\right) ,k}=2G_{\left( 
{\rm E}\right) },\qquad G_{\left( {\rm E}\right) ,i}G_{\left( {\rm E}\right)
}^{ik^{\prime }}G_{\left( {\rm E}\right) ,k^{\prime }}=2G_{\left( {\rm E}%
\right) }  \label{aa2.3}
\end{equation}
where the following designations are used 
\begin{equation}
G_{\left( {\rm E}\right) ,i}\equiv \partial _{i}G_{\left( {\rm E}\right)
}\equiv \frac{\partial G_{\left( {\rm E}\right) }}{\partial x^{i}},\qquad
G_{\left( {\rm E}\right) ,i^{\prime }}\equiv \partial _{i^{\prime
}}G_{\left( {\rm E}\right) }\equiv \frac{\partial G_{\left( {\rm E}\right) }%
}{\partial x^{\prime i}}  \label{aa2.4}
\end{equation}
\begin{equation}
G_{\left( {\rm E}\right) ik^{\prime }}\equiv G_{\left( {\rm E}\right)
,ik^{\prime }}\equiv \partial _{i}\partial _{k^{\prime }}G_{\left( {\rm E}%
\right) }\equiv \frac{\partial ^{2}G_{\left( {\rm E}\right) }}{\partial
x^{i}\partial x^{\prime k}},\qquad G_{\left( {\rm E}\right) ik^{\prime
}}G_{\left( {\rm E}\right) }^{lk^{\prime }}=\delta _{k}^{l}  \label{aa2.5}
\end{equation}
\begin{equation}
g_{ik}=-\left[ G_{\left( {\rm E}\right) ik^{\prime }}\right] _{x^{\prime
}=x},\qquad g^{ik}=-\left[ G_{\left( {\rm E}\right) }^{ik^{\prime }}\right]
_{x^{\prime }=x},  \label{aa2.6}
\end{equation}

As far as all quantities in the relations (\ref{aa2.3}) are tensors with
respect to arbitrary coordinate transformations, the relations (\ref{aa2.3})
are valid in any curvilinear coordinate system.

As far as the mapping $\mu _{P}$ is one-to-one, one can use coordinates $%
x=\left\{ x^{1},x^{2},...,x^{n}\right\} $ for labelling of points in $V$.
Let $x=\mu _{P}\left( P\right) \in {\cal M}_{P,n}$ and $x^{\prime }=\mu
_{P}\left( P^{\prime }\right) \in {\cal M}_{P,n}$ be images of points $P\in
\Omega $ and $P^{\prime }\in \Omega $ respectively at the mapping $\mu _{P}$%
. At this mapping the symmetric component $G$ of the world function $\Sigma $
transforms as follows $G\left( P,P^{\prime }\right) \rightarrow G_{P}\left(
x,x^{\prime }\right) =G\left( \mu _{P}\left( P\right) ,\mu _{P}\left(
P^{\prime }\right) \right) ,$ and relations (\ref{aa2.3}) do not take place
for $G_{P}\left( x,x^{\prime }\right) $, in general. But if 
\begin{equation}
\lim\Sb x^{\prime }\rightarrow x,  \\ G_{P}\left( x,x^{\prime }\right) >0 
\endSb \frac{G_{P,i}\left( x,x^{\prime }\right) g^{ik}\left( x\right)
G_{P,k}\left( x,x^{\prime }\right) -2G_{P}\left( x,x^{\prime }\right) }{%
2G_{P}\left( x,x^{\prime }\right) }=0,  \label{aa2.7}
\end{equation}
one may speak that in vicinity of the point $x=\mu _{P}\left( P\right) $ the 
$\Sigma $-space $\left\{ G_{P},{\cal M}_{P,n}\right\} $ has a structure
close to that of the Euclidean space. If the relation (\ref{aa2.7}) takes
place for all point $P\in \Omega $, one may speak that the local structure
of the $\Sigma $-space $\left\{ G_{P},{\cal M}_{P,n}\right\} $ is close to
that of the Euclidean space.

Let the vector ${\bf u}_{xx^{\prime }}=\left\{ x,x^{\prime }\right\} =\mu
_{P}\left( \overrightarrow{PP}\right) $ in ${\cal M}_{P,n}$ be the image of
the vector $\overrightarrow{PP^{\prime }}$ in $\Omega $. Then $%
-G_{P,i}\equiv -G_{P,i}\left( x,x^{\prime }\right) $, $i=1,2,...n$ are
covariant coordinates of the vector ${\bf u}_{xx^{\prime }}$ at the
coordinate system $K.$ In the Euclidean space $E_{\left( {\rm P}\right) }$
one can introduce a linear space ${\cal L}_{\left( {\rm P}\right) }$ and
consider the vector ${\bf u}_{xx^{\prime }}$ to be a vector of ${\cal L}%
_{\left( {\rm P}\right) }$ with coordinates $-G_{P,i}$ at the coordinate
system $K$. The squared length $\left| {\bf u}_{xx^{\prime }}\right|
_{\left( {\rm E}\right) }^{2}$ of ${\bf u}_{xx^{\prime }}$ as a vector of $%
{\cal L}_{\left( {\rm P}\right) }$ is 
\begin{equation}
\left| {\bf u}_{xx^{\prime }}\right| _{\left( {\rm E}\right)
}^{2}=G_{P,i}g^{ik}G_{P,k}  \label{aa2.10}
\end{equation}
whereas the squared length $\left| {\bf PP}^{\prime }\right| ^{2}=\left( 
{\bf PP}^{\prime }.{\bf PP}^{\prime }\right) $ as a length of an original of 
${\bf u}_{xx^{\prime }}$ in $V=\left\{ \Sigma ,\Omega \right\} $ is equal to 
$2G\left( P,P^{\prime }\right) $. The numerator in relation (\ref{aa2.7})
describes this difference.

In the special case, when $\Sigma $-space $V=\left\{ \Sigma ,\Omega \right\} 
$ is a Riemannian space, the condition 
\begin{equation}
\left| {\bf u}_{xx^{\prime }}\right| _{\left( {\rm E}\right)
}^{2}=G_{P,i}g^{ik}G_{P,k}=2G=\left| {\bf PP}^{\prime }\right| ^{2}
\label{aa2.10a}
\end{equation}
takes place \cite{S60,R62}. In this case the Euclidean space $%
D_{P,n}=\left\{ G_{\left( {\rm E}\right) },{\cal M}_{P,n}\right\} $ is the
Euclidean space tangent to $V=\left\{ \Sigma ,\Omega \right\} $ at the point 
$P$. The mapping (\ref{aa2.1}) is a geodesic mapping, which transforms any
vector $\overrightarrow{PP^{\prime }}$ at the point $P$ in the Riemannian
space $V$ to the tangent vector ${\bf u}_{xx^{\prime }}${\bf \ }of the same
length in the tangent Euclidean space $D_{P,n}$. In this special case the
quantities $-G_{P,i}\equiv -G_{P,i}\left( x,x^{\prime }\right) $, $%
i=1,2,...n $ may be considered to be coordinates of the vector$\ 
\overrightarrow{PP^{\prime }}$ in the coordinate system $K$. Numerator of
the expression (\ref{aa2.7}) vanishes identically. It means that the local
structure of the Riemannian space is close to that of the Euclidean space.

Any mapping (\ref{aa2.1}) determines the world function $\Sigma \left(
P,P^{\prime }\right) $ as a function 
\[
\Sigma _{\left( x\right) }\left( x,x^{\prime }\right) =\Sigma \left( \mu
_{x}^{-1}\left( x\right) ,\mu _{x}^{-1}\left( x^{\prime }\right) \right) 
\]
on the manifold ${\cal M}_{P,n}\times {\cal M}_{P,n}$, where $\mu
_{x}^{-1}\left( x\right) $ is the function reverse to $x=\mu _{P}\left(
P\right) $. If the local structure of the function $\Sigma _{\left( x\right)
}\left( x,x^{\prime }\right) $ given on the manifold ${\cal M}_{P,n}$ is
close to the Euclidean structure, one can differentiate the function $\Sigma
_{\left( x\right) }\left( x,x^{\prime }\right) $ with respect to $x$ and $%
x^{\prime }$.

One may apply another approach which appears to be more convenient. Let it
be possible to attribute $n+1$ real numbers $x=\left\{ x^{i}\right\} ,$ $%
\;i=0,1,...n$ to any point $P$ in such a way, that there be one-to-one
correspondence between the point $P$ and the set $x$ of $n+1$ coordinates $%
\left\{ x^{i}\right\} ,$ $\;i=0,1,...n$. All points $x$ form a set ${\cal M}%
_{n+1}$. Then the world function $\Sigma \left( P,P^{\prime }\right) $ is a
function $\Sigma \left( x,x^{\prime }\right) $%
\begin{equation}
\Sigma :\;\;\;{\cal M}_{n+1}\times {\cal M}_{n+1}\rightarrow {\Bbb R},\qquad
\Sigma \left( x,x\right) =0,\qquad \forall x\in {\cal M}_{n+1}  \label{c3.1}
\end{equation}
of coordinates $x,x^{\prime }\in {\cal M}_{n+1}\subset {\Bbb R}^{n+1}$ of
points $P,P^{\prime }\in \Omega $. Two-point quantities ($\Sigma $-function
and their derivatives) are designed as a rule by capital characters.
One-point quantities are designed by small characters.

Let the function $\Sigma \left( x,x^{\prime }\right) $ be multiply
differentiable. Then the set ${\cal M}_{n+1}\subset {\Bbb R}^{n+1}$ may be
called the $\left( n+1\right) $th order manifold. One can differentiate $%
\Sigma \left( x,x^{\prime }\right) $ with respect to $x^{i}$ and with
respect to $x^{\prime i}$, forming two-point tensors. For instance, 
\begin{eqnarray*}
\Sigma _{,k}\left( x,x^{\prime }\right) &\equiv &\frac{\partial }{\partial
x^{k}}\Sigma \left( x,x^{\prime }\right) ,\qquad \Sigma _{,k^{\prime
}}\left( x,x^{\prime }\right) \equiv \frac{\partial }{\partial x^{\prime k}}%
\Sigma \left( x,x^{\prime }\right) \\
\Sigma _{,kl^{\prime }}\left( x,x^{\prime }\right) &\equiv &\Sigma
_{,l^{\prime }k}\left( x,x^{\prime }\right) \equiv \frac{\partial ^{2}}{%
\partial x^{k}\partial x^{\prime l}}\Sigma \left( x,x^{\prime }\right) ,
\end{eqnarray*}
are two-point tensors. Here indices after comma mean differentiation with
respect to $x^{k}$, if the index $k$ has not a prime and differentiation
with respect to $x^{\prime k}$, if the index $k$ has a prime. The first
argument of the two-point quantity is denoted by unprimed variable, whereas
the second one is denoted by primed one. Primed indices relate to the second
argument of the two-point quantity, whereas the unprimed ones relate to the
first argument.

$\Sigma _{k}\equiv \Sigma _{,k}=\Sigma _{,k}\left( x,x^{\prime }\right) $ is
a vector at the point $x$ and a scalar at the point $x^{\prime }.$ Vice
versa $\Sigma _{k^{\prime }}\equiv \Sigma _{,k^{\prime }}=\Sigma
_{,k^{\prime }}\left( x,x^{\prime }\right) $ is a vector at the point $%
x^{\prime }$ and a scalar at the point $x$. The quantity $\Sigma
_{,kl^{\prime }}=\Sigma _{,kl^{\prime }}\left( x,x^{\prime }\right) $ is a
vector at the point $x$ and a vector at the point $x^{\prime }$. Other
derivatives are not tensors. For instance, $\Sigma _{,kl}\left( x,x^{\prime
}\right) \equiv \Sigma _{,lk}\left( x,x^{\prime }\right) \equiv \frac{%
\partial ^{2}}{\partial x^{k}\partial x^{l}}\Sigma \left( x,x^{\prime
}\right) $ is a scalar at the point $x^{\prime }$, but it is not a tensor at
the point $x.$

To construct tensors of higher rank by means of differentiation, let us
introduce covariant derivatives. Let $\Sigma _{,kl^{\prime }}\equiv \Sigma
_{kl^{\prime }}\equiv \Sigma _{l^{\prime }k}$ and $\det ||\Sigma
_{kl^{\prime }}||\neq 0$. The quantity $\Sigma _{kl^{\prime }}$ will be
referred to as covariant fundamental metric tensor. One can introduce also
contravariant fundamental metric tensor $\Sigma ^{ik^{\prime }}\equiv \Sigma
^{k^{\prime }i}$, defining it by the relation 
\begin{equation}
\Sigma ^{ik^{\prime }}\Sigma _{lk^{\prime }}=\delta _{l}^{i}\qquad \Sigma
^{i^{\prime }k}\Sigma _{l^{\prime }k}=\delta _{l^{\prime }}^{i^{\prime }}
\label{c3.4}
\end{equation}
Let us note that the quantity 
\begin{equation}
\tilde{\Gamma}_{kl}^{i}\left( x,x^{\prime }\right) \equiv \Sigma
^{is^{\prime }}\Sigma _{,kls^{\prime }},\qquad \Sigma _{,kls^{\prime
}}\equiv \frac{\partial ^{3}\Sigma }{\partial x^{k}\partial x^{l}\partial
x^{\prime s}}  \label{c3.6}
\end{equation}
is a scalar at the point $x^{\prime }$ and the Christoffel symbol at the
point $x$. Vice versa, the quantity 
\begin{equation}
\tilde{\Gamma}_{k^{\prime }l^{\prime }}^{i^{\prime }}\left( x,x^{\prime
}\right) \equiv \Sigma ^{si^{\prime }}\Sigma _{,k^{\prime }l^{\prime
}s},\qquad \Sigma _{,k^{\prime }l^{\prime }s}\equiv \frac{\partial
^{3}\Sigma }{\partial x^{\prime k}\partial x^{\prime l}\partial x^{s}}
\label{c3.7}
\end{equation}
is a scalar at the point $x$ and a Christoffel symbol at the point $%
x^{\prime }$.

In the same way one can introduce two other Christoffel symbols on the basis
of the function $G$%
\begin{eqnarray}
\Gamma _{kl}^{i}\left( x,x^{\prime }\right) &\equiv &G^{is^{\prime
}}G_{,kls^{\prime }},\qquad G_{,kls^{\prime }}\equiv \frac{\partial ^{3}G}{%
\partial x^{k}\partial x^{l}\partial x^{\prime s}},\qquad G^{is^{\prime
}}G_{,ks^{\prime }}=\delta _{k}^{i}  \label{c3.7a} \\
\Gamma _{k^{\prime }l^{\prime }}^{i^{\prime }}\left( x,x^{\prime }\right)
&\equiv &G^{si^{\prime }}G_{,k^{\prime }l^{\prime }s},\qquad G_{,k^{\prime
}l^{\prime }s}\equiv \frac{\partial ^{3}G}{\partial x^{\prime k}\partial
x^{\prime l}\partial x^{s}}  \label{c3.7b}
\end{eqnarray}

Using Christoffel symbols (\ref{c3.6}) - (\ref{c3.7b}), one can introduce
two covariant derivatives $\tilde{\nabla}_{i}^{x^{\prime }}$, $\nabla
_{i}^{x^{\prime }}$ with respect to $x^{i}$ and two covariant derivatives $%
\tilde{\nabla}_{i^{\prime }}^{x}$, $\nabla _{i^{\prime }}^{x}$ with respect
to $x^{\prime i}.$ For instance, the quantities 
\begin{eqnarray}
\Sigma _{ik} &\equiv &\tilde{\nabla}_{k}^{x^{\prime }}\tilde{\nabla}%
_{i}^{x^{\prime }}\Sigma \equiv \Sigma _{,i||k}=\Sigma _{,ik}-\tilde{\Gamma}%
_{ik}^{s}\left( x,x^{\prime }\right) \Sigma _{,s}=\Sigma _{,ik}-\Sigma
^{ls^{\prime }}\Sigma _{,iks^{\prime }}\Sigma _{,l}  \label{c3.8} \\
G_{ik} &\equiv &\nabla _{k}^{x^{\prime }}\nabla _{i}^{x^{\prime }}G\equiv
G_{,i|k}=G_{,ik}-\Gamma _{ik}^{s}\left( x,x^{\prime }\right)
G_{,s}=G_{,ik}-G^{ls^{\prime }}G_{,iks^{\prime }}G_{,l}  \label{c3.8a}
\end{eqnarray}
are scalars at the point $x^{\prime }$ and a second rank tensors at the
point $x$. Here two vertical strokes denote covariant derivative with the
Christoffel symbol $\tilde{\Gamma}_{ik}^{s}$, and one vertical stroke denote
covariant derivative with the Christoffel symbol $\Gamma _{ik}^{s}$. In the
same way one obtains 
\begin{eqnarray}
\Sigma _{i^{\prime }k^{\prime }} &\equiv &\tilde{\nabla}_{k^{\prime }}^{x}%
\tilde{\nabla}_{i^{\prime }}^{x}\Sigma \equiv \Sigma _{,i^{\prime
}||k^{\prime }}=\Sigma _{,i^{\prime }k^{\prime }}-\tilde{\Gamma}_{i^{\prime
}k^{\prime }}^{s^{\prime }}\Sigma _{,s^{\prime }}=\Sigma _{,i^{\prime
}k^{\prime }}-\Sigma ^{l^{\prime }s}\Sigma _{,i^{\prime }k^{\prime }s}\Sigma
_{,l^{\prime }}  \label{c3.9} \\
G_{i^{\prime }k^{\prime }} &\equiv &\nabla _{k^{\prime }}^{x}\nabla
_{i^{\prime }}^{x}G\equiv G_{,i^{\prime }|k^{\prime }}=G_{,i^{\prime
}k^{\prime }}-\Gamma _{i^{\prime }k^{\prime }}^{s^{\prime }}G_{,s^{\prime
}}=G_{,i^{\prime }k^{\prime }}-G^{l^{\prime }s}G_{,i^{\prime }k^{\prime
}s}G_{,l^{\prime }}  \label{c3.9a}
\end{eqnarray}

Covariant derivatives $\tilde{\nabla}_{k}^{x^{\prime }}$, $\tilde{\nabla}%
_{i}^{x^{\prime }}$ with respect to $x$ commute, as well as $\nabla
_{k}^{x^{\prime }}$, $\nabla _{i}^{x^{\prime }}$, i.e. 
\begin{equation}
\left( \tilde{\nabla}_{k}^{x^{\prime }}\tilde{\nabla}_{i}^{x^{\prime }}-%
\tilde{\nabla}_{i}^{x^{\prime }}\tilde{\nabla}_{k}^{x^{\prime }}\right)
T_{ml^{\prime }}^{sp^{\prime }}\equiv 0,\qquad \left( \nabla _{k}^{x^{\prime
}}\nabla _{i}^{x^{\prime }}-\nabla _{i}^{x^{\prime }}\nabla _{k}^{x^{\prime
}}\right) T_{ml^{\prime }}^{sp^{\prime }}\equiv 0  \label{c3.10}
\end{equation}
where $T_{ml^{\prime }}^{sp^{\prime }}$ is an arbitrary tensor at points $x$
and $x^{\prime }$. Unprimed indices are associated with the point $x$, and
primed ones with the point $x^{\prime }$. The covariant derivatives commute,
because the Riemann-Christoffel curvature tensors $\tilde{R}_{i.kl}^{.j}$, $%
R_{i.kl}^{.j}$ constructed respectively of Christoffel symbols $\tilde{\Gamma%
}_{il}^{s}$ and $\Gamma _{il}^{s}$ vanish identically 
\begin{eqnarray}
\tilde{R}_{i.lm}^{.s} &\equiv &\tilde{\Gamma}_{il,m}^{s}-\tilde{\Gamma}%
_{im,l}^{s}+\tilde{\Gamma}_{il}^{j}\tilde{\Gamma}_{jm}^{s}-\tilde{\Gamma}%
_{im}^{j}\tilde{\Gamma}_{jl}^{s}\equiv 0  \label{c3.11} \\
R_{i.lm}^{.s} &\equiv &\Gamma _{il,m}^{s}-\Gamma _{im,l}^{s}+\Gamma
_{il}^{j}\Gamma _{jm}^{s}-\Gamma _{im}^{j}\Gamma _{jl}^{s}\equiv 0
\label{c3.11a}
\end{eqnarray}
One can test the identity (\ref{c3.11}), substituting (\ref{c3.6}) into (\ref
{c3.11}).

Covariant derivatives $\tilde{\nabla}_{k^{\prime }}^{x}$, $\tilde{\nabla}%
_{i^{\prime }}^{x}$ with respect to $x^{\prime }$ commute as well as $\nabla
_{k^{\prime }}^{x}$, $\nabla _{i^{\prime }}^{x}$. Commutativity of covariant
derivatives $\tilde{\nabla}_{i}^{x^{\prime }}$, $\tilde{\nabla}%
_{k}^{x^{\prime }}$ with respect to $x$ for all values of $x^{\prime }$
means that the covariant derivative $\tilde{\nabla}_{i}^{x^{\prime }}$, $%
\tilde{\nabla}_{k}^{x^{\prime }}$ are covariant derivatives in some flat
spaces $\tilde{E}_{x^{\prime }}$. The same is valid for covariant
derivatives $\nabla _{i}^{x^{\prime }}$, $\nabla _{k}^{x^{\prime }}$ which
are covariant derivatives in the flat spaces $E_{x^{\prime }}$. The spaces $%
\tilde{E}_{x^{\prime }},$ $E_{x^{\prime }}$ are associated with the $\Sigma $%
-spaces $V=\left\{ \Sigma ,{\cal M}_{n+1}\right\} $ and $V_{{\rm s}}=\left\{
G,{\cal M}_{n+1}\right\} $ respectively, given on ${\cal M}_{n+1}$ by means
of the world function $\Sigma $ and its symmetric part $G$. Any of two-point
invariant quantities $\Sigma $ and $G$ with nonvanishing determinants $\det
||\Sigma _{ik^{\prime }}||\neq 0,$ and $\det ||G_{ik^{\prime }}||\neq 0$
realize two sets of mappings. For instance, the quantity $\Sigma $ generates
mappings $V=\left\{ \Sigma ,{\cal M}_{n+1}\right\} \rightarrow \tilde{E}%
_{x^{\prime }}$, $\;V\rightarrow \tilde{E}_{x}$ The two-point quantity $G$
generates also two sets of mappings $V_{{\rm s}}=\left\{ G,{\cal M}%
_{n+1}\right\} \rightarrow E_{x^{\prime }}$, $\;V_{{\rm s}}\rightarrow E_{x}$%
. Mappings of any set are labelled by points $x$ or $x^{\prime }$ of the
manifold ${\cal M}_{n+1}$. In the case of $G$ both sets of mappings $%
V\rightarrow E_{x^{\prime }}$ and $\;V\rightarrow E_{x}$ coincide, but in
the case of $\Sigma $ the sets $V\rightarrow \tilde{E}_{x^{\prime }}$ and $%
\;V\rightarrow \tilde{E}_{x}$ are different, in general.

It is easy to see that 
\[
\tilde{\nabla}_{s}^{x^{\prime }}\Sigma _{ik^{\prime }}=\Sigma _{ik^{\prime
}||s}=\Sigma _{,ik^{\prime }s}-\tilde{\Gamma}_{is}^{p}\Sigma _{pk^{\prime
}}=\Sigma _{,ik^{\prime }s}-\Sigma ^{pl^{\prime }}\Sigma _{isl^{\prime
}}\Sigma _{pk^{\prime }}\equiv 0 
\]

The covariant derivatives have the following properties 
\begin{equation}
\tilde{\nabla}_{k}^{x^{\prime }}t_{l^{\prime }}^{i^{\prime }}\left(
x^{\prime }\right) =0,\qquad \nabla _{k}^{x^{\prime }}t_{l^{\prime
}}^{i^{\prime }}\left( x^{\prime }\right) =0,\qquad \tilde{\nabla}%
_{k^{\prime }}^{x}t_{l}^{i}\left( x\right) =0\qquad \nabla _{k^{\prime
}}^{x}t_{l}^{i}\left( x\right) =0  \label{c3.12}
\end{equation}
\begin{equation}
\tilde{\nabla}_{s}^{x^{\prime }}\Sigma _{ik^{\prime }}=\Sigma _{ik^{\prime
}||s}=0\qquad \Sigma _{ik^{\prime }||s^{\prime }}=0,\qquad G_{ik^{\prime
}|s}=0\qquad G_{ik^{\prime }|s^{\prime }}=0,  \label{c3.13}
\end{equation}
where $t_{l^{\prime }}^{i^{\prime }}\left( x^{\prime }\right) $ is an
arbitrary tensor at the point $x^{\prime }$, and $t_{l}^{i}\left( x\right) $
is an arbitrary tensor at the point $x$.

The mappings $V_{{\rm s}}\rightarrow E_{x^{\prime }}$, $\;V_{{\rm s}%
}\rightarrow E_{x}$ are associated with the mappings (\ref{aa2.1}). It means
that description of $\Sigma $-space, given on a manifold, involves the set
of mappings (\ref{aa2.1}) and some other mappings of the $\Sigma $-space
onto Euclidean spaces. These mappings can serve as a powerful tool for
description of the $\Sigma $-space properties. Let us note in this
connection, that the Riemannian space may be considered to be a set of
infinitesimal pieces of Euclidean spaces glued in some way between
themselves. The way of gluing determines the character of the Riemannian
space in the sense, that different ways of gluing generate different
Riemannian spaces. The way of gluing is determined by the difference of the
metric tensor at the point $x$ and at the narrow point $x+dx$, where it has
the forms $g_{ik}\left( x\right) $ and $g_{ik}\left( x+dx\right) $
respectively. The metric tensor depends on a point and on a choice of the
coordinate system. It is rather difficult one to separate dependence on the
way of gluing from that on the choice of the coordinate system. Nevertheless
the procedure of separation has been well developed. It leads to the
curvature tensor, which is an indicator of the way of gluing.

In the case of the $\Sigma $-space one considers a set of {\it finite
Euclidean spaces} \ $E_{x^{\prime }}\;$(ins\-tead of its infinitesimal
pieces) and a set of mappings $\Sigma \rightarrow E_{x^{\prime }}$. Here the
''way of gluing'' is determined by the dependence of mapping on the
parameter $x^{\prime }$ and does not depend on a choice of the coordinate
system. This circumstance simplifies investigation. Differentiating the
mappings with respect to parameters $x^{\prime }$, one derives local
characteristics of the ''way of gluing'', which are modifications of the
curvature tensor. For instance, considering commutators of derivatives $%
\tilde{\nabla}_{i}^{x^{\prime }}$ and $\tilde{\nabla}_{i^{\prime }}^{x}$,
one can introduce two-point curvature tensor for the $\Sigma $-space, as it
have been made for the Riemannian space \cite{R62,R64}. We shall see this in
the sixth section.

Let $\tilde{G}_{(x^{\prime })ik}$ be the metric tensor in the Euclidean
space $\tilde{E}_{x^{\prime }}$ at the point $x$. Then the Christoffel
symbol $\tilde{\Gamma}_{kl}^{i}=\Sigma ^{is^{\prime }}\Sigma _{,kls^{\prime
}}$ in the space $\tilde{E}_{x^{\prime }}$ can be written in the form

\begin{equation}
\tilde{\Gamma}_{kl}^{i}=\Sigma ^{is^{\prime }}\Sigma _{,kls^{\prime }}=\frac{%
1}{2}\tilde{G}_{(x^{\prime })}^{im}\left( \tilde{G}_{(x^{\prime })km,l}+%
\tilde{G}_{(x^{\prime })lm,k}-\tilde{G}_{(x^{\prime })kl,m}\right)
\label{c3.14}
\end{equation}
where $\tilde{G}_{(x^{\prime })}^{im}$ are contravariant components of the
metric tensor $\tilde{G}_{(x^{\prime })ik}$.

Let us consider the set of equations (\ref{c3.14}) as a system of linear
differential equations for determination of the metric tensor components $%
\tilde{G}_{(x^{\prime })ik}$, which is supposed to be symmetric. Solution of
this system has the form 
\begin{equation}
\tilde{G}_{(x^{\prime })ik}=\Sigma _{ip^{\prime }}\tilde{g}_{\left(
x^{\prime }\right) }^{p^{\prime }q^{\prime }}\Sigma _{kq^{\prime }}
\label{c3.15}
\end{equation}
where $\tilde{g}_{\left( x^{\prime }\right) }^{p^{\prime }q^{\prime }}=%
\tilde{g}_{\left( x^{\prime }\right) }^{p^{\prime }q^{\prime }}\left(
x^{\prime }\right) $ is some symmetric tensor at the point $x^{\prime }$.
This fact can be tested by a direct substitution of (\ref{c3.15}) in (\ref
{c3.14}). Taking the relation (\ref{c3.15}) at the coinciding points $%
x=x^{\prime }$ and denoting coincidence of points $x$ and $x^{\prime }$ by
means of square brackets, one obtains from (\ref{c3.15}) 
\begin{equation}
\tilde{g}_{\left( x^{\prime }\right) }^{p^{\prime }q^{\prime }}\left(
x^{\prime }\right) =\left[ \Sigma ^{lp^{\prime }}\right] _{x^{\prime }}\left[
\tilde{G}_{(x^{\prime })lm}\right] _{x^{\prime }}\left[ \Sigma ^{mq^{\prime
}}\right] _{x^{\prime }}  \label{c3.16}
\end{equation}
or 
\begin{equation}
\tilde{G}_{(x^{\prime })ik}=\Sigma _{ip^{\prime }}\left[ \Sigma ^{lp^{\prime
}}\right] _{x^{\prime }}\left[ \tilde{G}_{(x^{\prime })lm}\right]
_{x^{\prime }}\left[ \Sigma ^{mq^{\prime }}\right] _{x^{\prime }}\Sigma
_{kq^{\prime }}  \label{c3.17}
\end{equation}
The equation (\ref{c3.17}) can be written in the form 
\begin{eqnarray}
\tilde{G}_{(x^{\prime })ik}\left( x,x^{\prime }\right) &=&\tilde{P}%
_{(x^{\prime })}{}_{i.}^{.l^{\prime }}\tilde{G}_{(x^{\prime })lm}\left(
x^{\prime },x^{\prime }\right) \tilde{P}_{(x^{\prime })}{}_{k.}^{.m^{\prime
}},  \label{c3.18} \\
\tilde{P}_{(x^{\prime })}{}_{k.}^{.m^{\prime }} &\equiv &\tilde{P}%
_{(x^{\prime })}{}_{k.}^{.m^{\prime }}\left( x,x^{\prime }\right) \equiv
\Sigma _{kq^{\prime }}\left( x,x^{\prime }\right) \Sigma ^{mq^{\prime
}}\left( x^{\prime },x^{\prime }\right) \equiv \Sigma _{kq^{\prime }}\left[
\Sigma ^{mq^{\prime }}\right] _{x^{\prime }}  \label{c3.19}
\end{eqnarray}

The relation (\ref{c3.18}) means that the metric tensor$\ \tilde{G}%
_{(x^{\prime })ik}$ of the Euclidean space $\tilde{E}_{x^{\prime }}$ at the
point $x$ can be obtained as a result of the parallel transport of the
metric tensor from the point $x^{\prime }$ in $\tilde{E}_{x^{\prime }}$ by
means of the parallel transport tensor $\tilde{P}_{(x^{\prime
})}{}_{k.}^{.m^{\prime }}$. The parallel transport of the vector $%
b_{k^{\prime }}$ from the point $x^{\prime }$ to the point $x$ is defined by
the relation 
\[
b_{k}=\tilde{P}_{(x^{\prime })}{}_{k.}^{.m^{\prime }}b_{m^{\prime }}. 
\]
The parallel transport tensor has evident properties 
\begin{eqnarray}
\tilde{\nabla}_{i}^{x^{\prime }}\tilde{P}_{(x^{\prime })}{}_{k.}^{.m^{\prime
}} &\equiv &\tilde{P}_{(x^{\prime })}{}_{k.||i}^{.m^{\prime }}=0,
\label{c3.20} \\
\left[ \tilde{P}_{(x^{\prime })}{}_{k.}^{.m^{\prime }}\right] _{x^{\prime }}
&\equiv &\tilde{P}_{(x^{\prime })}{}_{k.}^{.m^{\prime }}\left( x^{\prime
},x^{\prime }\right) =\delta _{k^{\prime }}^{m^{\prime }}.  \label{c3.21}
\end{eqnarray}

In the same way one can obtain the parallel transport tensor $\tilde{P}%
_{(x)}{}_{k^{\prime }.}^{.m}$ in the Euclidean space $\tilde{E}_{x}$%
\begin{equation}
\tilde{P}_{(x)}{}_{k^{\prime }.}^{.m}\equiv \Sigma _{k^{\prime }q}\left[
\Sigma ^{mq^{\prime }}\right] _{x}\equiv \Sigma _{k^{\prime }q}\left(
x,x^{\prime }\right) \Sigma ^{mq^{\prime }}\left( x,x\right)  \label{c3.22}
\end{equation}
describing a parallel transport from the point $x$ to the point $x^{\prime }$
in $\tilde{E}_{x}$.

In the same way one can obtain the parallel transport tensors $P_{(x^{\prime
})}{}_{k.}^{.m^{\prime }}$ and $P_{(x)}{}_{k^{\prime }.}^{.m}$ respectively
in Euclidean spaces $E_{x^{\prime }}$ and $E_{x}$%
\begin{equation}
P_{(x^{\prime })}{}_{k.}^{.m^{\prime }}\equiv G_{kq^{\prime }}\left[
G^{mq^{\prime }}\right] _{x^{\prime }},\qquad P_{(x)}{}_{k^{\prime
}.}^{.m}\equiv G_{k^{\prime }q}\left[ G^{mq^{\prime }}\right] _{x}
\label{c3.23}
\end{equation}

Thus, the world function $\Sigma $ of the $\Sigma $-space $V=\left\{ \Sigma ,%
{\cal M}_{n+1}\right\} $ and its symmetric component $G$ determine Euclidean
spaces $\tilde{E}_{x^{\prime }}$, $\tilde{E}_{x}$, $E_{x^{\prime }} $, $%
E_{x} $, mappings of $V$ on them and the parallel transport of vectors and
tensors in these Euclidean spaces independently of that, whether or not the $%
\Sigma $-space $V=\left\{ \Sigma ,{\cal M}_{n+1}\right\} $ is degenerate in
the sense of definitions \ref{d3.12a} -- \ref{d3.14}.

\section{Derivatives of the world function at coincidence of points $x$ and $%
x^{\prime }$.}

Let us represent the world function in the form 
\begin{eqnarray}
\Sigma \left( x,x^{\prime }\right) &=&G\left( x,x^{\prime }\right) +A\left(
x,x^{\prime }\right) ,  \label{c5.1} \\
G\left( x,x^{\prime }\right) &=&G\left( x^{\prime },x\right) ,\qquad A\left(
x,x^{\prime }\right) =-A\left( x^{\prime },x\right)  \label{c5.2}
\end{eqnarray}
Here $G$ is the symmetric component of the world function $\Sigma $, and $A$
is its antisymmetric component.

Let us expand the function $G$ and $A$ with respect to powers of $\xi
^{i}=x^{i}-x^{\prime i}$. Taking into account the symmetry relations (\ref
{c5.1}), (\ref{c5.2}), one obtains 
\begin{eqnarray}
G\left( x,x^{\prime }\right) &=&\frac{1}{2}g_{ik}\left( x^{\prime }\right)
\xi ^{i}\xi ^{k}+\frac{1}{6}g_{ikl}\left( x^{\prime }\right) \xi ^{i}\xi
^{k}\xi ^{l}+\frac{1}{24}g_{iklm}\left( x^{\prime }\right) \xi ^{i}\xi
^{k}\xi ^{l}\xi ^{m}+...  \label{c5.3} \\
&=&\frac{1}{2}g_{ik}\left( x\right) \xi ^{i}\xi ^{k}-\frac{1}{6}%
g_{ikl}\left( x\right) \xi ^{i}\xi ^{k}\xi ^{l}+\frac{1}{24}g_{iklm}\left(
x\right) \xi ^{i}\xi ^{k}\xi ^{l}\xi ^{m}+...  \label{c5.4}
\end{eqnarray}

\begin{eqnarray}
A\left( x,x^{\prime }\right) &=&a_{i}\left( x^{\prime }\right) \xi ^{i}+%
\frac{1}{2}a_{ik}\left( x^{\prime }\right) \xi ^{i}\xi ^{k}+\frac{1}{6}%
a_{ikl}\left( x^{\prime }\right) \xi ^{i}\xi ^{k}\xi ^{l}+...  \label{c5.5}
\\
&=&a_{i}\left( x\right) \xi ^{i}-\frac{1}{2}a_{ik}\left( x\right) \xi
^{i}\xi ^{k}+\frac{1}{6}a_{ikl}\left( x\right) \xi ^{i}\xi ^{k}\xi ^{l}-...
\label{c5.6}
\end{eqnarray}
In relations (\ref{c5.3}) and (\ref{c5.5}) the functions $G\left(
x,x^{\prime }\right) $ and $A\left( x,x^{\prime }\right) $ are expanded at
the point $x^{\prime }$. In the relations (\ref{c5.4}) and (\ref{c5.6}) one
has the same expansions after transposition $x\leftrightarrow x^{\prime }$.

Differentiating relations (\ref{c5.3}) - (\ref{c5.6}) with respect to $x$
and $x^{\prime }$ and setting $x=x^{\prime }$ thereafter, one obtains
relations between the expansion coefficients and expressions for derivatives
of functions $\Sigma $, $G$, $A$ at the limit of coincidence $x=x^{\prime }$.

After calculations one obtains 
\begin{eqnarray}
g_{ikl} &=&\frac{1}{2}\left( g_{ik,l}+g_{li,k}+g_{kl,i}\right) ,\qquad
g_{ik,l}\equiv g_{ik,l}\left( x\right) \equiv \frac{\partial }{\partial x^{l}%
}g_{ik}\left( x\right)  \label{c5.7} \\
a_{ik} &=&\frac{1}{2}\left( a_{i,k}+a_{k,i}\right) ,\qquad a_{i,k}\equiv
a_{i,k}\left( x\right) \equiv \frac{\partial }{\partial x^{k}}a_{i}\left(
x\right)  \label{c5.8}
\end{eqnarray}
\[
a_{iklm}=\frac{1}{2}\left(
a_{ikl,m}+a_{klm,i}+a_{lmi,k}+a_{mik,l}-a_{ik,lm}-a_{lm,ik}\right) 
\]
Coefficients $g_{ik},$ $a_{i}$, $a_{ikl}$ are arbitrary and symmetric with
respect to transposition of indices. Using square brackets for designation
of coincidence $x=x^{\prime }$ and relations (\ref{c5.7}), (\ref{c5.8}), one
obtains 
\begin{eqnarray}
\left[ G_{,i}\left( x,x^{\prime }\right) \right] &\equiv &\left[ G_{,i}%
\right] =\left[ G_{,i^{\prime }}\right] =0,  \nonumber \\
\left[ G_{,ik}\left( x,x^{\prime }\right) \right] &\equiv &\left[ G_{,ik}%
\right] =\left[ G_{,i^{\prime }k^{\prime }}\right] =g_{ik},  \label{c5.9} \\
\left[ G_{,i^{\prime }k}\left( x,x^{\prime }\right) \right] &\equiv &\left[
G_{,i^{\prime }k}\right] =-g_{ik}  \nonumber
\end{eqnarray}

\begin{eqnarray}
\left[ G_{,ikl}\left( x,x^{\prime }\right) \right] &\equiv &\left[ G_{,ikl}%
\right] =\left[ G_{,i^{\prime }k^{\prime }l^{\prime }}\right] =\frac{1}{2}%
\left( g_{ik,l}+g_{li,k}+g_{kl,i}\right)  \nonumber \\
\left[ G_{,ikl^{\prime }}\left( x,x^{\prime }\right) \right] &\equiv &\left[
G_{,ikl^{\prime }}\right] =\frac{1}{2}\left(
g_{ik,l}-g_{li,k}-g_{kl,i}\right)  \label{c5.10} \\
\left[ G_{,ik^{\prime }l^{\prime }}\left( x,x^{\prime }\right) \right]
&\equiv &\left[ G_{,ik^{\prime }l^{\prime }}\right] =\frac{1}{2}\left(
g_{kl,i}-g_{li,k}-g_{ik,l}\right)  \nonumber
\end{eqnarray}
\begin{eqnarray}
\left[ G_{,iklm}\left( x,x^{\prime }\right) \right] &\equiv &\left[ G_{,iklm}%
\right] =\left[ G_{,i^{\prime }k^{\prime }l^{\prime }m^{\prime }}\right]
=g_{iklm}  \nonumber \\
\left[ G_{,iklm^{\prime }}\left( x,x^{\prime }\right) \right] &\equiv &\left[
G_{,iklm^{\prime }}\right] =\left[ G_{,i^{\prime }k^{\prime }l^{\prime }m}%
\right] =-g_{iklm}+g_{ikl,m}  \label{c5.10a} \\
\left[ G_{,ikl^{\prime }m^{\prime }}\left( x,x^{\prime }\right) \right]
&\equiv &\left[ G_{,ikl^{\prime }m^{\prime }}\right]
=g_{iklm}-g_{ikl,m}-g_{ikm,l}+g_{ik,ml}  \nonumber
\end{eqnarray}
\begin{equation}
\left[ A_{,i}\left( x,x^{\prime }\right) \right] \equiv \left[ A_{,i}\right]
=a_{i},\qquad \left[ A_{,i^{\prime }}\left( x,x^{\prime }\right) \right]
\equiv \left[ A_{,i^{\prime }}\right] =-a_{i}  \label{c5.11}
\end{equation}
\begin{eqnarray}
\left[ A_{,ik}\left( x,x^{\prime }\right) \right] &\equiv &\left[ A_{,ik}%
\right] =a_{ik}=\frac{1}{2}\left( a_{i,k}+a_{k,i}\right)  \nonumber \\
\left[ A_{,ik^{\prime }}\left( x,x^{\prime }\right) \right] &\equiv &\left[
A_{,ik^{\prime }}\right] =\frac{1}{2}\left( a_{i,k}-a_{k,i}\right)
\label{c5.12} \\
\left[ A_{,i^{\prime }k^{\prime }}\left( x,x^{\prime }\right) \right]
&\equiv &\left[ A_{,i^{\prime }k^{\prime }}\right] =-a_{ik}=-\frac{1}{2}%
\left( a_{i,k}+a_{k,i}\right)  \nonumber
\end{eqnarray}
\begin{eqnarray}
\left[ A_{,ikl}\left( x,x^{\prime }\right) \right] &\equiv &\left[ A_{,ikl}%
\right] =a_{ikl}  \nonumber \\
\left[ A_{,ikl^{\prime }}\left( x,x^{\prime }\right) \right] &\equiv &\left[
A_{,ikl^{\prime }}\right] =\frac{1}{2}a_{i,kl}+\frac{1}{2}a_{k,il}-a_{ikl}
\label{c5.12a} \\
\left[ A_{,ik^{\prime }l^{\prime }}\left( x,x^{\prime }\right) \right]
&\equiv &\left[ A_{,ik^{\prime }l^{\prime }}\right] =-\frac{1}{2}a_{l,ik}-%
\frac{1}{2}a_{k,il}+a_{ikl}  \nonumber
\end{eqnarray}
\begin{eqnarray}
\left[ A_{,iklm}\left( x,x^{\prime }\right) \right] &\equiv &\left[ A_{,iklm}%
\right] =a_{iklm}  \nonumber \\
\left[ A_{,iklm^{\prime }}\left( x,x^{\prime }\right) \right] &\equiv &\left[
A_{,iklm^{\prime }}\right] =-a_{iklm}+a_{ikl,m}  \nonumber \\
\left[ A_{,ikl^{\prime }m^{\prime }}\left( x,x^{\prime }\right) \right]
&\equiv &\left[ A_{,ikl^{\prime }m^{\prime }}\right]
=a_{iklm}-a_{ikl,m}-a_{ikm,l}+a_{ik,lm}  \label{c5.13a} \\
\left[ A_{,ik^{\prime }l^{\prime }m^{\prime }}\left( x,x^{\prime }\right) %
\right] &\equiv &\left[ A_{,ik^{\prime }l^{\prime }m^{\prime }}\right]
=-a_{iklm}-a_{klm,i}  \nonumber \\
\left[ A_{,i^{\prime }k^{\prime }l^{\prime }m^{\prime }}\left( x,x^{\prime
}\right) \right] &\equiv &\left[ A_{,i^{\prime }k^{\prime }l^{\prime
}m^{\prime }}\right] =-a_{iklm}  \nonumber
\end{eqnarray}

The first order coefficient $a_{k}\left( x\right) $ is a covariant vector at
the point $x$. The second order coefficients $g_{ik}\left( x\right) $ is the
second rank covariant tensors at the point $x.$ The second order coefficient 
$a_{ik}\left( x\right) $ and the third order coefficients $g_{ikl}\left(
x\right) $, $a_{ikl}\left( x\right) $ are not tensors, in general. The law
of their transformation at the coordinate transformation is more complicated.

According to (\ref{c5.1}), (\ref{c5.9}) and (\ref{c5.11}) 
\begin{equation}
\left[ \Sigma _{,i}\right] =a_{i}\left( x\right) ,\qquad \left[ \Sigma
_{,i^{\prime }}\right] =-a_{i}\left( x\right)  \label{c5.14}
\end{equation}
According to (\ref{c5.1}), (\ref{c5.9}) and (\ref{c5.12})\label{01} 
\begin{eqnarray}
\left[ \Sigma _{,ik}\right] &=&\sigma _{\left( {\rm f}\right) ik}=g_{ik}+%
\frac{1}{2}\left( a_{i,k}+a_{k,i}\right)  \nonumber \\
\left[ \Sigma _{,i^{\prime }k^{\prime }}\right] &=&\sigma _{\left( {\rm p}%
\right) ik}=g_{ik}-\frac{1}{2}\left( a_{i,k}+a_{k,i}\right)  \label{c5.15} \\
\left[ \Sigma _{ik^{\prime }}\right] &\equiv &\left[ \Sigma _{,ik^{\prime }}%
\right] =-\tilde{g}_{ik}=-g_{ik}+\frac{1}{2}\left( a_{i,k}-a_{k,i}\right) 
\nonumber
\end{eqnarray}
For the value $\left[ \Sigma ^{ik^{\prime }}\right] $ of the quantity $%
\Sigma ^{ik^{\prime }}$ one obtains 
\begin{equation}
\left[ \Sigma ^{ik^{\prime }}\right] =-\tilde{g}^{ik},  \label{c5.16}
\end{equation}
where $\tilde{g}^{ik}$ is determined by the relation 
\begin{equation}
\tilde{g}^{il}\tilde{g}_{ik}=\tilde{g}^{il}\left( g_{ik}-\frac{1}{2}\left(
a_{i,k}-a_{k,i}\right) \right) =\delta _{k}^{l}  \label{c5.17}
\end{equation}
The quantity $a_{i}$ is a one-point vector, and $g_{ik}$ is a one-point
tensor. Then it follows from (\ref{c5.15}), (\ref{c5.17}), that $\tilde{g}%
_{ik}$ and $\tilde{g}^{ik}$ are also a one-point tensors, whereas $\sigma
_{\left( {\rm f}\right) ik}$ and $\sigma _{\left( {\rm p}\right) ik}$ are
not tensors, in general.

For the quantities $\left[ \Gamma _{kl}^{i}\right] ,$ $\left[ \Gamma
_{k^{\prime }l^{\prime }}^{i^{\prime }}\right] $ and $\left[ \tilde{\Gamma}%
_{kl}^{i}\right] ,$ $\left[ \tilde{\Gamma}_{k^{\prime }l^{\prime
}}^{i^{\prime }}\right] $ one obtains 
\begin{eqnarray}
\left[ \Gamma _{kl}^{i}\right] &=&\left[ \Gamma _{k^{\prime }l^{\prime
}}^{i^{\prime }}\right] =\gamma _{kl}^{i}\left( x\right) =\frac{1}{2}%
g^{si}\left( g_{ks,l}+g_{sl,k}-g_{lk,s}\right) ,\qquad g^{ik}g_{lk}=\delta
_{l}^{i}  \label{c5.18a} \\
\left[ \tilde{\Gamma}_{kl}^{i}\right] &=&\tilde{\gamma}_{\left( {\rm f}%
\right) }{}_{kl}^{i}=\tilde{g}^{is}g_{ps}\left( \gamma _{kl}^{p}+\beta
_{kl}^{p}\right) ,  \label{c5.18} \\
\left[ \tilde{\Gamma}_{k^{\prime }l^{\prime }}^{i^{\prime }}\right] &=&%
\tilde{\gamma}_{\left( {\rm p}\right) }{}_{kl}^{i}=\tilde{g}%
^{si}g_{ps}\left( \gamma _{kl}^{p}-\beta _{kl}^{p}\right)  \label{c5.18b}
\end{eqnarray}
where 
\begin{eqnarray}
\gamma _{kl}^{i}\left( x\right) &=&\frac{1}{2}g^{si}\left(
g_{ks,l}+g_{sl,k}-g_{lk,s}\right) ,\qquad g^{ik}g_{lk}=\delta _{l}^{i}
\label{c5.20} \\
\beta _{kl}^{i}\left( x\right) &=&g^{si}\left( -\frac{1}{2}\left(
a_{k,ls}+a_{l,ks}\right) +a_{kls}\right)  \label{c5.21}
\end{eqnarray}
Here $\gamma _{kl}^{i}\left( x\right) $ is the Christoffel symbol for the
symmetric case, when $A\left( x,x^{\prime }\right) \equiv 0$.

Note that the tensors $\tilde{g}^{ik},$ $\tilde{g}_{ik}$ are not symmetric
with respect to transposition indices, in general, whereas $\tilde{\gamma}%
_{\left( {\rm f}\right) }{}_{kl}^{i}=\left[ \tilde{\Gamma}_{kl}^{i}\right] $%
, $\;\tilde{\gamma}_{\left( {\rm p}\right) }{}_{kl}^{i}=\left[ \tilde{\Gamma}%
_{k^{\prime }l^{\prime }}^{i^{\prime }}\right] $, $\gamma _{kl}^{i}$ and $%
\beta _{kl}^{i}$ are symmetric with respect to transposition of indices $k$
and $l$. Besides it follows from (\ref{c5.18}), (\ref{c5.18b}), that 
\begin{equation}
\beta _{kl}^{i}\left( x\right) =g^{si}\left( -\frac{1}{2}\left(
a_{k,ls}+a_{l,ks}\right) +a_{kls}\right) =\frac{1}{2}g^{ip}\tilde{g}%
_{ps}\left( \tilde{\gamma}_{\left( {\rm f}\right) }{}_{kl}^{s}-\tilde{\gamma}%
_{\left( {\rm p}\right) }{}_{kl}^{s}\right)  \label{c5.22}
\end{equation}
It means that the quantity $\beta _{kl}^{i}$ is a one-point tensor, because
a difference of two Christoffel symbols $\tilde{\gamma}_{\left( {\rm f}%
\right) }{}_{kl}^{s}-\tilde{\gamma}_{\left( {\rm p}\right) }{}_{kl}^{s}$ is
a tensor.

\section{Curvature tensors}

In the Riemannian geometry the Riemann-Christoffel curvature tensor $\tilde{r%
}_{\left( {\rm q}\right) }{}_{s.ik}^{.l}$ is defined as a commutator of
covariant derivatives $\tilde{D}_{\left( {\rm q}\right) i}$ with the
Christoffel symbol $\tilde{\gamma}_{\left( {\rm q}\right) }{}_{kl}^{i}$ 
\[
\left( \tilde{D}_{\left( {\rm q}\right) i}\tilde{D}_{\left( {\rm q}\right)
k}-\tilde{D}_{\left( {\rm q}\right) k}\tilde{D}_{\left( {\rm q}\right)
i}\right) t_{s}=\tilde{r}_{\left( {\rm q}\right) }{}_{s.ik}^{.l}t_{l} 
\]
where $\tilde{D}_{\left( {\rm q}\right) i}$ is the usual covariant
derivative in the Riemannian space with the Christoffel symbol $\tilde{\gamma%
}_{\left( {\rm q}\right) }{}_{si}^{l}$, 
\begin{equation}
\tilde{r}_{\left( {\rm q}\right) }{}_{s.ik}^{.l}=\tilde{\gamma}_{\left( {\rm %
q}\right) }{}_{si,k}^{l}-\tilde{\gamma}_{\left( {\rm q}\right)
}{}_{sk,i}^{l}+\tilde{\gamma}_{\left( {\rm q}\right) }{}_{si}^{p}\tilde{%
\gamma}_{\left( {\rm q}\right) }{}_{pk}^{l}-\tilde{\gamma}_{\left( {\rm q}%
\right) }{}_{sk}^{p}\tilde{\gamma}_{\left( {\rm q}\right) }{}_{pi}^{l}
\label{f6.2}
\end{equation}
and $t_{l}$ is an arbitrary vector at the point $x$. Index $q$ runs the
values $p$ and $f.$

In the $\Sigma $-space one can consider commutator of covariant derivatives $%
\tilde{\nabla}_{i}^{x^{\prime }}$ and $\tilde{\nabla}_{k^{\prime }}^{x}$
with respect to $x^{i}$ and $x^{\prime k}$ respectively. Calculations gives 
\begin{eqnarray}
&&\left( \tilde{\nabla}_{i}^{x^{\prime }}\tilde{\nabla}_{s^{\prime }}^{x}-%
\tilde{\nabla}_{s^{\prime }}^{x}\tilde{\nabla}_{i}^{x^{\prime }}\right)
T_{k\ldots l^{\prime }{\ldots }}^{j\ldots m^{\prime }{\ldots }}  \nonumber \\
&=&\tilde{F}_{ika^{\prime }s^{\prime }}\Sigma ^{ba^{\prime }}T_{b\ldots
l^{\prime }{\ldots }}^{j\ldots m^{\prime }{\ldots }}+\ldots -\Sigma
^{ja^{\prime }}\tilde{F}_{iba^{\prime }s^{\prime }}T_{k\ldots l^{\prime
}\ldots }^{b\ldots m^{\prime }{\ldots }}  \nonumber \\
&&-\ldots +\Sigma ^{am^{\prime }}\tilde{F}_{iab^{\prime }s^{\prime
}}T_{k\ldots l^{\prime }\ldots }^{j\ldots b^{\prime }\ldots }+\ldots -\tilde{%
F}_{ial^{\prime }s^{\prime }}\Sigma ^{ab^{\prime }}T_{k\ldots b^{\prime
}\ldots }^{j\ldots m^{\prime }\ldots }-\ldots  \label{f6.3}
\end{eqnarray}
where $T_{k\ldots l^{\prime }{\ldots }}^{j\ldots m^{\prime }{\ldots }}$ is a
two-point tensor. $\tilde{F}$-tensor, defined by the relation 
\begin{equation}
\tilde{F}_{ilk^{\prime }j^{\prime }}\equiv \Sigma _{iq^{\prime }}\tilde{%
\Gamma}_{k^{\prime }j^{\prime }||l}^{q^{\prime }}=\Sigma _{pj^{\prime }}%
\tilde{\Gamma}_{il||k^{\prime }}^{p}=\Sigma _{,ilj^{\prime }\parallel
k^{\prime }}=\Sigma _{,ilj^{\prime }k^{\prime }}-\Sigma _{,sj^{\prime
}k^{\prime }}\Sigma ^{sm^{\prime }}\Sigma _{,ilm^{\prime }}  \label{f6.4}
\end{equation}
is a two-point analog of the one-point curvature tensor $%
r_{slik}=g_{lp}r_{s.ik}^{.p}$. To test that the quantity (\ref{f6.4}) is a
tensor, let us represent it in one of two forms 
\begin{equation}
\tilde{F}_{ilk^{\prime }j^{\prime }}\equiv \Sigma _{pj^{\prime }}\tilde{%
\Gamma}_{il||k^{\prime }}^{p}=\Sigma _{pj^{\prime }}\left( \tilde{\Gamma}%
_{il}^{p}-\tilde{\gamma}_{\left( {\rm f}\right) }{}_{il}^{p}\right)
_{||k^{\prime }}  \label{f6.5}
\end{equation}
\begin{equation}
\tilde{F}_{ilk^{\prime }j^{\prime }}\equiv \Sigma _{iq^{\prime }}\tilde{%
\Gamma}_{k^{\prime }j^{\prime }||l}^{q^{\prime }}=\Sigma _{iq^{\prime
}}\left( \tilde{\Gamma}_{k^{\prime }j^{\prime }}^{q^{\prime }}-\tilde{\gamma}%
_{\left( {\rm p}\right) }{}_{k^{\prime }j^{\prime }}^{q^{\prime }}\right)
_{||l}  \label{f6.5a}
\end{equation}

As far as the difference 
\begin{equation}
\tilde{Q}_{\left( {\rm f}\right) }{}_{il}^{p}=\tilde{\gamma}_{\left( {\rm f}%
\right) }{}_{il}^{p}-\tilde{\Gamma}_{il}^{p}\qquad \tilde{Q}_{\left( {\rm p}%
\right) }{}_{k^{\prime }j^{\prime }}^{q^{\prime }}=\tilde{\gamma}_{\left( 
{\rm p}\right) }{}_{k^{\prime }j^{\prime }}^{q^{\prime }}-\tilde{\Gamma}%
_{k^{\prime }j^{\prime }}^{q^{\prime }}  \label{f6.6}
\end{equation}
of two Christoffel symbols is a tensor, it follows from (\ref{f6.5}) and (%
\ref{f6.6}) that $\tilde{F}_{ilk^{\prime }j^{\prime }}$ is a tensor. $\tilde{%
F}$-tensor can be presented as a result of covariant differentiation of the $%
\Sigma $-function. Indeed 
\begin{eqnarray}
\tilde{Q}_{\left( {\rm f}\right) }{}_{il}^{s} &=&-\Sigma ^{sj^{\prime
}}(\Sigma _{,ilj^{\prime }}-\tilde{\gamma}_{\left( {\rm f}\right)
}{}_{il}^{m}\Sigma _{mj^{\prime }})=-\Sigma ^{sj^{\prime }}\tilde{D}_{\left( 
{\rm f}\right) l}\Sigma _{ij^{\prime }}=-\Sigma ^{sj^{\prime }}\tilde{D}%
_{\left( {\rm f}\right) l}\tilde{D}_{\left( {\rm p}\right) j^{\prime
}}\Sigma _{i}  \nonumber \\
&=&-\Sigma ^{sj^{\prime }}\tilde{D}_{\left( {\rm p}\right) j^{\prime }}%
\tilde{D}_{\left( {\rm f}\right) l}\Sigma _{i}=-\Sigma ^{sj^{\prime }}\left( 
\tilde{D}_{\left( {\rm f}\right) l}\Sigma _{i}\right) _{||j^{\prime
}}=-\left( \Sigma ^{sj^{\prime }}\tilde{D}_{\left( {\rm f}\right) l}\Sigma
_{i}\right) _{||j^{\prime }}  \label{f6.7}
\end{eqnarray}
\begin{eqnarray}
\tilde{Q}_{\left( {\rm p}\right) }{}_{k^{\prime }j^{\prime }}^{q^{\prime }}
&=&-\Sigma ^{sq^{\prime }}\left( \Sigma _{,sk^{\prime }j^{\prime }}-\tilde{%
\gamma}_{\left( {\rm p}\right) }{}_{k^{\prime }j^{\prime }}^{m^{\prime
}}\Sigma _{sm^{\prime }}\right) =-\Sigma ^{sq^{\prime }}\tilde{D}_{\left( 
{\rm p}\right) j^{\prime }}\Sigma _{sk^{\prime }}  \nonumber \\
&=&-\Sigma ^{sq^{\prime }}\tilde{D}_{\left( {\rm p}\right) j^{\prime }}%
\tilde{D}_{\left( {\rm f}\right) s}\Sigma _{k^{\prime }}=-\Sigma
^{sq^{\prime }}\tilde{D}_{\left( {\rm f}\right) s}\tilde{D}_{\left( {\rm p}%
\right) j^{\prime }}\Sigma _{k^{\prime }}=-\left( \Sigma ^{sq^{\prime }}%
\tilde{D}_{\left( {\rm p}\right) j^{\prime }}\Sigma _{k^{\prime }}\right)
_{||s}  \label{f6.7a}
\end{eqnarray}
Then according to (\ref{f6.5}) -- (\ref{f6.7a}), one obtains

\begin{equation}
\tilde{F}_{ilk^{\prime }j^{\prime }}=\left( \tilde{D}_{\left( {\rm f}\right)
l}\Sigma _{i}\right) _{||j^{\prime }||k^{\prime }}=\left( \tilde{D}_{\left( 
{\rm p}\right) j^{\prime }}\Sigma _{k^{\prime }}\right) _{||i||l}
\label{f6.8}
\end{equation}

The commutator of covariant derivatives $\nabla _{i}^{x^{\prime }}$ and $%
\nabla _{k^{\prime }}^{x}$, connected with the symmetric component $G$ of
the world function, has the property 
\[
T_{k\ldots l^{\prime }{\ldots |s}^{\prime }|i}^{j\ldots m^{\prime }{\ldots }%
}-T_{k\ldots l^{\prime }{\ldots }|i|s^{\prime }}^{j\ldots m^{\prime }{\ldots 
}}=F_{ika^{\prime }s^{\prime }}G^{ba^{\prime }}T_{b\ldots l^{\prime }{\ldots 
}}^{j\ldots m^{\prime }{\ldots }}+\ldots -G^{ja^{\prime }}F_{iba^{\prime
}s^{\prime }}T_{k\ldots l^{\prime }\ldots }^{b\ldots m^{\prime }{\ldots }} 
\]
\begin{equation}
-\ldots +G^{am^{\prime }}F_{iab^{\prime }s^{\prime }}T_{k\ldots l^{\prime
}\ldots }^{j\ldots b^{\prime }\ldots }+\ldots -F_{ial^{\prime }s^{\prime
}}G^{ab^{\prime }}T_{k\ldots b^{\prime }\ldots }^{j\ldots m^{\prime }\ldots
}-\ldots  \label{f6.16}
\end{equation}
where the curvature $F$-tensor has the form 
\begin{equation}
F_{ilk^{\prime }j^{\prime }}\equiv G_{pj^{\prime }}\Gamma _{il|k^{\prime
}}^{p}=G_{pj^{\prime }}\left( \Gamma _{il}^{p}-\gamma _{il}^{p}\right)
_{|k^{\prime }}=G_{i;l|k^{\prime }|j^{\prime }}  \label{f6.17}
\end{equation}
Here ($;$) denotes the usual covariant derivative with the Christoffel
symbol $\gamma _{kl}^{i}$, and the $Q$-tensor is written as follows

\begin{equation}
Q_{kl}^{i}=\gamma _{kl}^{i}-\Gamma _{kl}^{i}=-G^{sj^{\prime
}}G_{i;l|j^{\prime }}  \label{f6.18}
\end{equation}

Let us discover a connection between the $F$-tensor at the coincidence limit 
$\left[ F_{ilk^{\prime }j^{\prime }}\right] $ and the curvature tensor $%
r_{s.ik}^{.l}$, constructed of the Christoffel symbols $\gamma _{si}^{l}$ by
means of formula (\ref{f6.2}) 
\begin{equation}
r_{s.ik}^{.l}=\gamma _{si,k}^{l}-\gamma _{sk,i}^{l}+\gamma _{si}^{p}\gamma
_{pk}^{l}-\gamma _{sk}^{p}\gamma _{pi}^{l}.  \label{f6.18a}
\end{equation}
Let us take into account that 
\begin{equation}
\left[ \left( \nabla _{i}^{x^{\prime }}+\nabla _{i^{\prime }}^{x}\right)
T_{sp^{\prime }}\left( x,x^{\prime }\right) \right]
_{x}=D_{i}t_{sp}=t_{sp;i},\qquad t_{sp}\equiv \left[ T_{sp^{\prime }}\left(
x,x^{\prime }\right) \right] _{x}  \label{f6.9}
\end{equation}
Then, using (\ref{f6.16}), (\ref{f6.17}), one obtains 
\begin{eqnarray}
&&\left[ \left( \nabla _{i}^{x^{\prime }}+\nabla _{i^{\prime }}^{x}\right)
\left( \nabla _{k}^{x^{\prime }}+\nabla _{k^{\prime }}^{x}\right)
T_{s}\left( x,x^{\prime }\right) -\left( \nabla _{k}^{x^{\prime }}+\nabla
_{k^{\prime }}^{x}\right) \left( \nabla _{i}^{x^{\prime }}+\nabla
_{i^{\prime }}^{x}\right) T_{s}\left( x,x^{\prime }\right) \right] _{x} 
\nonumber \\
&=&\left( D_{i}D_{k}-D_{k}D_{i}\right) t_{s}=\left[ \Gamma _{is|k^{\prime
}}^{a}-\Gamma _{ks|i^{\prime }}^{a}\right] _{x}t_{a}  \label{f6.10}
\end{eqnarray}
where $D_{i}$ is the usual covariant derivative in the Riemannian space with
the Christoffel symbol $\gamma _{kl}^{i}=\left[ \Gamma _{kl}^{i}\right] _{x}$

On the other hand the relation 
\begin{equation}
\left( D_{i}D_{k}-D_{k}D_{i}\right) t_{s}=r_{s.ik}^{.l}t_{l}  \label{f6.10a}
\end{equation}
takes place. Comparison of relations (\ref{f6.10}) and (\ref{f6.10a}) gives 
\begin{equation}
r_{s.ik}^{.l}=\left[ \Gamma _{is|k^{\prime }}^{l}-\Gamma _{ks|i^{\prime
}}^{l}\right] _{x}=\left[ \Gamma _{is,k^{\prime }}^{l}-\Gamma _{ks,i^{\prime
}}^{l}\right] _{x}=-g^{lp}f_{ispk}+g^{lp}f_{kspi}  \label{f6.13}
\end{equation}
where 
\begin{equation}
f_{ispk}\equiv \left[ F_{isp^{\prime }k^{\prime }}\right] _{x}=\left[
G_{lp^{\prime }}\Gamma _{is|k^{\prime }}^{l}\right] _{x}=-g_{lp}\left[
\Gamma _{is|k^{\prime }}^{l}\right] _{x}  \label{f6.14}
\end{equation}
According to (\ref{f6.17}) the one-point tensor $f_{ispk}$ is symmetric with
respect to transposition indices $i\leftrightarrow s$ and $p\leftrightarrow
k $ separately. 
\begin{equation}
f_{ispk}=f_{sipk},\qquad f_{ispk}=f_{iskp}  \label{f6.14a}
\end{equation}
Equation (\ref{f6.13}) can be written in the form 
\begin{equation}
g_{lp}r_{s.ik}^{.l}=-f_{ispk}+f_{kspi}  \label{f6.15}
\end{equation}
The metric tensor $g_{ik}$ is symmetric, and $f_{ispk}$ has the following
symmetry properties 
\begin{equation}
f_{ispk}=f_{sipk},\qquad f_{ispk}=f_{iskp},\qquad f_{ispk}=f_{pkis}
\label{f6.20}
\end{equation}
One can obtain connection of the type (\ref{f6.13}), (\ref{f6.14}) between
the $\tilde{F}$-tensor and the Riemannian -Christoffel curvature tensor in
the case of nonsymmetric T-geometry. Taking into account (\ref{f6.5}),
evident identity 
\begin{equation}
\frac{\partial }{\partial x^{k}}\left[ \tilde{\Gamma}_{il}^{p}\right]
_{x}\equiv \left[ \tilde{\Gamma}_{il,k}^{p}\right] _{x}+\left[ \tilde{\Gamma}%
_{il,k^{\prime }}^{p}\right] _{x}  \label{f6.21}
\end{equation}
and using relations (\ref{c5.15}), (\ref{c5.18}), one obtains 
\begin{equation}
\left[ \tilde{F}_{ilk^{\prime }j^{\prime }}\right] _{x}=-\tilde{g}_{pj}\left[
\tilde{\Gamma}_{il,k^{\prime }}^{p}\right] _{x}=-\tilde{g}_{pj}\left( \tilde{%
\gamma}_{\left( {\rm f}\right) }{}_{il,k}^{p}-\left[ \tilde{\Gamma}%
_{il,k}^{p}\right] _{x}\right)  \label{f6.22}
\end{equation}
Let us take into account identity 
\begin{equation}
\Sigma ^{pr^{\prime }}\Sigma _{,sr^{\prime }k}+\left( \Sigma ^{pr^{\prime
}}\right) _{,k}\Sigma _{sr^{\prime }}=0  \label{f6.23}
\end{equation}
obtained by differentiation of (\ref{c3.4}). Then using relations (\ref
{c5.15}), (\ref{c5.18}), one obtains from (\ref{f6.22}) 
\begin{eqnarray*}
\left[ \tilde{F}_{ilk^{\prime }j^{\prime }}\right] _{x} &=&-\tilde{g}_{pj}%
\tilde{\gamma}_{\left( {\rm f}\right) }{}_{il,k}^{p}+\tilde{g}_{pj}\left[
\Sigma ^{pq^{\prime }}\Sigma _{,ilkq^{\prime }}-\Sigma ^{pr^{\prime }}\Sigma
_{,sr^{\prime }k}\Sigma ^{sq^{\prime }}\Sigma _{,ilq^{\prime }}\right] _{x}
\\
&=&-\tilde{g}_{pj}\tilde{\gamma}_{\left( {\rm f}\right) }{}_{il,k}^{p}+%
\tilde{g}_{pj}\left[ -\tilde{g}^{pq}\Sigma _{,ilkq^{\prime }}-\tilde{\Gamma}%
_{sk}^{p}\tilde{\Gamma}_{il}^{s}\right] _{x}
\end{eqnarray*}
\begin{equation}
\tilde{f}_{ilkj}\equiv \left[ \tilde{F}_{ilk^{\prime }j^{\prime }}\right]
_{x}=-\left[ \Sigma _{,ilkj^{\prime }}\right] _{x}-\tilde{g}_{pj}\left( 
\tilde{\gamma}_{\left( {\rm f}\right) }{}_{il,k}^{p}+\tilde{\gamma}_{\left( 
{\rm f}\right) }{}_{il}^{s}\tilde{\gamma}_{\left( {\rm f}\right)
}{}_{sk}^{p}\right)  \label{f6.24}
\end{equation}
Alternating with respect to indices $k,l$, one obtains 
\begin{equation}
\tilde{f}_{ilkj}-\tilde{f}_{iklj}=\tilde{g}_{pj}\tilde{r}_{\left( {\rm f}%
\right) }{}_{i.kl}^{.p}  \label{f6.25}
\end{equation}
where $\tilde{r}_{\left( {\rm f}\right) }{}_{i.kl}^{.s}$ is the
Riemann-Christoffel curvature tensor, constructed on the base of the
Christoffel symbol $\gamma _{\left( {\rm f}\right) }{}_{ik}^{s}\ $ 
\begin{equation}
\tilde{r}_{\left( {\rm f}\right) }{}_{i.kl}^{.p}=\tilde{\gamma}_{\left( {\rm %
f}\right) }{}_{ik,l}^{p}-\tilde{\gamma}_{\left( {\rm f}\right)
}{}_{il,k}^{p}+\tilde{\gamma}_{\left( {\rm f}\right) }{}_{ik}^{s}\tilde{%
\gamma}_{\left( {\rm f}\right) }{}_{sl}^{p}-\tilde{\gamma}_{\left( {\rm f}%
\right) }{}_{il}^{p}\tilde{\gamma}_{\left( {\rm f}\right) }{}_{sk}^{k}
\label{f6.26}
\end{equation}

In the same way one can express $\tilde{f}_{ilkj}-\tilde{f}_{iklj}$ via the
Riemann-Christoffel curvature tensor $\tilde{r}_{\left( {\rm p}\right)
}{}_{i.kl}^{.s}$, constructed on the base of the Christoffel symbol $\tilde{%
\gamma}_{\left( {\rm p}\right) }{}_{ik}^{s}$ 
\begin{equation}
\tilde{f}_{ilkj}-\tilde{f}_{iklj}=\tilde{g}_{lp}\tilde{r}_{\left( {\rm p}%
\right) }{}_{k.ij}^{.p}  \label{f6.27}
\end{equation}
where 
\begin{equation}
\tilde{r}_{\left( {\rm p}\right) }{}_{i.kl}^{.p}=\tilde{\gamma}_{\left( {\rm %
p}\right) }{}_{ik,l}^{p}-\tilde{\gamma}_{\left( {\rm p}\right)
}{}_{il,k}^{p}+\tilde{\gamma}_{\left( {\rm p}\right) }{}_{ik}^{s}\tilde{%
\gamma}_{\left( {\rm p}\right) }{}_{sl}^{p}-\tilde{\gamma}_{\left( {\rm p}%
\right) }{}_{il}^{p}\tilde{\gamma}_{\left( {\rm p}\right) }{}_{sk}^{k}
\label{f6.28}
\end{equation}

To obtain representation (\ref{f6.27}), let us use another representation (%
\ref{f6.5a}) 
\[
\tilde{F}_{ilk^{\prime }j^{\prime }}=\Sigma _{iq^{\prime }}\tilde{\Gamma}%
_{k^{\prime }j^{\prime },l}^{q^{\prime }} 
\]
of the $\tilde{F}$-tensor, which differs from the representation (\ref{f6.5}%
) by a change $x\leftrightarrow x^{\prime }$. Producing the same operations (%
\ref{f6.22}) -- (\ref{f6.24}), one obtains (\ref{f6.27}) instead of (\ref
{f6.25}).

Note that relations (\ref{f6.25}) and (\ref{f6.27}) are different, because
the tensor $\tilde{g}_{ik}$ is not symmetric. In (\ref{f6.25}) summation is
produced over the first index, whereas in (\ref{f6.27}) it is produced over
the second index. In the symmetric T-geometry, when $\tilde{g}_{ik}$ is
symmetric, three expressions (\ref{f6.13}) (\ref{f6.25}) and (\ref{f6.27})
coincide.

There are two essentially different cases of asymmetric T-geometry:

1. Rough antisymmetry, when the field $a_{i}\left( x\right) \neq 0$. In this
case at small distances $x-x^{\prime }$ the field $a_{i}\left( x\right) $
dominates, and the world function is determined by the linear form 
\[
\Sigma \left( x,x^{\prime }\right) =a_{i}\left( x\right) \left(
x^{i}-x^{\prime i}\right) +... 
\]
In this case the antisymmetry is the main phenomenon at small distances.

2. Fine antisymmetry, when the field $a_{i}\left( x\right) \equiv 0$. In
this case the antisymmetric effects are described by the field $a_{ikl}$. At
small distances $x-x^{\prime }$ the symmetric structure dominates, and the
world function is determined by the quadratic form 
\[
\Sigma \left( x,x^{\prime }\right) =\frac{1}{2}g_{ik}\left( x\right) \left(
x^{i}-x^{\prime i}\right) \left( x^{k}-x^{\prime k}\right) +... 
\]
as in the symmetric T-geometry. In this case the antisymmetric effects may
be considered as corrections to gravitational effects. This corrections may
be essential at large distances $\xi ^{i}=x^{i}-x^{\prime i}$, when the form 
$\frac{1}{6}a_{ikl}\xi ^{i}\xi ^{k}\xi ^{l}$ becomes of the same order as
the form $\frac{1}{2}g_{ik}\xi ^{i}\xi ^{k}$,

The asymmetric T-geometry with fine antisymmetry is simpler, because it is
rather close to the usual symmetric T-geometry.

\section{Gradient lines on the manifold in the case of fine antisymmetry $%
a_{i}\equiv 0$}

Let us consider a one-dimensional line ${\cal L}_{({\rm f})}$, passing
through points $x^{\prime }$ and $x^{\prime \prime }.$ This line is defined
by the relations 
\begin{equation}
{\cal L}_{({\rm f})}:\;\;\;\Sigma _{,i^{\prime }}\left( x,x^{\prime }\right)
=\tau \Sigma _{,i^{\prime }}\left( x^{\prime \prime },x^{\prime }\right)
=\tau b_{i^{\prime }},\qquad i=0,1,...n  \label{c6.1}
\end{equation}
Let us suppose that $\det ||\Sigma _{,i^{\prime }k}\left( x,x^{\prime
}\right) ||\neq 0$. Then $n+1$ equations (\ref{c6.1}) can be resolved with
respect to $x$ in the form 
\begin{equation}
{\cal L}_{({\rm f})}:\;\;x^{i}=x^{i}\left( \tau \right) ,\qquad i=0,1,...n
\label{c6.2}
\end{equation}
where $\tau $ is a parameter along the line ${\cal L}_{({\rm f})}$. As it
follows from (\ref{c6.1}), this line passes through the point $x^{\prime }$
at $\tau =0$ and through the point $x^{\prime \prime }$ at $\tau =1$. Such a
line will be referred to as gradient line (curve) from the future. Let us
derive differential equation for the gradient curve ${\cal L}_{({\rm f})}$.

Differentiating (\ref{c6.1}) with respect to $\tau $, one obtains 
\begin{equation}
\Sigma _{,ki^{\prime }}\left( x,x^{\prime }\right) \frac{dx^{k}}{d\tau }%
=\Sigma _{,i^{\prime }}\left( x^{\prime \prime },x^{\prime }\right)
=b_{i^{\prime }},\qquad i=0,1,...n  \label{c6.3}
\end{equation}
Differentiating once more, one obtains 
\begin{equation}
\Sigma _{,ki^{\prime }}\left( x,x^{\prime }\right) \frac{d^{2}x^{k}}{d\tau
^{2}}+\Sigma _{,kli^{\prime }}\left( x,x^{\prime }\right) \frac{dx^{k}}{%
d\tau }\frac{dx^{l}}{d\tau }=0,\qquad i=0,1,...n  \label{c6.4}
\end{equation}
Using relation (\ref{c3.6}), one can write equations (\ref{c6.4}) in the
form 
\begin{equation}
\frac{d^{2}x^{i}}{d\tau ^{2}}+\tilde{\Gamma}_{kl}^{i}\left( x,x^{\prime
}\right) \frac{dx^{k}}{d\tau }\frac{dx^{l}}{d\tau }=0,\qquad i=0,1,...n
\label{c6.5}
\end{equation}
The equation (\ref{c6.5}) may be interpreted as an equation for a geodesic
in some $(n+1)$-dimensional Euclidean space with the Christoffel symbol $%
\tilde{\Gamma}_{kl}^{i}\left( x,x^{\prime }\right) $. This geodesic passes
through the points $x^{\prime }$ and $x^{\prime \prime }$.

Let the points $x^{\prime }$ and $x^{\prime \prime }$ be infinitesimally
close. Then equation (\ref{c6.5}) can be written in the form 
\begin{equation}
{\cal L}_{({\rm f})}:\;\;\frac{d^{2}x^{i}}{d\tau ^{2}}+\tilde{\gamma}%
_{\left( {\rm f}\right) }{}_{kl}^{i}\frac{dx^{k}}{d\tau }\frac{dx^{l}}{d\tau 
}=0,\qquad i=0,1,...n  \label{c6.6}
\end{equation}
where $\tilde{\gamma}_{\left( {\rm f}\right) }{}_{kl}^{i}=\tilde{\gamma}%
_{\left( {\rm f}\right) }{}_{kl}^{i}\left( x\right) =\left[ \tilde{\Gamma}%
_{kl}^{i}\right] _{x}$. Dividing the gradient line ${\cal L}_{({\rm f})}$
into infinitesimal segments and writing equations (\ref{c6.5}) in the form (%
\ref{c6.6}) on each segment, one obtains that the gradient line ${\cal L}_{(%
{\rm f})}$ is described by the equations (\ref{c6.6}) everywhere.

The equation (\ref{c6.6}) does not contain a reference to the point $%
x^{\prime }$, and any gradient line (\ref{c6.1}), (\ref{c6.2}) is to satisfy
this equation.

In the case of fine antisymmetry, when $a_{i}\equiv 0$, the equation (\ref
{c6.6}) can be written in other form. Using relations (\ref{c5.18}), and
taking into account that $a_{i}\equiv 0$, one obtains instead of (\ref{c6.6}%
) 
\begin{equation}
{\cal L}_{({\rm f})}:\qquad \frac{d^{2}x^{i}}{d\tau ^{2}}+\left( \gamma
_{kl}^{i}+\beta _{kl}^{i}\right) \frac{dx^{k}}{d\tau }\frac{dx^{l}}{d\tau }%
=0,\qquad a_{i}\equiv 0  \label{c6.6b}
\end{equation}
where 
\begin{eqnarray}
\gamma _{kl}^{i} &=&\gamma _{kl}^{i}\left( x\right) =\frac{1}{2}g^{si}\left(
g_{ks,l}+g_{sl,k}-g_{lk,s}\right) ,  \label{c6.6c} \\
\beta _{kl}^{i} &=&\beta _{kl}^{i}\left( x\right) =g^{si}a_{kls}
\label{c6.6d}
\end{eqnarray}
If $a_{ikl}=0$, the equations (\ref{c6.6b}) may be considered to be the
equations for a geodesic in a Riemannian space with the metric tensor $%
g_{ik} $.

In the case of rough antisymmetry, when $a_{i}\neq 0$, equations (\ref{c6.5}%
), (\ref{c6.6}) also describe a gradient line, but this line cannot be
described by equation (\ref{c6.1}). The fact is that at $a_{i}\neq 0,$ rhs
of (\ref{c6.1}) vanishes at $\tau \rightarrow 0,$ whereas lhs of (\ref{c6.1}%
) does not vanish, in general. Equation (\ref{c6.1}) stops to be valid at $%
\tau \rightarrow 0$.

Now let us consider another type of gradient line ${\cal L}_{({\rm p})}$,
passing through the points $x$ and $x^{\prime \prime }$. Let the gradient
line ${\cal L}_{({\rm p})}$ be described by the equations (It is supposed
again that $a_{i}\equiv 0$) 
\begin{equation}
{\cal L}_{({\rm p})}:\qquad \Sigma _{,i}\left( x,x^{\prime }\right) =\tau
\Sigma _{,i}\left( x,x^{\prime \prime }\right) =\tau b_{i},\qquad i=0,1,...n
\label{c6.7}
\end{equation}
which determine 
\begin{equation}
{\cal L}_{({\rm p})}:\qquad x^{\prime i}=x^{\prime i}\left( \tau \right)
,\qquad i=0,1,...n  \label{c6.8}
\end{equation}
Equation (\ref{c6.7}) distinguishes from the equation (\ref{c6.1}) only in
transposition of the first and second arguments of the world function $%
\Sigma \left( x,x^{\prime }\right) .$ The gradient line ${\cal L}_{({\rm p}%
)} $, determined by the relation (\ref{c6.7}), may be referred to as the
gradient line from the past. Manipulating with the equation (\ref{c6.7}) in
the same way as with (\ref{c6.1}), one obtains instead of (\ref{c6.6}) 
\begin{equation}
{\cal L}_{({\rm p})}:\qquad \frac{d^{2}x^{i}}{d\tau ^{2}}+\tilde{\gamma}%
_{\left( {\rm p}\right) }{}_{kl}^{i}\frac{dx^{k}}{d\tau }\frac{dx^{l}}{d\tau 
}=0,\qquad i=0,1,...n  \label{c6.8a}
\end{equation}

In the case of fine antisymmetry, when $a_{i}\equiv 0$, the equation (\ref
{c6.8a}) can be written in other form. Using relation (\ref{c5.18b}), and
taking into account that $a_{i}\equiv 0$, one obtains instead of (\ref{c6.8a}%
) 
\begin{equation}
{\cal L}_{({\rm p})}:\qquad \frac{d^{2}x^{i}}{d\tau ^{2}}+\left( \gamma
_{kl}^{i}-\beta _{kl}^{i}\right) \frac{dx^{k}}{d\tau }\frac{dx^{l}}{d\tau }%
=0,\qquad a_{i}\equiv 0  \label{c6.9a}
\end{equation}
where $\gamma _{kl}^{i}$ and $\beta _{kl}^{i}$ are defined by the relations (%
\ref{c6.6c}), (\ref{c6.6d}).

In the case of symmetric T-geometry, when $a_{ikl}=0$ and $\beta _{kl}^{i}=0$%
, differential equations (\ref{c6.6b}) and (\ref{c6.9a}) respectively for
gradient line ${\cal L}_{({\rm f})}$ and for gradient line ${\cal L}_{({\rm p%
})}$ coincide.

In the case of asymmetric T-geometry the quantities $\left[ \Gamma _{kl}^{i}%
\right] _{x}$ and $\left[ \Gamma _{k^{\prime }l^{\prime }}^{i^{\prime }}%
\right] _{x}$ do not coincide, in general. In this case the equations (\ref
{c6.6}) and (\ref{c6.9a}) determine, in general, different gradient curves,
passing through the same points $x^{\prime }$ and $x^{\prime \prime }$.
Differential equations for the gradient curves ${\cal L}_{({\rm p})}$ and $%
{\cal L}_{({\rm f})}$ differ in the sign of the ''antisymmetric force'' 
\begin{equation}
\beta _{kl}^{i}\frac{dx^{k}}{d\tau }\frac{dx^{l}}{d\tau }=g^{si}\left( -%
\frac{1}{2}a_{k,ls}-\frac{1}{2}a_{l,ks}+a_{kls}\right) \frac{dx^{k}}{d\tau }%
\frac{dx^{l}}{d\tau }  \label{c6.10}
\end{equation}

Finally, one can introduce the neutral gradient line ${\cal L}_{({\rm n})}$,
defining it by the relations 
\begin{equation}
{\cal L}_{({\rm n})}:\;\;G_{,i^{\prime }}\left( x,x^{\prime }\right) =\tau
G_{,i^{\prime }}\left( x^{\prime \prime },x^{\prime }\right) =\tau
b_{i^{\prime }},\qquad i=0,1,...n  \label{c6.10a}
\end{equation}
which determine 
\begin{equation}
{\cal L}_{({\rm n})}:\;\;x^{i}=x^{i}\left( \tau \right) ,\qquad i=0,1,...n
\label{c6.10b}
\end{equation}
Equation (\ref{c6.10a}) distinguishes from the equation (\ref{c6.1}) only in
replacement of the world function $\Sigma \left( x,x^{\prime }\right) $ by
its symmetric component $G\left( x,x^{\prime }\right) .$ Manipulating with
the equation (\ref{c6.10a}) in the same way as with (\ref{c6.1}), one
obtains instead of (\ref{c6.6b}) 
\begin{equation}
{\cal L}_{({\rm n})}:\qquad \frac{d^{2}x^{i}}{d\tau ^{2}}+\gamma _{kl}^{i}%
\frac{dx^{k}}{d\tau }\frac{dx^{l}}{d\tau }=0  \label{c6.10c}
\end{equation}
where the ''antisymmetric force'' is absent.

The gradient lines (\ref{c6.1}) and (\ref{c6.7}) are insensitive to
transformation of the world function of the form 
\begin{equation}
\Sigma \rightarrow \tilde{\Sigma}=f\left( \Sigma \right) ,\qquad \left|
f^{\prime }\left( \Sigma \right) \right| >0  \label{c6.10d}
\end{equation}
where $f$ is an arbitrary function, because for determination of the
gradient line only direction of the gradient $\Sigma _{i}$ or $\Sigma
_{i^{\prime }}$ is important, but not its module. Indeed, after substitution
of $\tilde{\Sigma}$ from (\ref{c6.10d}) in (\ref{c6.1}) one obtains the
equation 
\begin{equation}
\Sigma _{,i^{\prime }}\left( x,x^{\prime }\right) =\tau ^{\prime }\Sigma
_{,i^{\prime }}\left( x^{\prime \prime },x^{\prime }\right) ,\qquad
i=0,1,...n,\qquad \tau ^{\prime }=\tau \frac{f^{\prime }\left( \Sigma \left(
x^{\prime \prime },x^{\prime }\right) \right) }{f^{\prime }\left( \Sigma
\left( x,x^{\prime }\right) \right) }  \label{c6.10e}
\end{equation}
which describes the same gradient line, but with another parametrization.

\section{$\protect\sigma $-immanent description of gradient lines.}

The particle motion in the space-time geometry should be described $\sigma $%
-immanently, i.e. without a reference to a coordinate system and to a
manifold. The motion of a particle is described by its world tube, which is
a broken tube 
\begin{equation}
{\cal T}_{\left( {\rm f}\right) {\rm br}}=\bigcup_{i}{\cal T}_{({\rm f}%
)[P_{i}P_{i+1}]},  \label{c6.11}
\end{equation}
consisting of segments ${\cal T}_{({\rm f})[P_{i}P_{i+1}]},$ defined by the
relation (\ref{f2.6}). Each segment ${\cal T}_{({\rm f})[P_{i}P_{i+1}]}$ is
associated with the vector $\overrightarrow{P_{i}P_{i+1}}$ of the particle
momentum, having the length $\mu $%
\begin{equation}
\mu =|\overrightarrow{P_{i}P_{i+1}}|=\sqrt{2G\left( P_{i},P_{i+1}\right) }
\label{c6.12}
\end{equation}
The length $\mu $ is supposed to be similar for all segments. It is
associated with the particle mass $m$%
\begin{equation}
m=b\mu ,\qquad b=\text{const}  \label{c6.14}
\end{equation}
where $b$ is some universal constant. The broken tube (\ref{c6.11}), (\ref
{c6.12}) will be referred to as the future broken tube.

In the asymmetric T-geometry there is the past broken tube 
\begin{equation}
{\cal T}_{({\rm p}){\rm br}}=\bigcup_{i}{\cal T}_{({\rm p})[P_{i}P_{i+1}]},
\label{c6.15}
\end{equation}
consisting of segments ${\cal T}_{({\rm p})[P_{i}P_{i+1}]}$, defined by the
relation (\ref{f2.6}) at $q=p$.

There is a neutral broken tube, defined by the relation 
\begin{equation}
{\cal T}_{({\rm n}){\rm br}}=\bigcup_{i}{\cal T}_{({\rm n})[P_{i}P_{i+1}]},
\label{c6.16a}
\end{equation}
and the relation (\ref{f2.6}) at $q=n$. In the symmetric T-geometry
definitions of all broken tubes ${\cal T}_{({\rm f}){\rm br}}$, ${\cal T}_{(%
{\rm p}){\rm br}}$, ${\cal T}_{({\rm n}){\rm br}}$ coincide.

The gradient line ${\cal L}_{\left( {\rm f}\right) }$ may be defined as a
broken tube (\ref{c6.11}) with the segment length $\mu \rightarrow 0.$ If
segment $\overrightarrow{P_{i}P_{i+1}}$ is fixed, the next segment $%
\overrightarrow{P_{i+1}P_{i+2}}$ is defined by its end point $P_{i+2}$,
which provides extremum to the expression $\sqrt{2\Sigma \left(
P_{i},P_{i+2}\right) }$ under additional constraints 
\[
\sqrt{2\Sigma \left( P_{i},P_{i+1}\right) }=\mu ,\qquad \sqrt{2\Sigma \left(
P_{i+1},P_{i+2}\right) }=\mu . 
\]

The gradient line ${\cal L}_{\left( {\rm p}\right) }$ may be defined as a
broken tube (\ref{c6.15}) with the segment length $\mu \rightarrow 0.$ If
segment $\overrightarrow{P_{i}P_{i+1}}$ is fixed, the next segment $%
\overrightarrow{P_{i+1}P_{i+2}}$ is defined by its end point $P_{i+2}$,
which provides extremum to the expression $\sqrt{2\Sigma \left(
P_{i+2},P_{i}\right) }$ under additional constraints 
\[
\sqrt{2\Sigma \left( P_{i+1},P_{i}\right) }=\mu ,\qquad \sqrt{2\Sigma \left(
P_{i+2},P_{i+1}\right) }=\mu . 
\]

If extremum is achieved at one point $P_{i+2}$, then position of $P_{i+2}$
is determined uniquely. In this case the broken tubes ${\cal T}_{({\rm f})%
{\rm br}}$ and ${\cal T}_{({\rm p}){\rm br}}$ turn to the gradient lines at $%
\mu \rightarrow 0$.

Thus, using this procedure, one can define gradient lines $\sigma $%
-immanently. In the space-time geometry a gradient line is associated with a
world line of a free massless particle, whereas the broken tubes with fixed
length $\left| {\cal T}_{[P_{i}P_{i+1}]}\right| =\mu $ of segments ${\cal T}%
_{[P_{i}P_{i+1}]}$ are associated with world tubes of massive particles.

If the particle is free, adjacent segments are parallel, i.e. for the broken
tube (\ref{c6.11}) $\overrightarrow{P_{i}P_{i+1}}\uparrow \uparrow _{({\rm f}%
)}\overrightarrow{P_{i+1}P_{i+2}}$,\ \ $i=...0,1,2,...$and 
\begin{equation}
\left| \overrightarrow{P_{i}P_{i+1}}\right| \cdot \left| \overrightarrow{%
P_{i+1}P_{i+2}}\right| -\left( \overrightarrow{P_{i}P_{i+1}}.\overrightarrow{%
P_{i+1}P_{i+2}}\right) =0  \label{c6.17}
\end{equation}
For the broken tubes ${\cal T}_{({\rm p}){\rm br}}$, ${\cal T}_{({\rm n})%
{\rm br}}$, defined respectively by the relations (\ref{c6.15}), and (\ref
{c6.16a}), the condition of parallelism has the form $\overrightarrow{%
P_{i}P_{i+1}}\uparrow \uparrow _{({\rm p})}\overrightarrow{P_{i+1}P_{i+2}}$%
,\ \ $i=...0,1,2,...$ or 
\begin{equation}
\left| \overrightarrow{P_{i}P_{i+1}}\right| \left| \overrightarrow{%
P_{i+1}P_{i+2}}\right| -\left( \overrightarrow{P_{i+1}P_{i+2}}.%
\overrightarrow{P_{i}P_{i+1}}\right) =0  \label{c6.18}
\end{equation}
and $\overrightarrow{P_{i}P_{i+1}}\uparrow \uparrow _{({\rm n})}%
\overrightarrow{P_{i+1}P_{i+2}}$,\ \ $i=...0,1,2,...$ or 
\begin{equation}
\left| \overrightarrow{P_{i}P_{i+1}}\right| \left| \overrightarrow{%
P_{i+1}P_{i+2}}\right| -\sqrt{\left( \overrightarrow{P_{i+1}P_{i+2}}.%
\overrightarrow{P_{i}P_{i+1}}\right) \left( \overrightarrow{P_{i}P_{i+1}}.%
\overrightarrow{P_{i+1}P_{i+2}}\right) }=0  \label{c6.19}
\end{equation}

In general segments ${\cal T}_{({\rm f})[P_{i}P_{i+1}]}$ and ${\cal T}_{(%
{\rm p})[P_{i}P_{i+1}]}$ are inswept surfaces. Under some conditions they
degenerate into a one-dimensional line, and the first order tubes ${\cal T}%
_{({\rm f}){\rm br}}$ and ${\cal T}_{({\rm p}){\rm br}}$ degenerate
correspondently into gradient lines ${\cal L}_{\left( {\rm f}\right) }$ and $%
{\cal L}_{\left( {\rm p}\right) }.$ Let us investigate, when this fact takes
place.

\section{Conditions of degeneration of the neutral first order tube}

In general, the asymmetric T-geometry is nondegenerate geometry, even if it
is given on $n$-dimensional manifold ${\cal M}_{n}$. In general case the
first order tube ${\cal T}_{({\rm n}){\bf x}^{\prime }{\bf x}^{\prime \prime
}}$, passing through points $x^{\prime }$ and $x^{\prime \prime }$, does not
coincide with gradient line ${\cal L}_{\left( {\rm f}\right) }$ or ${\cal L}%
_{\left( {\rm p}\right) },$ passing through the points $x^{\prime }$ and $%
x^{\prime \prime }$ and defined by the relations (\ref{c6.1}) and (\ref{c6.7}%
) respectively.

Let us investigate, under what conditions the neutral first order tube $%
{\cal T}_{({\rm n}){\bf x}^{\prime }{\bf x}^{\prime \prime }}$ degenerates
into gradient line ${\cal L}_{\left( {\rm f}\right) }$ or ${\cal L}_{\left( 
{\rm p}\right) }$. The first order tube ${\cal T}_{({\rm n}){\bf P}_{0}{\bf P%
}_{1}}$, passing through the points $P_{0}=\left\{ x\right\} $, $%
\;P_{1}=\left\{ x^{\prime }\right\} $ is defined by the relation 
\begin{equation}
F_{2}\left( P_{0},P_{1},R\right) =\left| 
\begin{array}{cc}
\left( \overrightarrow{P_{0}P_{1}}.\overrightarrow{P_{0}P_{1}}\right) & 
\left( \overrightarrow{P_{0}P_{1}}.\overrightarrow{P_{0}R}\right) \\ 
\left( \overrightarrow{P_{0}R}.\overrightarrow{P_{0}P_{1}}\right) & \left( 
\overrightarrow{P_{0}R}.\overrightarrow{P_{0}R}\right)
\end{array}
\right| =0  \label{c7.1}
\end{equation}
where $R=\left\{ x^{\prime }+dx^{\prime }\right\} $ is a running point on
the tube ${\cal T}_{P_{0}P_{1}}$.

One has the following expansion for the scalar $\Sigma $-products $\left( 
\overrightarrow{P_{0}P_{1}}.\overrightarrow{P_{0}R}\right) $ and $\left( 
\overrightarrow{P_{0}R}.\overrightarrow{P_{0}P_{1}}\right) $%
\[
\left( \overrightarrow{P_{0}P_{1}}.\overrightarrow{P_{0}R}\right) =\Sigma
\left( P_{1},P_{0}\right) +\Sigma \left( P_{0},R\right) -\Sigma \left(
P_{1},R\right) 
\]
\begin{eqnarray}
&=&\Sigma \left( x^{\prime },x\right) +\Sigma \left( x,x^{\prime
}+dx^{\prime }\right) -\Sigma \left( x^{\prime },x^{\prime }+dx^{\prime
}\right)  \nonumber \\
&=&2G+\left( \Sigma _{,i^{\prime }}-\left[ \Sigma _{,i^{\prime }}\right]
_{x^{\prime }}\right) dx^{\prime i^{\prime }}+\frac{1}{2}\left( \Sigma
_{,i^{\prime }k^{\prime }}-\left[ \Sigma _{,i^{\prime }k^{\prime }}\right]
_{x^{\prime }}\right) dx^{\prime i^{\prime }}dx^{\prime k^{\prime }}
\label{c7.2}
\end{eqnarray}
\[
\left( \overrightarrow{P_{0}R}.\overrightarrow{P_{0}P_{1}}\right) =\Sigma
\left( P_{0},P_{1}\right) +\Sigma \left( R,P_{0}\right) -\Sigma \left(
R,P_{1}\right) 
\]
\begin{eqnarray}
&=&\Sigma \left( x,x^{\prime }\right) +\Sigma \left( x^{\prime }+dx^{\prime
},x\right) -\Sigma \left( x^{\prime }+dx^{\prime },x^{\prime }\right) 
\nonumber \\
&=&2G+\left( G_{,i^{\prime }}-A_{,i^{\prime }}-\left[ \Sigma _{,i}\right]
_{x^{\prime }}\right) dx^{\prime i^{\prime }}+\frac{1}{2}\left(
G_{,i^{\prime }k^{\prime }}-A_{,i^{\prime }k^{\prime }}-\left[ \Sigma _{,ik}%
\right] _{x^{\prime }}\right) dx^{\prime i}dx^{\prime k}  \label{c7.3}
\end{eqnarray}
where unprimed indices are associated with the first argument and the primed
ones with the second argument of the $\Sigma $-function. 
\begin{equation}
\left( \overrightarrow{P_{0}P_{1}}.\overrightarrow{P_{0}P_{1}}\right) =2G
\label{c7.4}
\end{equation}
\begin{equation}
\left( \overrightarrow{P_{0}R}.\overrightarrow{P_{0}R}\right) =2G\left(
x,x^{\prime }+dx^{\prime }\right) =2G+2G_{,i^{\prime }}dx^{\prime i^{\prime
}}+G_{,i^{\prime }k^{\prime }}dx^{\prime i^{\prime }}dx^{\prime k^{\prime }}
\label{c7.5}
\end{equation}
Using relations (\ref{c7.2}) -- (\ref{c7.5}) and (\ref{c5.14}), (\ref{c5.15}%
), one reduces equation (\ref{c7.1}) to the form 
\begin{equation}
\left( \left( G_{,i^{\prime }}-A_{,i^{\prime }}-a_{i}\left( x^{\prime
}\right) \right) dx^{\prime i^{\prime }}\right) \left( \left( G_{,l^{\prime
}}+A_{,l^{\prime }}+a_{l}\left( x^{\prime }\right) \right) dx^{\prime
^{\prime }}\right) =2Gg_{l^{\prime }m^{\prime }}\left( x^{\prime }\right)
dx^{\prime l^{\prime }}dx^{\prime m^{\prime }}  \label{c7.6}
\end{equation}

We suppose that the world function is such, that the tube ${\cal T}_{{\bf x}%
^{\prime }{\bf x}}$ degenerates to a line. Then the solution of (\ref{c7.6})
has either the form 
\begin{equation}
dx^{\prime k^{\prime }}=g^{k^{\prime }l^{\prime }}\left( x^{\prime }\right)
\left( G_{,i^{\prime }}-A_{,i^{\prime }}-a_{i^{\prime }}\left( x^{\prime
}\right) \right) d\tau  \label{c7.9}
\end{equation}
or the form 
\begin{equation}
dx^{\prime k^{\prime }}=g^{k^{\prime }l^{\prime }}\left( x^{\prime }\right)
\left( G_{,k^{\prime }}+A_{,k^{\prime }}+a_{k^{\prime }}\left( x^{\prime
}\right) \right) d\tau  \label{c7.10}
\end{equation}
where $d\tau $ is an infinitesimal parameter. The relations (\ref{c7.9}), (%
\ref{c7.10}) describe a one-dimensional lines in vicinity of the point $%
x^{\prime }$. Both expressions (\ref{c7.9}), (\ref{c7.10}) are solutions of
the equation (\ref{c7.6}), provided the relation 
\begin{equation}
\left( G_{,i^{\prime }}-A_{,i^{\prime }}-a_{i^{\prime }}\left( x^{\prime
}\right) \right) g^{i^{\prime }k^{\prime }}\left( x^{\prime }\right) \left(
G_{,k^{\prime }}+A_{,k^{\prime }}+a_{k^{\prime }}\left( x^{\prime }\right)
\right) =2G  \label{c7.11}
\end{equation}
is fulfilled.

There is only one solution, provided solutions (\ref{c7.9}) and (\ref{c7.10}%
) coincide. It means that 
\begin{equation}
A_{,k^{\prime }}+a_{k^{\prime }}\left( x^{\prime }\right) =0  \label{c7.11a}
\end{equation}
and the condition (\ref{c7.11}) transforms to the equation 
\begin{equation}
G_{,i^{\prime }}g^{i^{\prime }k^{\prime }}\left( x^{\prime }\right)
G_{,k^{\prime }}=2G  \label{c7.12}
\end{equation}
This is well known equation for the world function of a Riemannian space 
\cite{S60}. The Riemannian geometry is locally degenerate in the sense of
definition \ref{d3.14}, and the equation (\ref{c7.12}) describes this
property of Riemannian geometry.

According to expansion (\ref{c5.5}) the condition (\ref{c7.11a}) may take
place, only if 
\begin{equation}
A\left( x,x^{\prime }\right) =a_{i}\left( x^{i}-x^{\prime i}\right) ,\qquad
a_{i}=\text{const}  \label{c7.13}
\end{equation}
In this case the quantities (\ref{f2.2ee}) vanish, i.e. 
\begin{equation}
\eta _{{\rm f}}=A\left( x,x^{\prime }\right) +A\left( x^{\prime },y\right)
+A\left( y,x\right) =0,\qquad \forall x,x^{\prime },y\in {\cal M}_{n}
\label{c7.14}
\end{equation}
and the first order tubes are similar in symmetric and nonsymmetric
geometries.

For the case of the past first order tube and the future one the conditions
of degeneration are also rather rigid. In this case instead of (\ref{c7.6})
one obtains two conditions 
\begin{equation}
4G\left( A_{,i^{\prime }}+a_{i^{\prime }}\right) dx^{\prime i^{\prime }}=0
\label{c7.15}
\end{equation}
\begin{equation}
\left( 2G\left( A_{,i^{\prime }k^{\prime }}-\left[ \Sigma _{,i^{\prime
}k^{\prime }}\right] _{x^{\prime }}\right) -G_{,i^{\prime }}G_{,k^{\prime
}}\right) dx^{\prime i^{\prime }}dx^{\prime k^{\prime }}=0  \label{c7.16}
\end{equation}
In the case (\ref{c7.13}) the condition (\ref{c7.15}) is fulfilled, and (\ref
{c7.16}) reduces to (\ref{c7.12}).

If the first order tube is nondegenerate in symmetric T-geometry, it cannot
degenerate after addition of antisymmetric component, because the local
degeneration condition (\ref{c7.12}) remains to be not fulfilled.

Thus, practically any antisymmetric component of the world function destroys
degeneration of the neutral first order tube. If one connects quantum
effects with the first order tube degeneration \cite{R91}, one concludes
that the possible asymmetry of the space-time geometry is connected with
quantum effects.

\section{Examples of the first order tubes}

To imagine the possible corollaries of asymmetry in T-geometry, let us
construct the first order tube ${\cal T}_{P_{0}P_{1}}$ in the $\Sigma $%
-space. Let us consider $\Sigma $-space on the $4$-dimensional manifold with
the world function 
\begin{eqnarray}
\Sigma \left( x,x^{\prime }\right) &=&a_{i}\xi ^{i}+\frac{1}{2}g_{ki}\xi
^{i}\xi ^{k},\qquad a_{i}=b_{i}\left( 1+\alpha f\left( \xi ^{2}\right)
\right) ,  \label{cc3.1} \\
\xi ^{i} &=&x^{i}-x^{\prime i},\qquad \xi ^{2}\equiv g_{ki}\xi ^{i}\xi ^{k},%
\text{\qquad }\alpha ,b_{i},g_{ik}=\text{const,}  \nonumber
\end{eqnarray}
where $f$ is some function of $\xi ^{2}$ and summation is made over
repeating indices from $0$ to $3$. One can interpret the relation (\ref
{cc3.1}) as an Euclidean space with a linear structure $a_{i}\xi ^{i}$ given
on it. Such a $\Sigma $-space is uniform, but not isotropic, because there
is a vector $a_{i}$, describing some preferable direction in the $\Sigma $%
-space.

In the given case the characteristic quantity (\ref{f2.2ee}) has the form 
\begin{equation}
\eta _{{\rm f}}=\alpha \left( -b_{i}x^{i}f\left( x^{2}\right)
+b_{i}x^{\prime i}f\left( x^{\prime 2}\right) +b_{i}\xi ^{i}f\left( \xi
^{2}\right) \right)  \label{cc3.1a}
\end{equation}
It does not depend on the constant component of the vector $a_{i}.$ Then
according to (\ref{f2.3}) - (\ref{f2.1}) the shape of the tube does not
depend on the constant component of the vector $a_{i}$. If $\alpha =0$ and $%
a_{i}=$const, shape of all first order tubes is the same, as in the case of
symmetric T-geometry, when $a_{i}=0$. In other words, the shape of the first
order tubes is insensitive to the space-time anisotropy, described by the
vector field $a_{i}=$const. We omit the constant component of the function $%
f $ and consider the cases, when its variable part has the form 
\begin{eqnarray}
1 &:&\quad f\left( \xi ^{2}\right) =\xi ^{2},\qquad 2:\quad f\left( \xi
^{2}\right) =\frac{1}{1+\beta \xi ^{2}},\;\;\;\;\beta =\text{const}
\label{cc3.1b} \\
\xi ^{2} &\equiv &g_{ik}\xi ^{i}\xi ^{k},\qquad \xi ^{i}\equiv
x^{i}-x^{\prime i}  \nonumber
\end{eqnarray}
In the first case the antisymmetric structure is essential at large
distances $\xi =x-x^{\prime }$. In the second case the antisymmetric
structure vanishes at large $\xi $.

Now let us construct the first order neutral tube ${\cal T}_{x^{\prime }y}$,
determined by two points $x^{\prime }$ and $y$. For simplicity the
coordinate system is chosen in such a way, that $x^{\prime }=0$. The
equation (\ref{f2.3}) determining the shape of the tube has the form 
\begin{equation}
\left| 
\begin{array}{cc}
2G\left( x^{\prime },y\right) & \left( {\bf x}^{\prime }{\bf y.x}^{\prime }%
{\bf x}\right) \\ 
\left( {\bf x}^{\prime }{\bf x}^{\prime }{\bf .x}^{\prime }{\bf y}\right) & 
2G\left( x^{\prime },x\right)
\end{array}
\right| =4G\left( x^{\prime },y\right) G\left( x^{\prime },x\right) -\left( 
{\bf x}^{\prime }{\bf y.x}^{\prime }{\bf x}\right) \left( {\bf x}^{\prime }%
{\bf x.x}^{\prime }{\bf y}\right) =0  \label{cc3.2}
\end{equation}
In the first case, when $f\left( \xi ^{2}\right) =\xi ^{2}$, calculation
gives for (\ref{cc3.2}) 
\begin{equation}
\left( x_{i}y^{i}\right) ^{2}-x^{2}y^{2}=\alpha ^{2}\left[ \left(
x_{i}y^{i}\right) \left( -2\left( b_{k}y^{k}\right) +2\left(
b_{k}x^{k}\right) \right) -\left( b_{i}x^{i}\right) y^{2}+\left(
b_{i}y^{i}\right) x^{2}\right] ^{2}  \label{cc3.4}
\end{equation}
where $x^{2}=x^{i}x_{i}$, $y^{2}=y_{i}y^{i}$. In the case, when the metric
tensor $g_{ik}$ is the metric tensor of the proper Euclidean space, $%
x^{2}y^{2}\geq \left( x_{i}y^{i}\right) ^{2}$, and the equation (\ref{cc3.4}%
) has an interesting solution, only if $\alpha =0$. Then 
\begin{equation}
x^{2}y^{2}=\left( x_{i}y^{i}\right) ^{2},\;\;\;\;{\cal T}_{0y}=\left\{
x\left| \bigwedge_{i=0}^{i=3}x^{i}=y^{i}\tau \right. \right\} ,\qquad \alpha
=0  \label{cc3.5}
\end{equation}
In the case $\alpha \neq 0$, the first order tube ${\cal T}_{0y}$
degenerates to the set of basic points $\left\{ 0,y\right\} $, because
substitution of $x^{i}=y^{i}\tau $ in the square bracket in (\ref{cc3.4})
shows that the bracket vanishes only at $\tau =0$ or $\tau =1$. Thus, in the
case of proper Euclidean metric tensor $g_{ik}$ the first order tube shape
does not depend on $a_{i},$ provided $a_{i}=$const.

Let us consider a more interesting case, when the metric tensor $g_{ik}$ of $%
\Sigma $-space is the Minkowski one. Then $x^{2}y^{2}\leq \left(
x_{i}y^{i}\right) ^{2},$ provided the ${\bf 0y}$ is timelike $\left| {\bf 0y}%
\right| ^{2}=2G\left( 0,y\right) >0$. In this case the equation (\ref{cc3.4}%
) has the solution (\ref{cc3.5}), if $\alpha =0$.

If $\alpha \neq 0,$ let us introduce new variables $z,w$ by means of
relations 
\begin{equation}
x^{i}=\zeta ^{i}z\left| y\right| ,\qquad y^{i}=\eta ^{i}\left| y\right|
,\qquad \left| y\right| =\sqrt{y_{i}y^{i}},\qquad w=\left( \zeta ^{i}\eta
_{i}\right) ,  \label{cc3.5b}
\end{equation}
where $\eta ^{i}$ and $\zeta ^{i}$ are unit vectors 
\begin{equation}
\zeta ^{i}\zeta _{i}=\nu =\pm 1,\qquad \eta _{i}\eta ^{i}=1  \label{cc3.5d}
\end{equation}
Let us suppose that the vector $b_{i},$ which determines the antisymmetric
structure, is a unit timelike vector, and the vector $y^{i}$, determining
the tube, is chosen in such a way, that 
\begin{equation}
b_{i}=\varepsilon \eta _{i},\;\;\varepsilon =\pm 1  \label{cc3.5c}
\end{equation}
Then equation (\ref{cc3.4}) takes the form 
\begin{equation}
0=z^{2}w^{2}\left[ \left( 1-\frac{\nu }{w^{2}}\right) -\alpha ^{2}\left|
y\right| ^{2}\left( \left( 2+\frac{\nu }{w^{2}}\right) zw-3\right) ^{2}%
\right]  \label{cc3.5e}
\end{equation}
One can see that equation (\ref{cc3.5e}) depends essentially only on
combinations $\tau =zw$ and $s=\nu /w^{2}$ of variables $z$, $w$. Instead of
variables $z,w$ let us introduce variables $r,\tau $, connected with $z,w$
by means of relations. 
\[
\tau =zw,\qquad r=z\sqrt{w^{2}-\nu }=\tau \sqrt{1-s},\qquad s=\frac{\nu }{%
w^{2}} 
\]
Then the equation (\ref{cc3.5e}) transforms to the form 
\begin{equation}
z^{2}\left[ \frac{r^{2}}{\tau ^{2}}-\alpha ^{2}\left| y\right| ^{2}\left(
\left( 3-\frac{r^{2}}{\tau ^{2}}\right) \tau -3\right) ^{2}\right] =0
\label{cc3.6}
\end{equation}
Trivial solution $z=0$ of equation (\ref{cc3.6}) describes vector ${\bf 0x}=%
{\bf 0}.$ It means that the coordinate origin belongs to the tube. The tube $%
{\cal T}_{0y}$ is a result of all possible rotation of the curve $r=r\left(
\tau \right) $ around the vector ${\bf 0y}$. This curve is described by its
coordinate $\tau \left| y\right| $ along the vector ${\bf 0y}$ and
coordinate $r\left| y\right| $ in the orthogonal direction. Solution of
equation (\ref{cc3.6}) has the form 
\begin{equation}
r=\pm \frac{1}{2\alpha \left| y\right| }\left( -1\pm \sqrt{\left( 1+12\alpha
^{2}\left| y\right| ^{2}\tau \left( \tau -1\right) \right) }\right)
\label{cc3.8}
\end{equation}
Any section $\tau =$const of the three-dimensional surface ${\cal T}_{0y}$
form two (or zero) spheres, whose radii $r=r\left( \tau \right) $ are
determined by the relation (\ref{cc3.8}). Equation (\ref{cc3.8}) gives four
values of $r$ for any value of $\tau $, but only two of them are essential,
because radii $\ r$ and $-r$ describe the same surface. The two surfaces,
determined by (\ref{cc3.8}) cross between themselves at $\tau =0$ and $\tau
=1$, when $r=0$.

It follows from (\ref{cc3.8}) that 
\[
\lim_{\tau \rightarrow \infty }\frac{r}{\tau }=\pm \sqrt{3}, 
\]
It means the tube ${\cal T}_{0y}$ is infinite only in spacelike directions.
In the timelike directions the tube size is restricted.

In the vicinity of the vector ${\bf 0y}$, generating the tube, the shape of
the tube depends on interrelation between the intensity of the antisymmetry,
described by the constant $\alpha $, and the length of the vector ${\bf 0y.}$
$\alpha $ appears in the equation (\ref{cc3.8}) only in the combination $%
g=\alpha \left| y\right| $. In any case, when $\alpha \neq 0$, the tube $%
{\cal T}_{0y}$ does not degenerate into a one-dimensional curve. The tube $%
{\cal T}_{0y}$ is symmetric with respect to the reflection $\tau \rightarrow
1-\tau $. (See figures 1,2). If the antisymmetric structure is strong
enough, and $\alpha \left| y\right| >1/\sqrt{3}$, the tube ${\cal T}_{0y}$
is empty in its center in the sense that intersection of ${\cal T}_{0y}$
with the plane $\tau =0.5$ is empty. If $\alpha \left| y\right| <1/\sqrt{3}$
intersection of ${\cal T}_{0y}$ with the plane $\tau =0.5$ forms two
concentric spheres of radii 
\begin{equation}
r_{1}=\frac{3\alpha \left| y\right| }{2\left( \sqrt{\left( 1-3\alpha
^{2}\left| y\right| ^{2}\right) }+1\right) },\qquad r_{2}=\frac{3\alpha
\left| y\right| }{2\left( 1-\sqrt{\left( 1-3\alpha ^{2}\left| y\right|
^{2}\right) }\right) }  \label{cc3..7}
\end{equation}
If $\alpha \left| y\right| \ll 1$, one of radii is small $r_{1}=0.75\alpha
\left| y\right| $ and another one is large $r_{2}=1/\left( \alpha \left|
y\right| \right) $.

Let us consider now the second case (\ref{cc3.1b}), when the antisymmetric
structure is essential only at small distances. In this case one obtains
instead of equation (\ref{cc3.5e}) 
\begin{equation}
\sqrt{\left( w^{2}-\nu \right) }z\left| y\right| =\varepsilon \alpha \left( 
\frac{1}{1+\beta \left| y\right| ^{2}}-\frac{zw}{1+\beta \nu z^{2}\left|
y\right| ^{2}}+\frac{wz-1}{1+\beta \left( \nu z^{2}-2\varepsilon zw+1\right)
\left| y\right| ^{2}}\right)  \label{cc3.9}
\end{equation}
where the same designations (\ref{cc3.5b}) - (\ref{cc3.5c}) are used. At
large $z$ the equation (\ref{cc3.9}) transforms to the equation 
\begin{equation}
\sqrt{\left( w^{2}-\nu \right) }=\frac{\varepsilon \alpha }{z\left| y\right|
\left( 1+\beta \left| y\right| ^{2}\right) }  \label{cc3.10}
\end{equation}
Let $\nu =1$, and vector ${\bf 0x}$ is timelike Then at $z\rightarrow \infty 
$%
\begin{eqnarray}
r &=&z\sqrt{\left( w^{2}-1\right) }=\frac{\varepsilon \alpha }{\left|
y\right| \left( 1+\beta \left| y\right| ^{2}\right) }  \label{cc3.11} \\
\tau &=&zw=\sqrt{z^{2}+\frac{\alpha ^{2}}{\left| y\right| ^{2}\left( 1+\beta
\left| y\right| ^{2}\right) ^{2}}}\rightarrow \infty  \label{cc3.12}
\end{eqnarray}
It means that the tube is unlimited in the timelike direction ${\bf 0y}$ and
has a finite radius. If $\nu =-1$ (vector ${\bf 0x}$ is spacelike), the
quantity $z$ is restricted $z^{2}<\frac{\alpha ^{2}}{\left| y\right|
^{2}\left( 1+\beta \left| y\right| ^{2}\right) ^{2}}$, as it follows from (%
\ref{cc3.10}). The tube is limited in spacelike directions.

Thus, the local antisymmetric structure produces only local perturbation of
the tube shape. At the timelike infinity this perturbation reduces to a
nonvanishing radius of the tube. Any asymmetry of the world function
generates nondegeneracy of T-geometry. On the other hand, nondegeneracy of
T-geometry is connected with the particle mass geometrization and with
quantum effects \cite{R91}.

\section{Concluding remarks.}

The main goal of the nonsymmetric T-geometry development is its possible
application as a space-time geometry, especially as a space-time geometry of
microcosm. Approach and methods of T-geometry distinguish from those of the
Riemannian (pseudo-Riemannian) geometry, which is used now as a space-time
geometry. The Riemannian geometry imposes on the space-time geometry a
series of constraints. These restrictions are generated by methods used at
the description of the Riemannian geometry. Let us list some of them.

1. The continuity of space-time. This is very fine property which cannot be
tested by a direct experiment. T-geometry is insensitive to continuity, and
it is free of this constraint.

2. The Riemannian geometry is a geometry with fixed dimension. It is very
difficult to imagine a geometry with variable dimension in the scope of the
Riemannian geometry. Such a problem is absent in T-geometry.

3. For the Riemannian geometry description one uses a coordinate system.
Considering all possible coordinate systems, one can truncate this
restriction essentially. But the Riemannian geometry does not admit a
coordinate-free description. It is connected with the circumstance, that in
the scope of the Riemannian CG one does not construct axiomatics for any
special Riemannian geometry. Axiomatics of the standard (Euclidean or
pseudo-Euclidean ) geometry is used. Relations (and axiomatics) of any
special geometry are obtained as a result of the standard geometry
deformation. The coordinate system is a means of such a deformation
description, and it cannot be omitted in the Riemannian geometry. T-geometry
uses another method of description of the standard geometry deformation, and
T-geometry admits a coordinate-free description.

4. The Riemannian geometry uses a concept of a curve, which is essentially a
method of the Riemannian geometry description. The curve is considered
conventionally to be a geometrical object (but not as a method of the
geometry description), and nobody is interested in separation of properties
of geometry from properties imported by a use of the description in terms of
curves. In particular, in the Riemannian geometry the absolute parallelism
is absent, in general. Parallelism of two vectors at remote points is
established by means of a reference to a curve, along which the vector was
transported parallelly. In other words, geometrical property of parallelism
of two vectors is formulated in terms of the method of description, and
nobody knows how to remove this dependence on the methods of description.
T-geometry is free of this defect. The concept of a curve is not used at the
T-geometry construction. There is an absolute parallelism in T-geometry.

5. T-geometry uses a special geometrical language, which contains only
concepts immanent to the geometry in itself ($\Sigma $-function and finite
subspaces). One does not need to eliminate means of the geometry description.

6. This geometrical language admits one to consider and to investigate
effectively such a situation, when the future and the past are not
geometrically equivalent.

7. The means of the Riemannian geometry description suppress such an
important property of geometry as nondegeneracy. As a corollary the particle
mass geometrization appears to be impossible in the scope of Riemannian
geometry. Geometrization of the particle mass is important, when the mass of
a particle is unknown and must be determined from some geometrical
relations. It may appear to be important for determination of the mass
spectrum of elementary particles. T-geometry admits geometrization of the
particle mass.

8. Consideration of nondegeneracy and geometrization of the particle mass
have admitted one to make the important step in understanding of the
microcosm space-time geometry. One succeeded in explanation of
non-relativistic quantum effects as geometrical effects, generated by
nondegeneracy of the space-time geometry. There is a hope that asymmetry of
the space-time geometry will admit one to explain important characteristics
of elementary particles.


\begin{thebibliography}{99}
\bibitem{R00}  Yu.~A.~Rylov, Metric space: classification of finite
subspaces instead of constraints on metric. {\it Proceedings on analysis and
geometry}, Novosibirsk, Publishing House of Mathematical institute, 2000.
pp. 481-504, (in Russian). English version: e-print math.MG/9905111

\bibitem{R01}  Yu.~A.~Rylov, Description of metric space as a classification
of its finite subspaces. {\it Funndamental'aya i Prikladnaya Matematika,}
(in print). (in Russian).

\bibitem{R02}  Yu. A.~Rylov, Geometry without topology as a new conception
of geometry. {\it Int. J. Math. Math. Sci. }(in print), eprint {\it %
math.MG/0103002.}

\bibitem{S60}  J. L.~Synge, {\it Relativity: The General Theory},
North-Holland, Amsterdam, 1960.

\bibitem{T59}  V. A.~Toponogov, Riemannian spaces with curvature restricted
below, {\it Uspechi\ Matematicheskich\ Nauk}{.} {\bf 14}, 87, (1959). (in
Russian).

\bibitem{ABN86}  A. D.~Alexandrov, V. N.~Berestovski, I.~G.~Nikolayev,
Generalized Riemannian spaces, {\it Uspechi Matematicheskich Nauk,} {\bf 41}%
, iss.3, 1, (1986). (in Russian).

\bibitem{BGP92}  Yu.~Burago, M.~Gromov, G.~Perelman, Alexandrov spaces with
curvatures restricted below, {\it Uspechi Matematicheskich Nauk,} {\bf 47},
iss. 2, 3, (1992). (in Russian).

\bibitem{M28}  K.~Menger, Untersuchen \"{u}ber allgemeine Metrik, {\it %
Mathematische Annalen,} {\bf 100}, 75-113, (1928).

\bibitem{R90}  Yu. A.~Rylov, Extremal properties of Synge's world function
and discrete geometry. {\it J. Math. Phys.} {\bf 31}, 2876-2890, (1990).

\bibitem{R62}  Yu.~A.~Rylov, Description of Riemannian space by means of
finite interval. {\it Izvestiya Vysshich Uchebnych Zavedeniyi,} ser.
fis.mat. 132, (1962). (in Russian).

\bibitem{R64}  Yu. A.~Rylov, Relative gravitational field and conservation
laws in general relativity. {\it Ann. Phys.} (Leipzig) {\bf 12}, 329, (1964).

\bibitem{R91}  Yu. A.~Rylov, Non-Riemannian model of space-time responsible
for quantum effects. {\it J. Math. Phys.} {\bf 32}, 2092-2098, (1991).
\end{thebibliography}
\end{document}